\documentclass{article}
\usepackage{standalone}
\usepackage{tikz,caption}
\usepackage{amsfonts,amsmath}
\usepackage{graphicx}

\usetikzlibrary{
	knots,
	hobby,
	decorations.pathreplacing,
	shapes.geometric,
	calc
}

\usepackage{amsthm}

\newtheorem{theorem}{Theorem}
\newtheorem{proposition}{Proposition}
\newtheorem*{corollary}{Corollary}
\newtheorem*{remark}{Remark}

\theoremstyle{definition}
\newtheorem{definition}{Definition}

\usepackage[
backend=biber,
style=alphabetic,
]{biblatex}
\begin{document}

\title{Triangular decomposition of $SL_3$ skein algebras}

\author{Vijay Higgins}

\maketitle

\begin{abstract}
We give an $SL_3$ analogue of the triangular decomposition of the Kauffman bracket stated skein algebras described by Le. To any punctured bordered surface, we associate an $SL_3$ stated skein algebra which contains the $SL_3$ skein algebra of closed webs. These algebras admit natural algebra morphisms associated to the splitting of surfaces along ideal arcs. We give an explicit basis for the $SL_3$ stated skein algebra and show that the splitting morphisms are injective and describe their images. By splitting a surface along the edges of an ideal triangulation, we see that the $SL_3$ stated skein algebra of any ideal triangulable surface embeds into a tensor product of stated skein algebras of triangles.
\end{abstract}

\pagebreak

\section{Introduction}

The representation theory of a quantum group admits a diagrammatic calculus. Ribbon diagrams in the plane depict maps between representations, and skein relations among the diagrams correspond to relations among the maps. For the case of $SL_n$, it has been shown that such skein theories essentially completely describe the representation theory of the quantum group $U_q\mathfrak{sl}_n$ \cite{Kup96,CKM14}. Instead of restricting our attention to diagrams in the plane, we can associate to any surface a skein module which consists of formal linear combinations of diagrams on the surface modulo local skein relations. Such a skein module admits a natural algebra structure given by superimposing one diagram on top of another. For the $SL_n$ case, Sikora \cite{Sik05} has described a skein theory built from directed $n\text{-valent}$ ribbon graphs which is particularly amenable to the algebra structure since only one color of strand is used. For the $SL_2$ case, this skein theory coincides with that given by the Kauffman bracket skein relations and for the $SL_3$ case it coincides with that given by the Kuperberg bracket skein relations.

Since the skein algebras are defined as quotients of free modules which happen to admit natural algebra structures, it is difficult to study the algebra structures explicitly. In particular, it can be difficult to construct algebra morphisms whose domains are the skein algebra. Nevertheless, the $SL_2$ case of the Kauffman bracket skein algebra of surfaces is relatively well-studied since the diagrams are built out of curves on surfaces and various geometric and combinatorial techniques have been developed to handle them.

In \cite{BW11}, Bonahon and Wong defined an algebra embedding of the Kauffman bracket skein algebra into quantum Teichmuller space. This map is called the $\textit{quantum trace map}$ and is viewed as a quantization of the classical trace map. The quantum Teichmuller space is a certain quantum torus, an algebra whose presentation is constructed from an ideal triangulation of a surface. Inspired by the computations involved with checking that the quantum trace map is well-defined, Le developed in \cite{Le18} a triangular decomposition of the Kauffman bracket skein algebra by introducing the Kauffman bracket \textit{stated} skein algebra. The stated skein algebra is a finer version of the regular skein algebra with extra relations along the boundary allowing one to define splitting morphisms, which are algebra maps associated to the splitting of a surface along an ideal arc. The extra boundary relations are consistent with the regular skein algebra relations, so no information is lost when passing to the stated skein algebra. The splitting morphisms are injective, so no information is lost by passing to the triangular decomposition. One can define an algebra map out of a skein algebra by defining maps on stated skein algebras of triangles and then precomposing with the triangular decomposition. Using this method, Le was able to reconstruct the quantum trace map.

Le and Costantino further developed the theory of the Kauffman bracket stated skein algebras of surfaces in \cite{CL19}. There, they highlighted the connections between the stated skein algebra, the quantized ring of coordinate functions $\mathcal{O}_q(SL_2)$, and the Reshitikhin-Turaev invariant of ribbon tangles. This perspective provides hints for how to construct a stated skein algebra for Lie groups other than $SL_2.$

A triangular decomposition for the algebra of functions on $G\text{-character}$ varieties of surfaces has been developed by Korinman in \cite{Kor19} for quite a general class of Lie groups $G.$ Furthermore, Korinman-Quesney and Costantino-Le showed in \cite{KQ19, CL19} that for the $SL_2$ case, the triangular decomposition of the stated skein algebra fits into an exact sequence which parallels the exact sequence associated to the triangular decomposition of the character variety. One can hope that skein algebras associated to other Lie groups $G$ have analogous triangular decompositions lining up with those of the character variety.

The goal of this paper is to define a stated version of the $SL_3$ skein algebra that admits a triangular decomposition analogous to the Kauffman bracket stated skein algebra. A successful definition should satisfy the items in the following wish list:

\begin{itemize}
	\item The stated skein algebra of the monogon should be 1-dimensional.
	\item The stated skein algebra of the bigon should be $\mathcal{O}_q(SL_3).$ It should have a Hopf algebra structure such that any surface with a boundary arc has a comodule-algebra structure over $\mathcal{O}_q(SL_3)$ associated to splitting off a bigon from the boundary arc.
	\item The stated skein algebra of the triangle should be a braided tensor product $\mathcal{O}_q(SL_3) \underset{-}{\otimes} \mathcal{O}_q(SL_3).$
	\item The regular skein algebra should embed into the stated skein algebra.
	\item There should exist injective splitting maps which embed the skein algebra of an ideal triangulable surface into a tensor product of stated skein algebras of triangles, and this triangular decomposition should fit into an exact sequence of the form described in \cite{Kor19, KQ19, CL19}. 
\end{itemize}

In this paper we give a definition of a stated skein algebra for $SL_3$ and show that it satisfies the above items.

To show the injectivity of the splitting map associated to an arbitrary ideal arc, we show that it suffices to prove the injectivity of the splitting map for an ideal arc that bounds a monogon. While this case is much simpler to consider, it is still difficult to prove. We use the list of reduction rules suggested by the $SL_3$ skein algebra relations as in  \cite{SW07} and expand the list to include relations along the boundary. Using these reduction rules, we apply the Diamond Lemma to find an explicit basis that helps us prove the injectivity of the splitting map for an ideal arc bounding a monogon. To generalize these results to other $G$, we desire similar reduction rules for their skein theories or else we hope to find a replacement for the role that the basis plays in this paper.

\subsection{Acknowledgments} The author would like to thank his advisor Stephen Bigelow and Ken Goodearl for helpful conversations. The author would also like to thank Stephen Bigelow, Wade Bloomquist, Thang Le, and Adam Sikora for looking at an earlier version of this paper.

\section{The $SL_3$ stated skein algebra}

\begin{definition}
	A $\textit{punctured bordered surface}$ is a pair $(\Sigma', \mathcal{P}),$ where $\Sigma'$ is a smooth compact oriented surface, possibly with boundary, and $\mathcal{P}$ is a collection of finitely many points of $\Sigma'.$ We require that each boundary component of $\Sigma'$ contains at least one point of $\mathcal{P}.$ We do not require $\Sigma'$ to be connected. We let $\Sigma=\Sigma'\setminus \mathcal{P}.$ To simplify notation, we also refer to the pair $(\Sigma', \mathcal{P})$ simply by $\Sigma.$ A $\textit{boundary arc}$ of $\Sigma$ is a connected component of $\partial \Sigma.$
\end{definition}

For a punctured bordered surface $\Sigma$, a $\textit{web}$ in $\Sigma \times (0,1)$ is an embedding of a directed ribbon graph $\Gamma$ such that each interior vertex of $\Gamma$ in $\mathring{\Sigma} \times (0,1)$ is a trivalent sink or a trivalent source. We allow $\Gamma$ to have univalent vertices, called $\textit{endpoints}$, contained in $\partial \Sigma \times (0,1)$ such that for each boundary arc $b$ of $\Sigma$ the vertices contained in $b \times (0,1)$ have distinct heights. We require the web to have a vertical framing with respect to the $(0,1)$ component and we require that strands that terminate in a univalent vertex are transverse to $\partial \Sigma.$

We consider isotopies of webs in the class of webs. In particular, our isotopies must preserve the height order of boundary points of webs for each boundary arc of $\Sigma.$

For a web $\Gamma$ a $\textit{state}$ is a function $s: \partial \Gamma \rightarrow \{-,0,+\}$. A \textit{stated web} is a web together with a state. We will make use of the order $-<0<+$ on the set $\{-,0,+\}$. For notational purposes, it will be convenient to sometimes add states together. By identifying the state $-$ with the integer $-1$ and the state $+$ with the integer $1$, we partially define an addition on the set $\{-,0,+\}$ whenever the answer is contained in the set as well. 

\begin{definition}
A web $\Gamma$ in $\Sigma \times (0,1)$ is in generic position if the projection $\pi: \Sigma \times(0,1) \rightarrow \Sigma$ restricts to an embedding of $\Gamma$ except for the possibility of transverse double points in the interior of $\Sigma.$ Each web is isotopic to a web in generic position. A \textit{stated diagram} $D$ of a generic stated web $\Gamma$ is the projection $\pi(\Gamma)$ along with the over/undercrossing information at each double point and the height order and states of the boundary points of $\Gamma.$ Web diagrams are isotopic if they are isotopic through an isotopy of the surface.
\end{definition}

As in \cite{Le18} it will be convenient for us to record the local height order of the boundary points of a web diagram by drawing an arrow along a portion of the boundary arc of $\Sigma$.

Let $\mathcal{R}$ be a commutative ring with identity and with an invertible element $q^{1/3}.$ The quantum integer $[n]$ denotes the Laurent polynomial defined by $[n]=\frac{q^{n}-q^{-n}}{q-q^{-1}}.$

\begin{definition}
The \textit{$SL_3$ stated skein algebra} $\mathcal{S}_q^{SL_3}(\Sigma)$ is the $\mathcal{R}\text{-module}$ freely spanned by isotopy classes of webs in $\Sigma \times (0,1)$ modulo the following relations.
\end{definition}
\textbf{Interior relations:}

\begin{equation}\tag{$\text{I}1\text{a}$} \begin{tikzpicture}[use Hobby shortcut, baseline=3ex]
\begin{knot}[
consider self intersections=true,
%  draft mode=crossings,
flip crossing=1,
ignore endpoint intersections=false]
\strand [->] (0,1)..(1/2,0);
\strand [->] (1/2,1)..(0,0);
\end{knot}
\end{tikzpicture}=q^{2/3}  \begin{tikzpicture}[baseline=3ex]
\draw [->] (0,1)--(0,0);
\draw [->] (1/2,1)--(1/2,0);
\end{tikzpicture} +q^{-3-1/3} \begin{tikzpicture}[baseline=3ex]
\draw [<-] (0,0)--(1/4,1/4);
\draw [<-] (1/2,0)--(1/4,1/4);
\draw (1/4,1/4)--(1/4,3/4);
\draw [-<] (1/4,3/4)--(1/4,1/2);
\draw (1/4,3/4)--(0,1);
\draw (1/4,3/4)--(1/2,1);
\draw [->] (0,1)--(1/8,7/8);
\draw [->] (1/2,1)--(3/8,7/8);
\end{tikzpicture}\end{equation}

\begin{equation}\tag{$\text{I}1\text{b}$} \begin{tikzpicture}[use Hobby shortcut, baseline=3ex]
\begin{knot}[
consider self intersections=true,
%  draft mode=crossings,
ignore endpoint intersections=false]
\strand [->] (0,1)..(1/2,0);
\strand [->] (1/2,1)..(0,0);
\end{knot}
\end{tikzpicture}=q^{-2/3}  \begin{tikzpicture}[baseline=3ex]
\draw [->] (0,1)--(0,0);
\draw [->] (1/2,1)--(1/2,0);
\end{tikzpicture} +q^{-3+1/3} \begin{tikzpicture}[baseline=3ex]
\draw [<-] (0,0)--(1/4,1/4);
\draw [<-] (1/2,0)--(1/4,1/4);
\draw (1/4,1/4)--(1/4,3/4);
\draw [-<] (1/4,3/4)--(1/4,1/2);
\draw (1/4,3/4)--(0,1);
\draw (1/4,3/4)--(1/2,1);
\draw [->] (0,1)--(1/8,7/8);
\draw [->] (1/2,1)--(3/8,7/8);
\end{tikzpicture}\end{equation}

\begin{equation}\tag{$\text{I}2$}
\begin{tikzpicture}[baseline=3ex]
\draw [->] (0,0)--(0,1/8);
\draw [-<] (0,1/8)--(0,1/2);
\draw [->] (0,0)--(0,7/8);
\draw (0,7/8)--(0,1);
\draw [-<] (1/2,0)--(1/2,1/8);
\draw [->] (1/2,0)--(1/2,1/2);
\draw [-<] (1/2,1/2)--(1/2,7/8);
\draw (1/2,0)--(1/2,1);
\draw [-<] (0,1/4)--(1/4,1/4);
\draw (0,1/4)--(1/2,1/4);
\draw [->] (0,3/4)--(1/4,3/4);
\draw (1/4,3/4)--(1/2,3/4);
\end{tikzpicture}=q^6( \begin{tikzpicture}[baseline=3ex]
\draw [->] (0,0)--(0,1/2);
\draw (0,0)--(0,1);
\draw [-<] (1/2,0)--(1/2,1/2);
\draw (1/2,0)--(1/2,1);
\end{tikzpicture} + \begin{tikzpicture}[baseline=3ex]
\draw [<-] plot [smooth] coordinates {(0,1)(1/4,3/4)(1/2,1)};
\draw [->] plot [smooth] coordinates {(0,0)(1/4,1/4)(1/2,0)};
\end{tikzpicture})
\end{equation}

\begin{equation}\tag{$\text{I}3$}
\begin{tikzpicture}[baseline=3ex]
\draw [-<] (1/4,0)--(1/4,1/8);
\draw (1/4,0)--(1/4,1/4);
\draw [-<] (1/4,3/4)--(1/4,7/8);
\draw (1/4,3/4)--(1/4,1);
\draw [->](1/4,1/4) to [out=135, in=-90] (0,1/2);
\draw (0,1/2) to [out=90, in =-135] (1/4,3/4);
\draw [->] (1/4,1/4) to [out=45, in=-90] (1/2,1/2);
\draw (1/2,1/2) to [out=90, in=-45] (1/4,3/4);
\end{tikzpicture}=-q^3[2] \hspace{.1cm} \begin{tikzpicture}[baseline=3ex]
\draw [->] (0,1)--(0,1/2);
\draw (0,1)--(0,0);
\end{tikzpicture}
\end{equation}

\begin{equation}\tag{$\text{I}4\text{a}$}
\begin{tikzpicture}[baseline=-4ex]
\draw [>-] (0,0) arc [radius=.5, start angle=90, end angle=465];
\end{tikzpicture}=[3]
\end{equation}

\begin{equation}\tag{$\text{I}4\text{b}$}
\begin{tikzpicture}[baseline=-4ex]
\draw [<-] (0,0) arc [radius=.5, start angle=90, end angle=465]; \end{tikzpicture}=[3]
\end{equation}

\textbf{Boundary relations:}

\begin{equation}\tag{$\text{B}1$}
\begin{tikzpicture}[baseline=.5ex]
\draw [line width=1.5, ->] (-1/4,0)--(3/4,0);
\draw [->] (1/4,0)--(1/4,1/2);
\draw (1/4,0)--(1/4,1);
\node [below] at (1/4,0) {$a+b$};
\end{tikzpicture}=(-1)^{a+b}q^{-1/3-(a+b)}\begin{tikzpicture}[baseline=.5ex]
\draw [line width=1.5, ->] (-1/4,0)--(3/4,0);
\draw (0,0)--(1/4,1/2);
\draw (1/2,0)--(1/4,1/2);
\draw (1/4,1/2)--(1/4,1);
\draw [-<] (0,0)--(1/8,1/4);
\draw [-<] (1/2,0)--(3/8,1/4);
\draw [-<] (1/4,1)--(1/4,3/4);
\node [below] at (0,-.07) {$a$};
\node [below] at (1/2,0) {$b$};
\end{tikzpicture} \text{(for $b>a$)}
\end{equation}

\begin{equation}\tag{$\text{B}2$}
\begin{tikzpicture}[baseline=.5ex]
\draw [line width =1.5, ->] (-1/4,0)--(3/4,0);
\draw (0,1)--(0,0);
\draw (1/2,1)--(1/2,0);
\draw [->] (0,1)--(0,1/2);
\draw [->] (1/2,1)--(1/2,1/2);
\node [below] at (0,0) {$b$};
\node [below] at (1/2,-.07) {$a$};
\end{tikzpicture}=q^{-1} \begin{tikzpicture}[baseline=.5ex]
\draw [line width =1.5, ->] (-1/4,0)--(3/4,0);
\draw (0,1)--(0,0);
\draw (1/2,1)--(1/2,0);
\draw [->] (0,1)--(0,1/2);
\draw [->] (1/2,1)--(1/2,1/2);
\node [below] at (0,-.07) {$a$};
\node [below] at (1/2,0) {$b$};
\end{tikzpicture}+q^{-3}\begin{tikzpicture}[baseline=.5ex]
\draw (0,0)--(1/4,1/4);
\draw [->] (1/4,1/4)--(3/8,1/8);
\draw (1/2,0)--(1/4,1/4);
\draw [->] (1/4,1/4)--(1/8,1/8);
\draw (1/4,1/4)--(1/4,3/4);
\draw [-<] (1/4,3/4)--(1/4,1/2);
\draw (1/4,3/4)--(0,1);
\draw (1/4,3/4)--(1/2,1);
\draw [->] (0,1)--(1/8,7/8);
\draw [->] (1/2,1)--(3/8,7/8);
\draw [line width=1.5, ->] (-1/4,0)--(3/4,0);
\node [below] at (0,-.07) {$a$};
\node [below] at (1/2,0) {$b$};
\end{tikzpicture} \text{(for $b>a$)}
\end{equation}

\begin{equation}\tag{$\text{B}3$}
\begin{tikzpicture}[baseline=.5ex]
\draw [line width=1.5, ->] (-1/4,0)--(3/4,0);
\draw (0,0)--(1/4,1/2);
\draw (1/2,0)--(1/4,1/2);
\draw (1/4,1/2)--(1/4,1);
\draw [-<] (0,0)--(1/8,1/4);
\draw [-<] (1/2,0)--(3/8,1/4);
\draw [-<] (1/4,1)--(1/4,3/4);
\node [below] at (0,0) {$a$};
\node [below] at (1/2,0) {$a$};
\end{tikzpicture}=0 \hspace{1cm} \text{(for any $a \in \{-,0,+\}$)}
\end{equation}

\begin{equation}\tag{$\text{B}4$}
\begin{tikzpicture}[baseline=.5ex]
\draw (0,1/2)--(0,0);
\draw [->] (0,1/2)--(0,1/4);
\draw (1/4,0)--(0,1/2);
\draw [->] (0,1/2)--(1/8,1/4);
\draw (-1/4,0)--(0,1/2);
\draw [->] (0,1/2)--(-1/8,1/4);
\draw [line width=1.5,->] (-1/2,0)--(1/2,0);
\node [below] at (-1/4,0) {$-$};
\node [below] at (0,0) {$0$};
\node [below] at (1/4,0) {$+$};
\end{tikzpicture}=q^{-2}\begin{tikzpicture}[baseline=.5ex]
\draw [line width=1.5,->] (-1/2,0)--(1/2,0);
\end{tikzpicture}
\end{equation}

The interior relations above hold for local diagrams contained in an embedded disk in $\Sigma.$ The boundary relations hold for local diagrams in a neighborhood of a point of $\partial \Sigma.$ The thicker line denotes a portion of a boundary arc while the thin lines belong to a web. The arrow along the boundary arc indicates the height order of that boundary arc. For example, in the diagram on the right side of relation (B1), the endpoint with the state $b$ has a greater height than the endpoint with the state $a$.

The module defined above admits a natural multiplication where the product $\Gamma_1\Gamma_2$ of two stated webs $\Gamma_1,\Gamma_2$ in $\Sigma \times(0,1)$ is given by isotoping $\Gamma_1$ so that it is contained in $\Sigma \times (1/2, 1),$ isotoping $\Gamma_2$ so that it is contained in $\Sigma \times (0,1/2)$, and then taking the union of these two stated webs in $\Sigma \times (0,1).$ This gives $\mathcal{S}_q^{SL_3}(\Sigma)$ an associative, unital $\mathcal{R}\text{-algebra}$ structure.

\section{Consequences of the defining relations}

\begin{proposition}
The following relations are consequences of the defining relations.
\end{proposition}

\begin{equation}\tag{a}
q^{-8/3}\begin{tikzpicture}[use Hobby shortcut, baseline=3ex]
\begin{knot}[
consider self intersections=true,
%  draft mode=crossings,
ignore endpoint intersections=false]
\strand [<-] plot [smooth] coordinates{(0,0)(0,1/4)(1/4,3/4)(1/2,1/2)(1/4,1/4)(0,3/4)(0,1)};
\end{knot}
\end{tikzpicture}= \begin{tikzpicture}[baseline=3ex]
\draw [->] (0,1)--(0,0);
\end{tikzpicture}=q^{8/3}\begin{tikzpicture}[use Hobby shortcut, baseline=3ex]
{\begin{knot}[
	consider self intersections=true,
	%  draft mode=crossings,
	ignore endpoint intersections=false]
	\strand [->] plot [smooth] coordinates {(0,1)(0,3/4)(1/4,1/4)(1/2,1/2)(1/4,3/4)(0,1/4)(0,0)};
	\end{knot}}
\end{tikzpicture}
\end{equation}

\begin{equation}\tag{b}
-q^{-4}\begin{tikzpicture}[use Hobby shortcut, baseline=3ex]
\begin{knot}[
consider self intersections=true,
%  draft mode=crossings,
ignore endpoint intersections=false]
\strand plot [smooth] coordinates {(0,1)(-1/4,3/4)(1/4,1/4)(1/4,0)};
\strand plot [smooth] coordinates {(0,1)(0,3/4)(-1/4,1/2)(-1/4,1/4)(0,0)};
\strand plot [smooth] coordinates {(0,1)(1/4,3/4)(1/4,1/2)(-1/4,0)};
\end{knot}
\draw [->] plot [smooth] coordinates {(0,1)(-1/4,3/4)(1/4,1/4)(1/4,0)};
\draw (0,1)--(1/14,13/14);
\draw [<-] (-1/4,0)--(-1/4+1/16,1/16);
\draw [<-] (0,0)--(-1/16,1/16);
\end{tikzpicture}=\begin{tikzpicture}[baseline=3ex]
\draw [->] (0,1)--(-1/4,0);
\draw [->] (0,1)--(0,0);
\draw [->] (0,1)--(1/4,0);
\end{tikzpicture}=-q^4 \begin{tikzpicture}[use Hobby shortcut, baseline=3ex]
\begin{knot}[
consider self intersections=true,
%  draft mode=crossings,
ignore endpoint intersections=false]
\strand plot [smooth] coordinates {(0,1)(1/4,3/4)(1/4,1/2)(-1/4,0)};
\strand plot [smooth] coordinates {(0,1)(0,3/4)(-1/4,1/2)(-1/4,1/4)(0,0)};
\strand plot [smooth] coordinates {(0,1)(-1/4,3/4)(1/4,1/4)(1/4,0)};
\end{knot}
\draw [->] plot [smooth] coordinates {(0,1)(1/4,3/4)(1/4,1/2)(-1/4,0)};
\draw [<-] (0,0)--(-1/16,1/16);
\draw [<-] (1/4,0)--(1/4,1/100);
\draw (0,1)--(-1/14,18/19);
\end{tikzpicture}
\end{equation}

\begin{equation}\tag{c}
\begin{tikzpicture}[baseline=.5ex]
\draw [->] (-1/2,0) arc [radius=.5, start angle=180, end angle=90];
\draw (-1/2,0) arc [radius=.5, start angle=180, end angle=0];
\draw [line width=1.5, ->] (-3/4,0)--(3/4,0);
\node [below] at (-1/2,-.07) {$a$};
\node [below] at (1/2,0) {$b$};
\end{tikzpicture}=-q^{-4/3}\delta_{a+b,0} \begin{tikzpicture}[baseline=.5ex]
\draw [line width=1.5, ->] (-3/4,0)--(3/4,0);
\end{tikzpicture}
\end{equation}

\begin{equation}\tag{d}
\begin{tikzpicture}[baseline=.5ex]
\draw [->] (1/2,1) arc [radius=.5, start angle=0, end angle=-90];
\draw (1/2,1) arc [radius=.5, start angle=0, end angle=-180];
\draw [line width=1.5, <-] (-3/4,0)--(3/4,0);
\end{tikzpicture}=-q^{-4/3}\sum_{a+b=0}\begin{tikzpicture}[baseline=.5ex]
\draw [->] (-1/2,0)--(-1/2,1/2);
\draw (-1/2,0)--(-1/2,1);
\draw [-<] (1/2,0)--(1/2,1/2);
\draw (1/2,0)--(1/2,1);
\draw [line width=1.5, <-] (-3/4,0)--(3/4,0);
\node [below] at (-1/2,0) {$b$};
\node [below] at (1/2,-.07) {$a$};
\end{tikzpicture}
\end{equation}

\begin{equation}\tag{e}
\begin{tikzpicture}[baseline=.5ex]
\draw [->] (1/2,0) arc [radius=.5, start angle=0, end angle =90];
\draw (1/2,0) arc [radius=.5, start angle=0, end angle=180];
\draw [line width=1.5, ->] (-3/4,0)--(3/4,0);
\node [below] at (-1/2,-.07) {$a$};
\node [below] at (1/2,0){$b$};
\end{tikzpicture}=-q^{-4/3}q^{2a}\delta_{a+b,0} \begin{tikzpicture}[baseline=.5ex]
\draw [line width=1.5, ->] (-3/4,0)--(3/4,0);
\end{tikzpicture}
\end{equation}

\begin{equation}\tag{f}
\begin{tikzpicture}[baseline=.5ex]
\draw [-<] (1/2,1) arc [radius=.5, start angle=0, end angle=-90];
\draw (1/2,1) arc [radius=.5, start angle=0, end angle=-180];
\draw [line width=1.5, <-] (-3/4,0)--(3/4,0);
\end{tikzpicture}=-q^{-4/3}\sum_{a+b=0}q^{2a}\begin{tikzpicture}[baseline=.5ex]
\draw [-<] (-1/2,0)--(-1/2,1/2);
\draw (-1/2,0)--(-1/2,1);
\draw [->] (1/2,0)--(1/2,1/2);
\draw (1/2,0)--(1/2,1);
\draw [line width=1.5, <-] (-3/4,0)--(3/4,0);
\node [below] at (-1/2,0) {$b$};
\node [below] at (1/2,-.07) {$a$};
\end{tikzpicture}
\end{equation}

\begin{equation}\tag{g}
\begin{tikzpicture}[baseline=.5ex]
\draw [->] (1/2,0) arc [radius=.5, start angle=0, end angle =90];
\draw (1/2,0) arc [radius=.5, start angle=0, end angle=180];
\draw [line width=1.5, <-] (-3/4,0)--(3/4,0);
\node [below] at (-1/2,-.07) {$a$};
\node [below] at (1/2,0){$b$};
\end{tikzpicture}=-q^{4/3}\delta_{a+b,0} \begin{tikzpicture}[baseline=.5ex]
\draw [line width=1.5, <-] (-3/4,0)--(3/4,0);
\end{tikzpicture}
\end{equation}

\begin{equation}\tag{h}
\begin{tikzpicture}[baseline=.5ex]
\draw [-<] (1/2,1) arc [radius=.5, start angle=0, end angle=-90];
\draw (1/2,1) arc [radius=.5, start angle=0, end angle=-180];
\draw [line width=1.5, ->] (-3/4,0)--(3/4,0);
\end{tikzpicture}=-q^{4/3}\sum_{a+b=0}\begin{tikzpicture}[baseline=.5ex]
\draw [-<] (-1/2,0)--(-1/2,1/2);
\draw (-1/2,0)--(-1/2,1);
\draw [->] (1/2,0)--(1/2,1/2);
\draw (1/2,0)--(1/2,1);
\draw [line width=1.5, ->] (-3/4,0)--(3/4,0);
\node [below] at (-1/2,0) {$b$};
\node [below] at (1/2,-.07) {$a$};
\end{tikzpicture}
\end{equation}

\begin{equation}\tag{i}
\begin{tikzpicture}[baseline=.5ex]
\draw [-<] (1/2,0) arc [radius=.5, start angle=0, end angle =90];
\draw (1/2,0) arc [radius=.5, start angle=0, end angle=180];
\draw [line width=1.5, <-] (-3/4,0)--(3/4,0);
\node [below] at (-1/2,-.07) {$a$};
\node [below] at (1/2,0){$b$};
\end{tikzpicture}=-q^{4/3}q^{2b}\delta_{a+b,0} \begin{tikzpicture}[baseline=.5ex]
\draw [line width=1.5, <-] (-3/4,0)--(3/4,0);
\end{tikzpicture}
\end{equation}

\begin{equation}\tag{j}
\begin{tikzpicture}[baseline=.5ex]
\draw [->] (1/2,1) arc [radius=.5, start angle=0, end angle=-90];
\draw (1/2,1) arc [radius=.5, start angle=0, end angle=-180];
\draw [line width=1.5, ->] (-3/4,0)--(3/4,0);
\end{tikzpicture}=-q^{4/3}\sum_{a+b=0}q^{2b}\begin{tikzpicture}[baseline=.5ex]
\draw [->] (-1/2,0)--(-1/2,1/2);
\draw (-1/2,0)--(-1/2,1);
\draw [-<] (1/2,0)--(1/2,1/2);
\draw (1/2,0)--(1/2,1);
\draw [line width=1.5, ->] (-3/4,0)--(3/4,0);
\node [below] at (-1/2,0) {$b$};
\node [below] at (1/2,-.07) {$a$};
\end{tikzpicture}
\end{equation}

\begin{equation}\tag{k}
\begin{tikzpicture}[baseline=.5ex]
\draw (0,1/2)--(0,0);
\draw [->] (0,1/2)--(0,1/4);
\draw (0,1/2)--(-1/2,0);
\draw [->] (0,1/2)--(-1/4,1/4);
\draw (0,1/2)--(1/2,0);
\draw [->] (0,1/2)--(1/4,1/4);
\draw [line width=1.5, ->] (-3/4,0)--(3/4,0);
\node [below] at (-1/2,0) {$\sigma_1$};
\node [below] at (0,0) {$\sigma_2$};
\node [below] at (1/2,0) {$\sigma_3$};
\end{tikzpicture}=\begin{cases}
q^{-2}(-q)^{l(\sigma)} \text{ if $\sigma=(\sigma_1,\sigma_2,\sigma_3)\in S_3$}\\
0 \text{ if $(\sigma_1,\sigma_2,\sigma_3) \notin S_3$}\\
\end{cases} \text{(same for sinks)}
\end{equation}

\begin{equation}\tag{l}
\begin{tikzpicture}[baseline=.5ex]
\draw (0,1/2)--(0,1);
\draw [->] (0,1/2)--(0,3/4);
\draw (0,1/2)--(-1/2,1);
\draw [->] (0,1/2)--(-1/4,3/4);
\draw (0,1/2)--(1/2,1);
\draw [->] (0,1/2)--(1/4,3/4);
\draw [line width=1.5, <-] (-3/4,0)--(3/4,0);
\end{tikzpicture}=q^{-2} \sum_{\sigma \in S^3} (-q)^{l(\sigma)} \begin{tikzpicture}[baseline=.5ex]
\draw (-1/2,0)--(-1/2,1);
\draw [->] (-1/2,0)--(-1/2,1/2);
\draw (0,0)--(0,1);
\draw [->] (0,0)--(0,1/2);
\draw (1/2,0)--(1/2,1);
\draw [->] (1/2,0)--(1/2,1/2);
\draw [line width=1.5, <-] (-3/4,0)--(3/4,0);
\node [below] at (-1/2,0) {$\sigma_3$};
\node [below] at (0,0) {$\sigma_2$};
\node [below] at (1/2,0) {$\sigma_1$};
\end{tikzpicture} \text{(same for sinks)}
\end{equation}

\begin{equation}\tag{m}
\begin{tikzpicture}[baseline=.5ex]
\draw (0,1/2)--(0,0);
\draw [->] (0,1/2)--(0,1/4);
\draw (0,1/2)--(-1/2,0);
\draw [->] (0,1/2)--(-1/4,1/4);
\draw (0,1/2)--(1/2,0);
\draw [->] (0,1/2)--(1/4,1/4);
\draw [line width=1.5, <-] (-3/4,0)--(3/4,0);
\node [below] at (-1/2,0) {$\sigma_3$};
\node [below] at (0,0) {$\sigma_2$};
\node [below] at (1/2,0) {$\sigma_1$};
\end{tikzpicture}=\begin{cases}
-q^{2}(-q)^{l(\sigma)} \text{ if $\sigma=(\sigma_1,\sigma_2,\sigma_3)\in S_3$}\\
0 \text{ if $(\sigma_1,\sigma_2,\sigma_3) \notin S_3$}\\
\end{cases} \text{(same for sinks)}
\end{equation}

\begin{equation}\tag{n}
\begin{tikzpicture}[baseline=.5ex]
\draw (0,1/2)--(0,1);
\draw [->] (0,1/2)--(0,3/4);
\draw (0,1/2)--(-1/2,1);
\draw [->] (0,1/2)--(-1/4,3/4);
\draw (0,1/2)--(1/2,1);
\draw [->] (0,1/2)--(1/4,3/4);
\draw [line width=1.5, ->] (-3/4,0)--(3/4,0);
\end{tikzpicture}=-q^{2} \sum_{\sigma \in S^3} (-q)^{l(\sigma)} \begin{tikzpicture}[baseline=.5ex]
\draw (-1/2,0)--(-1/2,1);
\draw [->] (-1/2,0)--(-1/2,1/2);
\draw (0,0)--(0,1);
\draw [->] (0,0)--(0,1/2);
\draw (1/2,0)--(1/2,1);
\draw [->] (1/2,0)--(1/2,1/2);
\draw [line width=1.5, ->] (-3/4,0)--(3/4,0);
\node [below] at (-1/2,0) {$\sigma_1$};
\node [below] at (0,0) {$\sigma_2$};
\node [below] at (1/2,0) {$\sigma_3$};
\end{tikzpicture} \text{(same for sinks)}
\end{equation}

In the notation above, we consider the permuation $(-,0,+)$ to be the identity permutation and $l(\sigma)$ denotes the length of the permutation $\sigma.$

\begin{proof}
Relations (a) and (b) follow from the defining interior relations.

The relations involving boundary orientations pointing to the right can be checked by reducing both sides according to the algorithm given by the Diamond Lemma described in Theorem 2.

The relations involving boundary orientations pointing to the left can be derived from those involving orientations pointing to the right by sliding the boundary points horizontally to reverse the height order and using the twisting relations (a) and (b).
\end{proof}

\section{The splitting morphism}

As in \cite{Le18}, our stated skein algebras of punctured bordered surfaces satisfy a compatibility with the gluing and splitting of surfaces. If $\Sigma$ is a punctured bordered surface and $a$ and $b$ are two boundary arcs of $\Sigma,$ we can obtain a new punctured bordered surface $\bar{\Sigma}=\Sigma/(a=b)$ by gluing the arcs $a$ and $b$ together in the way compatible with the orientation of $\Sigma.$ It is the reverse of this process that gives us an algebra morphism from $\mathcal{S}_q^{
SL_3}(\bar{\Sigma})$ to $\mathcal{S}_q^{SL_3}(\Sigma)$ associated with splitting the surface $\bar{\Sigma}$ along an ideal arc $c.$

\begin{definition}
If $\Sigma$ is a punctured bordered surface, an \textit{ideal arc} in $\Sigma$ is a proper embedding $c: (0,1) \rightarrow \mathring{\Sigma}$ such that its endpoints are (not necessarily distinct) points in the set of punctures, $\mathcal{P}.$
\end{definition}

Let $p:\Sigma \rightarrow \Sigma/(a=b)=:\bar{\Sigma}$ be the projection map associated to the gluing. Then $c:=p(a)=p(b)$ is an ideal arc. We will define a splitting morphism $$\Delta_c: \mathcal{S}^{SL_3}_q (\bar{\Sigma}) \rightarrow \mathcal{S}^{SL_3}_q (\Sigma)$$ by defining it on stated webs in $\bar{\Sigma}\times (0,1)$ and then checking that it is well-defined on $\mathcal{S}^{SL_3}_q(\bar{\Sigma}).$

For a stated web $(\Gamma,s)$ in $\bar{\Sigma}\times(0,1)$ we first isotope it so that $\Gamma$ intersects $c \times(0,1)$ transversely in points of distinct heights. By defining $p$ to act trivially on the $(0,1)$ factor, we can extend it to a map $p:\Sigma \times (0,1) \rightarrow \bar{\Sigma} \times (0,1).$ We then consider $p^{-1}(\Gamma),$ which is a web in $\Sigma \times (0,1).$ Except for the points of $p^{-1}(c \cap \Gamma),$ each boundary point of $p^{-1}(\Gamma)$ inherits a state from $\Gamma.$

We will say that $s'$ is an $\textit{admissible state}$ for $p^{-1}(\Gamma)$ if $s'(p^{-1}(x))=s(x)$ for all $x \in \partial \Gamma$ and if $y,z \in p^{-1}(\Gamma \cap c)$ then $s'(y)=s'(z).$

We define the splitting morphism on a stated web $(\Gamma,s)$ in $\bar{\Sigma} \times (0,1)$ by $$\Delta_c(\Gamma,s)=\sum_{\text{admissible } s'} (p^{-1}(\Gamma),s').$$

\begin{theorem}
\begin{itemize}
\item[(a)] The map $\Delta_c$ described above extends linearly to a well-defined algebra morphism $\Delta_c:\mathcal{S}_q^{SL_3}(\bar{\Sigma}) \rightarrow \mathcal{S}_q^{SL_3}(\Sigma)$.

\item[(b)] If $a$ and $b$ are two ideal arcs with disjoint interiors, then we have $$\Delta_a \circ \Delta_b=\Delta_b \circ \Delta_a.$$
\end{itemize}

\end{theorem}

As in \cite{Le18}, the map $\Delta_c$ is injective, but we will postpone a discussion of this fact until Section 8.

\begin{proof}
If $\Delta_c$ is well-defined, then the fact that it is an algebra morphism and that it satisfies the property given in part (b) of the Theorem 1 follows from the definition of the splitting morphism.

To check that it is well-defined, we first check that the effect of passing cups, caps, vertices, and crossings past the ideal arc $c$ commutes with the application of $\Delta_c.$ This will tell us that the splitting morphism is well-defined with respect to isotopies of diagrams. Cups and caps can slide past the arc because of relations (c)-(j) from above. To slide a vertex past the arc, we can first rotate the vertex, using the fact that cups and caps can slide past the arc, until it appears as in relations (k)-(n). Since crossings can be rewritten as a linear combination of cups, caps and vertices, this allows us to pass a crossing past the arc.

If strands intersecting $c \times (0,1)$ are isotoped vertically so as to alter their height order, then on a diagram this has the effect of a Reidemeister 2 move.  Since crossings can slide past $c$, we can isotope the disk containing the Reidemeister 2 move on the diagram past $c$ and then perform the move. This tells us that the splitting map is well-defined on isotopy classes of webs.

To check that the splitting morphism respects the defining relations of $\mathcal{S}_q^{SL_3}(\Sigma)$ we observe that if $c$ cuts through a disk or half disk appearing in one of the defining relations, we can isotope the diagram away from $c$ first and then apply the relation.
\end{proof}

\section{A basis for the stated skein algebra}

%As explained in \cite{SW07}, the Diamond Lemma has been successful in producing bases for various modules of diagrams on surfaces. It produces a basis for Kuperberg's $A_1,A_2,B_2,$ and $G_2$ webs. The case of $A_1$ is the case of the Kauffman bracket and Le 

If a module is defined as a quotient of a free module by a list of relations, and if each relation can be interpreted as a reduction rule that permits the replacement of one element by a linear combination of simpler elements, then the module is a good candidate for an attempted application of the Diamond Lemma to produce a basis. As explained in \cite{SW07}, the Diamond Lemma can accommodate modules built out of diagrams on surfaces and it has been successful in producing bases for webs on surfaces for the cases of Kuperberg's webs of type $A_1,A_2,B_2,$ and $G_2.$ In \cite{Le18}, Le organized the new boundary relations into reduction rules that are compatible with the reduction rules coming from the Kauffman bracket skein algebra and then applied the Diamond Lemma to find a basis. In this section, we will do the same for the $SL_3$ case.

We first summarize our goal. To apply the Diamond Lemma, we need to realize our skein module as a quotient of a free module by reduction rules that are terminal and locally confluent. The defining relations from Section 2 provide a starting point for a list of reduction rules. We will introduce a measure of complexity that allows us to say that the diagrams in the right side of each defining relation are simpler than the diagram on the left side. Using a reduction rule on a diagram $D$ replaces that diagram with a linear combination of simpler diagrams. We call any linear combination of diagrams obtained by applying a sequence of reduction rules to $D$ a $\textit{descendant}$ of $D$, and we call the diagrams appearing in the linear combination $\textit{descendant diagrams}$ of $D$. If there exists no infinite chain of descendant diagrams for $D$, then $D$ can be written as a linear combination of irreducible diagrams by repeatedly applying reduction rules to the diagram and to its descendants. If no diagram admits an infinite chain of descendant diagrams, then the reduction rules are called \textit{terminal} and this property implies that irreducible diagrams span our module. Sometimes more than one reduction rule will apply to a diagram. If there is always a common descendant for any two ways of reducing a diagram, then the reduction rules are called $\textit{locally confluent}.$ If the set of reduction rules are terminal and locally confluent, then the set of irreducible diagrams forms a basis for our module, by Theorem 2.3 in \cite{SW07}.

In anticipation of issues regarding local confluence, we need to introduce the following redundant relations:

\begin{equation}\tag{S}
\begin{tikzpicture}[baseline=.5ex]
\draw [dashed] (0,0) circle (.25);
\draw  (0,0) circle (.5);
\draw [>-] (.5,0) arc [radius=.5, start angle=0, end angle=180];
\draw (0,0) circle (.75);
\draw [<-] (.75,0) arc [radius=.75, start angle=0, end angle =180];
\draw [dashed] (0,0) circle (1);
\end{tikzpicture}=
\begin{tikzpicture}[baseline=.5ex]
\draw [dashed] (0,0) circle (.25);
\draw  (0,0) circle (.5);
\draw [<-] (.5,0) arc [radius=.5, start angle=0, end angle=180];
\draw (0,0) circle (.75);
\draw [>-] (.75,0) arc [radius=.75, start angle=0, end angle =180];
\draw [dashed] (0,0) circle (1);
\end{tikzpicture}
\end{equation}

\begin{equation}\tag{$C_k$}
\begin{tikzpicture}[baseline=.5ex]
\draw [line width=1.5, ->] (0,0)--(3.5,0);
\draw (1/4,0)--(3/4,1/2);
\draw [-<] (1/4,0)--(1/2,1/4);
\draw (3/4,0)--(3/4,1/2);
\draw [-<] (3/4,0)--(3/4,1/4);
\draw (5/4,0)--(5/4,1/2);
\draw [-<] (5/4,0)--(5/4,1/4);
\draw (2,0)--(2,1/2);
\draw [-<] (2,0)--(2,1/4);
\draw (2+1/2,0)--(2+1/2,1/2);
\draw [-<] (2+1/2,0)--(2+1/2,1/4);
\draw (3,0)--(2+1/2,1/2);
\draw [-<] (3,0)--(2+3/4,1/4);
\draw (3/4,1/2)--(2+1/2,1/2);
\draw (1,1/2)--(1,3/4);
\draw [->] (1,3/4)--(1,5/8);
\draw (1+1/2,1/2)--(1+1/2,3/4);
\draw [->] (1+1/2,3/4)--(1+1/2,5/8);
\draw (1+3/4,1/2)--(1+3/4,3/4);
\draw [->] (1+3/4,3/4)--(1+3/4,5/8);
\draw (2+1/4,1/2)--(2+1/4,3/4);
\draw [->] (2+1/4,3/4)--(2+1/4,5/8);
\node [below] at (1/4,0) {$-$};
\node [below] at (3/4,0) {$0$};
\node [below] at (5/4,0) {$0$};
\node [below] at (2,0) {$0$};
\node [below] at (2+1/2,0) {$0$};
\node [below] at (3,0) {$+$};
\draw [line width=1.25, dotted] (1.4,-1/4)--(1.9,-1/4);
\draw [line width=1.25, dotted] (1.4,1/4)--(1.9,1/4);
\draw [line width=1.25, dotted] (1.5,5/8)--(1.75,5/8);
\draw (3/4,3/4)--(3/4,5/4);
\draw (3/4,1)--(6/4,1);
\draw (7/4,1)--(2+1/2,1);
\draw (5/2,3/4)--(5/2,5/4);
\node at (13/8,1) {$k$};
\end{tikzpicture}=q^{3k-2}\begin{tikzpicture}[baseline=.5ex]
\draw [line width=1.5, ->] (1/2,0)--(2+3/4,0);
\draw (1,0)--(1,3/4);
\draw [-<] (1,0)--(1,3/8);
\draw (2+1/4,0)--(2+1/4,3/4);
\draw [-<] (2+1/4,0)--(2+1/4,3/8);
\node [below] at (1,0) {$0$};
\node [below] at (2+1/4,0) {$0$};
\draw [line width=1.25, dotted] (1+1/4,3/8)--(2,3/8);
\draw [line width=1.25, dotted] (1+1/4,-1/4)--(2,-1/4);
\draw (3/4,3/4)--(3/4,5/4);
\draw (3/4,1)--(6/4,1);
\draw (7/4,1)--(2+1/2,1);
\draw (5/2,3/4)--(5/2,5/4);
\node at (13/8,1) {$k$};
\end{tikzpicture}
\end{equation}

Relation (S) allows one to switch two circles of opposite orientations whenever the two circles bound an annulus. We see from \cite{SW07} that relation (S) will be necessary for our list of reduction rules to be confluent, as none of the left sides of the defining relations are applicable to the diagrams in (S) unless they happen to bound a disk. We borrow notation from \cite{FS20} to say that two circles that bound an annulus on the surface and are oriented inconsistently with the boundary of the annulus form a \textit{British highway}. For example, the two circles on the left side of the relation (S) form a British highway. The fact that we are using oriented surfaces allows us to declare the right side of (S) to be the more reduced side. The relation (S) will serve as a reduction rule that will decrease the number of British highways on any connected component that is not a torus. The torus provides an exception since parallel nontrivial circles will bound two distinct annuli. See the remark after Theorem 2 regarding this exception. 

\begin{proposition}
i) The relations (S) hold in $S_q^{SL_3}(\Sigma)$ for any annulus embedded in $\Sigma.$ 

ii) The relations ($C_k$) hold in $\mathcal{S}_q^{SL_3}(\Sigma)$ for all $k \geq 0.$ 
\end{proposition}

\begin{proof}
i) (S) represents an isotopy of webs in the thickened surface $\Sigma \times (0,1)$, so the relation holds in $S_q^{SL_3}(\Sigma).$

ii) We will proceed by induction on $k.$ ($C_0$) is the same as (B4), so the statement is true for $k=0.$

If $k>0$ we can apply the relation (j) to the horizontal bar to the right of the top left strand to yield

\begin{equation*}
-q^{4/3}\sum_{b \in \{-,0,+\}} q^{-2b} \begin{tikzpicture}[baseline=.5ex]
\draw [line width=1.5, ->] (-3/2,0)--(3.5,0);
\draw (1/4,0)--(3/4,1/2);
\draw [-<] (1/4,0)--(1/2,1/4);
\draw (3/4,0)--(3/4,1/2);
\draw [-<] (3/4,0)--(3/4,1/4);
\draw (5/4,0)--(5/4,1/2);
\draw [-<] (5/4,0)--(5/4,1/4);
\draw (2,0)--(2,1/2);
\draw [-<] (2,0)--(2,1/4);
\draw (2+1/2,0)--(2+1/2,1/2);
\draw [-<] (2+1/2,0)--(2+1/2,1/4);
\draw (3,0)--(2+1/2,1/2);
\draw [-<] (3,0)--(2+3/4,1/4);
\draw (3/4,1/2)--(2+1/2,1/2);
\draw (1,1/2)--(1,3/4);
\draw [->] (1,3/4)--(1,5/8);
\draw (1+1/2,1/2)--(1+1/2,3/4);
\draw [->] (1+1/2,3/4)--(1+1/2,5/8);
\draw (1+3/4,1/2)--(1+3/4,3/4);
\draw [->] (1+3/4,3/4)--(1+3/4,5/8);
\draw (2+1/4,1/2)--(2+1/4,3/4);
\draw [->] (2+1/4,3/4)--(2+1/4,5/8);
\node [below] at (1/4,0) {$b$};
\node [below] at (3/4,0) {$0$};
\node [below] at (5/4,0) {$0$};
\node [below] at (2,0) {$0$};
\node [below] at (2+1/2,0) {$0$};
\node [below] at (3,0) {$+$};
\draw [line width=1.25, dotted] (1.4,-1/4)--(1.9,-1/4);
\draw [line width=1.25, dotted] (1.4,1/4)--(1.9,1/4);
\draw [line width=1.25, dotted] (1.5,5/8)--(1.75,5/8);
\draw (3/4,3/4)--(3/4,5/4);
\draw (3/4,1)--(2+1/2,1);
\draw (5/2,3/4)--(5/2,5/4);
\node at (13/8,5/4) {$k-1$};
\draw (0,0)--(-1/2,1/2);
\draw (-1/2,1/2)--(-3/4,1/2);
\draw (-3/4,0)--(-3/4,1/2);
\draw (-5/4,0)--(-3/4,1/2);
\draw (-1/2,1/2)--(-1/2,3/4);
\node [below] at (-5/4,0) {$-$};
\node [below] at (-3/4,0) {$0$};
\node [below] at (-1/8,0) {$-b$};
\draw [->] (0,0)--(-1/4,1/4);
\draw [-<] (-3/4,0)--(-3/4,1/4);
\draw [-<] (-5/4,0)--(-1,1/4);
\draw [->] (-1/2,3/4)--(-1/2,5/8);
\end{tikzpicture}
\end{equation*}

When $b=0$ the right connected component of the diagram is zero by relation (B3). When $b=+$ we compute that the left portion of the diagram reduces to

\begin{align*}
\begin{tikzpicture}[baseline=.5ex]
\draw [line width=1.5, ->] (-3/2,0)--(1/4,0);
\draw (0,0)--(-1/2,1/2);
\draw (-1/2,1/2)--(-3/4,1/2);
\draw (-3/4,0)--(-3/4,1/2);
\draw (-5/4,0)--(-3/4,1/2);
\draw (-1/2,1/2)--(-1/2,3/4);
\node [below] at (-5/4,0) {$-$};
\node [below] at (-3/4,0) {$0$};
\node [below] at (-1/8,0) {$-$};
\draw [->] (0,0)--(-1/4,1/4);
\draw [-<] (-3/4,0)--(-3/4,1/4);
\draw [-<] (-5/4,0)--(-1,1/4);
\draw [->] (-1/2,3/4)--(-1/2,5/8);
\end{tikzpicture}&=-q^{-1/3+1}\begin{tikzpicture}[baseline=.5ex]
\draw [line width=1.5, ->] (-1,0)--(1,0);
\draw (-3/4,0)--(-1/2,1/2);
\draw (-1/4,0)--(-1/2,1/2);
\draw (-1/2,1/2)--(1/2,1/2);
\draw (1/4,0)--(1/2,1/2);
\draw (3/4,0)--(1/2,1/2);
\draw (0,1/2)--(0,3/4);
\node [below] at (-3/4,0) {$-$};
\node [below] at (-1/4,0) {$0$};
\node [below] at (1/4,0) {$-$};
\node [below] at (3/4,0) {$0$};
\draw [-<] (-3/4,0)--(-5/8,1/4);
\draw [-<] (-1/4,0)--(-3/8,1/4);
\draw [-<] (1/4,0)--(3/8,1/4);
\draw [-<] (3/4,0)--(5/8,1/4);
\draw [->] (0,3/4)--(0,5/8);
\end{tikzpicture}\\
&=-q^{-1/3+1}(q^{-1}\begin{tikzpicture}[baseline=.5ex]
\draw [line width=1.5, ->] (-1,0)--(1,0);
\draw (-3/4,0)--(-1/2,1/2);
\draw (-1/4,0)--(-1/2,1/2);
\draw (-1/2,1/2)--(1/2,1/2);
\draw (1/4,0)--(1/2,1/2);
\draw (3/4,0)--(1/2,1/2);
\draw (0,1/2)--(0,3/4);
\node [below] at (-3/4,0) {$-$};
\node [below] at (-1/4,0) {$-$};
\node [below] at (1/4,0) {$0$};
\node [below] at (3/4,0) {$0$};
\draw [-<] (-3/4,0)--(-5/8,1/4);
\draw [-<] (-1/4,0)--(-3/8,1/4);
\draw [-<] (1/4,0)--(3/8,1/4);
\draw [-<] (3/4,0)--(5/8,1/4);
\draw [->] (0,3/4)--(0,5/8);
\end{tikzpicture}+q^{-3}\begin{tikzpicture}[baseline=.5ex]
\draw [line width=1.5,->] (-1,0)--(1,0);
\draw (-3/4,0)--(-3/4,3/4);
\draw (-3/4,3/4)--(-1/2,5/4);
\draw (-1/4,0)--(0,1/4);
\draw (1/4,0)--(0,1/4);
\draw (0,1/4)--(0,1/2);
\draw (0,1/2)--(-1/2,5/4);
\draw (0,1/2)--(1/2,5/4);
\draw (3/4,0)--(3/4,3/4);
\draw (3/4,3/4)--(1/2,5/4);
\draw (-1/2,5/4)--(1/2,5/4);
\draw (0,5/4)--(0,3/2);
\node [below] at (-3/4,0) {$-$};
\node [below] at (-1/4,0) {$-$};
\node [below] at (1/4,0) {$0$};
\node [below] at (3/4,0) {$0$};
\draw [-<] (-3/4,0)--(-3/4,1/2);
\draw [-<] (3/4,0)--(3/4,1/2);
\draw [-<] (-1/4,0)--(-1/8,1/8);
\draw [-<] (1/4,0)--(1/8,1/8);
\draw [->] (-1/2,5/4)--(-1/4,7/8);
\draw [->] (1/2,5/4)--(1/4,7/8);
\draw [->] (0,6/4)--(0,11/8);
\end{tikzpicture})
\end{align*}

Both of the last terms reduce to $0$ using (B3) after applying (I2) to the second diagram.

When $b=-$ we are interested in computing

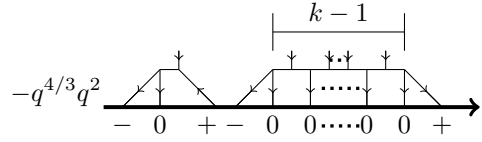
\begin{equation*}
-q^{4/3} q^{2} \begin{tikzpicture}[baseline=.5ex]
\draw [line width=1.5, ->] (-3/2,0)--(3.5,0);
\draw (1/4,0)--(3/4,1/2);
\draw [-<] (1/4,0)--(1/2,1/4);
\draw (3/4,0)--(3/4,1/2);
\draw [-<] (3/4,0)--(3/4,1/4);
\draw (5/4,0)--(5/4,1/2);
\draw [-<] (5/4,0)--(5/4,1/4);
\draw (2,0)--(2,1/2);
\draw [-<] (2,0)--(2,1/4);
\draw (2+1/2,0)--(2+1/2,1/2);
\draw [-<] (2+1/2,0)--(2+1/2,1/4);
\draw (3,0)--(2+1/2,1/2);
\draw [-<] (3,0)--(2+3/4,1/4);
\draw (3/4,1/2)--(2+1/2,1/2);
\draw (1,1/2)--(1,3/4);
\draw [->] (1,3/4)--(1,5/8);
\draw (1+1/2,1/2)--(1+1/2,3/4);
\draw [->] (1+1/2,3/4)--(1+1/2,5/8);
\draw (1+3/4,1/2)--(1+3/4,3/4);
\draw [->] (1+3/4,3/4)--(1+3/4,5/8);
\draw (2+1/4,1/2)--(2+1/4,3/4);
\draw [->] (2+1/4,3/4)--(2+1/4,5/8);
\node [below] at (1/4,0) {$-$};
\node [below] at (3/4,0) {$0$};
\node [below] at (5/4,0) {$0$};
\node [below] at (2,0) {$0$};
\node [below] at (2+1/2,0) {$0$};
\node [below] at (3,0) {$+$};
\draw [line width=1.25, dotted] (1.4,-1/4)--(1.9,-1/4);
\draw [line width=1.25, dotted] (1.4,1/4)--(1.9,1/4);
\draw [line width=1.25, dotted] (1.5,5/8)--(1.75,5/8);
\draw (3/4,3/4)--(3/4,5/4);
\draw (3/4,1)--(2+1/2,1);
\draw (5/2,3/4)--(5/2,5/4);
\node at (13/8,5/4) {$k-1$};
\draw (0,0)--(-1/2,1/2);
\draw (-1/2,1/2)--(-3/4,1/2);
\draw (-3/4,0)--(-3/4,1/2);
\draw (-5/4,0)--(-3/4,1/2);
\draw (-1/2,1/2)--(-1/2,3/4);
\node [below] at (-5/4,0) {$-$};
\node [below] at (-3/4,0) {$0$};
\node [below] at (-1/8,0) {$+$};
\draw [->] (0,0)--(-1/4,1/4);
\draw [-<] (-3/4,0)--(-3/4,1/4);
\draw [-<] (-5/4,0)--(-1,1/4);
\draw [->] (-1/2,3/4)--(-1/2,5/8);
\end{tikzpicture}
\end{equation*}

The right part of the diagram can be reduced by induction now while the left part of the diagram can be computed in the following manner:

\begin{align*}
\begin{tikzpicture}[baseline=.5ex]
\draw [line width=1.5, ->] (-3/2,0)--(1/4,0);
\draw (0,0)--(-1/2,1/2);
\draw (-1/2,1/2)--(-3/4,1/2);
\draw (-3/4,0)--(-3/4,1/2);
\draw (-5/4,0)--(-3/4,1/2);
\draw (-1/2,1/2)--(-1/2,3/4);
\node [below] at (-5/4,0) {$-$};
\node [below] at (-3/4,0) {$0$};
\node [below] at (-1/8,0) {$+$};
\draw [->] (0,0)--(-1/4,1/4);
\draw [-<] (-3/4,0)--(-3/4,1/4);
\draw [-<] (-5/4,0)--(-1,1/4);
\draw [->] (-1/2,3/4)--(-1/2,5/8);
\end{tikzpicture}&\stackrel{(B1)}{=} -q^{1/3-1} \begin{tikzpicture}[baseline=.5ex]
\draw [line width=1.5,->] (-1/2,0)--(1/2,0);
\draw (-1/4,0)--(0,1/2);
\draw [->] (-1/4,0)--(-1/8,1/4);
\draw (1/4,0)--(0,1/2);
\draw [->] (1/4,0)--(1/8,1/4);
\draw (0,1/2)--(0,3/4);
\draw [-<] (0,1/2)--(0,5/8);
\node [below] at (-1/4,0) {$-$};
\node [below] at (1/4,0) {$+$};
\end{tikzpicture} \\
&\stackrel{(h)}{=}-q^{1/3-1} (-q^{4/3}) \sum_{a \in \{-,0,+\}} \begin{tikzpicture}[baseline=.5ex]
\draw [line width=1.5,->] (-1,0)--(3/4,0);
\node [below] at (-3/4-1/8,0) {$-a$};
\draw (-3/4,0)--(-3/4,3/4);
\draw [-<] (-3/4,0)--(-3/4,3/8);
\draw (-1/2,0)--(0,1/2);
\draw [->] (-1/2,0)--(-1/4,1/4);
\draw (0,0)--(0,1/2);
\draw [->] (0,0)--(0,1/4);
\draw (1/2,0)--(0,1/2);
\draw [->] (1/2,0)--(1/4,1/4);
\node [below] at (-1/2,-.07) {$a$};
\node [below] at (0,0) {$-$};
\node [below] at (1/2,0) {$+$};
\end{tikzpicture}\\
&\stackrel{(k)}{=} -q^{1/3-1}(-q^{4/3})q^{-2}(-q) \begin{tikzpicture}[baseline=.5ex]
\draw [line width=1.5,->] (-1/4,0)--(1/4,0);
\draw (0,0)--(0,3/4);
\draw [-<] (0,0)--(0,3/8);
\node [below] at (0,0) {$0$};
\end{tikzpicture}
\end{align*}

This all reduces to

\begin{align*}
-q^{1/3-1}(-q^{4/3})q^{-2}(-q)(-q^{4/3}q^2)q^{3(k-1)-2}\begin{tikzpicture}[baseline=.5ex]
\draw [line width=1.5, ->] (1/2,0)--(2+3/4,0);
\draw (1,0)--(1,3/4);
\draw [-<] (1,0)--(1,3/8);
\draw (2+1/4,0)--(2+1/4,3/4);
\draw [-<] (2+1/4,0)--(2+1/4,3/8);
\node [below] at (1,0) {$0$};
\node [below] at (2+1/4,0) {$0$};
\draw [line width=1.25, dotted] (1+1/4,3/8)--(2,3/8);
\draw [line width=1.25, dotted] (1+1/4,-1/4)--(2,-1/4);
\draw (3/4,3/4)--(3/4,5/4);
\draw (3/4,1)--(6/4,1);
\draw (7/4,1)--(2+1/2,1);
\draw (5/2,3/4)--(5/2,5/4);
\node at (13/8,1) {$k$};
\end{tikzpicture}\\
= q^{3k-2}\begin{tikzpicture}[baseline=.5ex]
\draw [line width=1.5, ->] (1/2,0)--(2+3/4,0);
\draw (1,0)--(1,3/4);
\draw [-<] (1,0)--(1,3/8);
\draw (2+1/4,0)--(2+1/4,3/4);
\draw [-<] (2+1/4,0)--(2+1/4,3/8);
\node [below] at (1,0) {$0$};
\node [below] at (2+1/4,0) {$0$};
\draw [line width=1.25, dotted] (1+1/4,3/8)--(2,3/8);
\draw [line width=1.25, dotted] (1+1/4,-1/4)--(2,-1/4);
\draw (3/4,3/4)--(3/4,5/4);
\draw (3/4,1)--(6/4,1);
\draw (7/4,1)--(2+1/2,1);
\draw (5/2,3/4)--(5/2,5/4);
\node at (13/8,1) {$k$};
\end{tikzpicture}
\end{align*}

which concludes the proof by induction.

\end{proof}

For the rest of this section, we will assume any boundary arcs in our diagram have an orientation that matches the one appearing in the pictures of the defining boundary relations and that this orientation dictates the height order.

A univalent endpoint of a web diagram is a $\textit{bad endpoint}$ if the strand attached to the endpoint is oriented out of the boundary. For example, the endpoint in the picture on the left of relation (B1) is a bad endpoint while the two endpoints on the right of the relation are good. We say that a pair of two good endpoints on the same boundary arc with states $b$ and $a$ are a $\textit{bad pair}$ if $b>a$ but the endpoint with state $b$ is lower in the height order than the endpoint with state $a$. For example, the two endpoints on the left of (B2) form a bad pair, while the two endpoints in each diagram of the right side of the relation form a good pair. In the following, by the term $\textit{vertices}$ we mean only trivalent vertices of the web.

\begin{definition}
We define the $\textit{complexity}$ of a stated web diagram to be the tuple (\text{\#crossings}, \text{\#bad endpoints}, \text{\#bad pairs}, \text{\#vertices}, \text{\#connected components}, \text{\#British highways}) in $ \mathbb{Z}_{\geq 0}^6$. 
\end{definition} 

We use the lexigocraphic ordering on $\mathbb{Z}_{\geq 0}^6$ and note that each defining relation, each relation ($C_k$), and each relation (S) involve a single diagram on the left side of the equation while the right side of the equation contains only diagrams of strictly lower complexity than the one on the left side of the equation.

We say that a diagram contains a $\textit{reducible feature}$ if the left side of one of the relations (I1a)-(I4b), (B1)-(B4), ($C_k$), or (S) applies. If a diagram contains no reducible feature, we call such a diagram an $\textit{irreducible diagram}.$

\begin{theorem}
The set of isotopy classes of irreducible diagrams on $\Sigma$ forms a basis for $S_q^{SL_3}(\Sigma)$.
\end{theorem}

\begin{remark} If $\Sigma$ has a connected component that is a torus, we modify our notion of an irreducible diagram. By omitting the reduction rule (S) on any torus, the proof below can be modified to show that the remaining reduction rules will produce a basis consisting of the set of irreducible diagrams up to isotopy and circle flip moves (S) on any torus.
\end{remark}

\begin{proof}
We will apply the Diamond Lemma in much the same setup as \cite{Le18}. First, we claim that module freely spanned by isotopy classes of web diagrams with our chosen boundary orientations modulo the defining relations along with ($C_k$) and (S) yields a module isomorphic to $S_q^{SL_3}(\Sigma)$. To do this, one observes that ribbon Reidemeister moves RI, RII, and RIII and the fact that a strand can pass over or under a vertex all follow from the defining interior relations, as shown in \cite{Kup94}. The fact that ($C_k$) and (S) are redundant relations completes this part of the argument.

Next, we must verify that given a diagram $D,$ the process of iteratively applying the left sides of our relations to $D$ and to its descendants always terminates in a linear combination of irreducible diagrams. This is guaranteed by the fact that our reduction rules involve replacing a diagram by a linear combination of diagrams of strictly lower complexity in our lexicographic ordering, as in Theorem 2.2 of \cite{SW07}. Thus, the set of isotopy classes of irreducible diagrams span $\mathcal{S}_q^{SL_3}(\Sigma).$

To show that each diagram can be uniquely written as a linear combination of irreducible diagrams, we must show the local confluence of our relations. This is the reason that we had to include the redundant relations ($C_k$) and (S). We must check that if more than one relation is applicable to a diagram then we can reach a common descendant regardless of which relation we choose to apply. We use the same notion of the \textit{support} of a relation as \cite{Le18}. If two relations are applicable to a diagram, but their support is disjoint, then the applications of these relations commute with each other, and thus immediately reach a common descendant.

We must find local confluence for relations whose supports overlap nontrivially. If the two relations are both interior relations or (S), then we see by \cite{SW07} that they are locally confluent.

There is one possible way for a the support of an interior relation to intersect the support of a boundary relation: a square could be connected to the top of the relations ($C_k$) for some $k\geq 2$. The following diagram shows an example of an overlap of ($C_4$) and (I2).

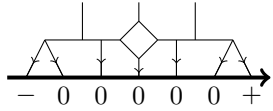
\begin{equation*}
\begin{tikzpicture}
\draw [->, line width=1.5] (0,0)--(14/4,0);
\draw  (1/4,0)--(1/2,1/2);
\draw [-<] (1/4,0)--(1/4+1/8,1/4);
\draw (3/4,0)--(1/2,1/2);
\draw [-<] (3/4,0)--(3/4-1/8,1/4);
\draw (5/4,0)--(5/4,1/2);
\draw [-<] (5/4,0)--(5/4,1/4);
\draw (7/4,0)--(7/4,1/4);
\draw [-<] (7/4,0)--(7/4,1/8);
\draw (9/4,0)--(9/4,1/2);
\draw [-<] (9/4,0)--(9/4,1/4);
\draw (11/4,0)--(3,1/2);
\draw [-<] (11/4,0)--(11/4+1/8,1/4);
\draw (13/4,0)--(3,1/2);
\draw [-<] (13/4,0)--(13/4-1/8,1/4);
\draw (1/2,1/2)--(3/2,1/2);
\draw (2,1/2)--(3,1/2);
\draw (3/2,1/2)--(7/4,1/4);
\draw (2,1/2)--(7/4,1/4);
\draw (1,1/2)--(1,1);
\draw (3/2,1/2)--(7/4,3/4);
\draw (2,1/2)--(7/4,3/4);
\draw (5/2,1/2)--(5/2,1);
\draw (7/4,3/4)--(7/4,1);
\node [below] at (1/4,0) {$-$};
\node [below] at (3/4,0) {$0$};
\node [below] at (5/4,0) {$0$};
\node [below] at (7/4,0) {$0$};
\node [below] at (9/4,0) {$0$};
\node [below] at (11/4,0) {$0$};
\node [below] at (13/4,0) {$+$};
\end{tikzpicture}
\end{equation*}

Such a situation will terminate at $0$ no matter which relation ($C_k$) or (I2) is applied first, as each resulting diagram will provide an opportunity to apply (B3).

Finally, we consider the cases of overlapping supports of the defining boundary relations and the additional relations ($C_k$). A first easy case is an overlap of (B3) with (B3), which must be of the following form:

\begin{equation*}
\begin{tikzpicture}[baseline=.5ex]
\draw (0,1/2)--(0,0);
\draw [->] (0,1/2)--(0,1/4);
\draw (1/4,0)--(0,1/2);
\draw [->] (0,1/2)--(1/8,1/4);
\draw (-1/4,0)--(0,1/2);
\draw [->] (0,1/2)--(-1/8,1/4);
\draw [line width=1.5,->] (-1/2,0)--(1/2,0);
\node [below] at (-1/4,0) {$a$};
\node [below] at (0,0) {$a$};
\node [below] at (1/4,0) {$a$};
\end{tikzpicture} \text{($a \in \{-,0,+\}$)} 
\end{equation*}

Applying (B3) to either the left triangle or the right triangle in the above diagram yields zero.

We see that the only other supports that can overlap are those of (B2) with any of (B2), (B3), (B4), and ($C_k$).

\underline{(B2) and (B2)}:

If (B2) overlaps with (B2): the overlap must be of the following form.
\begin{equation*}
\begin{tikzpicture}
\draw [line width=1.5, ->] (-1/2,0)--(1/2,0);
\draw (0,0)--(0,3/4);
\draw [-<] (0,0)--(0,3/8);
\draw (1/4,0)--(1/4,3/4);
\draw [-<] (1/4,0)--(1/4,3/8);
\draw (-1/4,0)--(-1/4,3/4);
\draw [-<] (-1/4,0)--(-1/4,3/8);
\node [below] at (-1/4,0) {$+$};
\node [below] at (0,0) {$0$};
\node [below] at (1/4,0) {$-$};
\end{tikzpicture}
\end{equation*}

If we first apply (B2) to the right two endpoints, and then we continue to apply (B2) until there are no longer any bad pairs we obtain:

\begin{align*}
q^{-3} \begin{tikzpicture}[baseline=.5ex]
\draw [line width=1.5, ->] (-1/2,0)--(1/2,0);
\draw (0,0)--(0,3/4);
\draw [-<] (0,0)--(0,3/8);
\draw (1/4,0)--(1/4,3/4);
\draw [-<] (1/4,0)--(1/4,3/8);
\draw (-1/4,0)--(-1/4,3/4);
\draw [-<] (-1/4,0)--(-1/4,3/8);
\node [below] at (-1/4,0) {$-$};
\node [below] at (0,0) {$0$};
\node [below] at (1/4,0) {$+$};
\end{tikzpicture}+2q^{-5}\begin{tikzpicture}[baseline=.5ex]
\draw [line width=1.5,->] (-1/2,0)--(1/2,0);
\draw (-1/4,0)--(-1/4,3/4);
\draw [-<] (-1/4,0)--(-1/4,3/8);
\draw (0,0)--(1/8,1/4);
\draw [-<] (0,0)--(1/16,1/8);
\draw (1/8,1/4)--(1/8,1/2);
\draw [->] (1/8,1/4)--(1/8,3/8);
\draw (1/4,0)--(1/8,1/4);
\draw [-<] (1/4,0)--(3/16,1/8);
\draw (1/8,1/2)--(0,3/4);
\draw [-<] (1/8,1/2)--(1/16,5/8);
\draw (1/8,1/2)--(1/4,3/4);
\draw [-<] (1/8,1/2)--(3/16,5/8);
\node [below] at (-1/4,0) {$-$};
\node [below] at (0,0) {$0$};
\node [below] at (1/4,0) {$+$};
\end{tikzpicture}+q^{-5}\begin{tikzpicture}[baseline=.5ex]
\draw [line width=1.5,->] (-1/4,0)--(3/4,0);
\draw (1/2,0)--(1/2,3/4);
\draw [-<] (1/2,0)--(1/2,3/8);
\draw (0,0)--(1/8,1/4);
\draw [-<] (0,0)--(1/16,1/8);
\draw (1/8,1/4)--(1/8,1/2);
\draw [->] (1/8,1/4)--(1/8,3/8);
\draw (1/4,0)--(1/8,1/4);
\draw [-<] (1/4,0)--(3/16,1/8);
\draw (1/8,1/2)--(0,3/4);
\draw [-<] (1/8,1/2)--(1/16,5/8);
\draw (1/8,1/2)--(1/4,3/4);
\draw [-<] (1/8,1/2)--(3/16,5/8);
\node [below] at (0,0) {$-$};
\node [below] at (1/4,0) {$0$};
\node [below] at (1/2,0) {$+$};
\end{tikzpicture}+q^{-7}\begin{tikzpicture}[baseline=.5ex]
\draw [line width=1.5, ->] (-1/4,0)--(3/4,0);
\draw (0,0)--(0,3/4);
\draw [-<] (0,0)--(0,3/8);
\draw (1/4,0)--(3/8,1/4);
\draw [-<] (1/4,0)--(5/16,1/8);
\draw (1/2,0)--(3/8,1/4);
\draw [-<] (1/2,0)--(7/16,1/8);
\draw (3/8,1/4)--(3/8,1/2);
\draw [->] (3/8, 1/4)--(3/8,3/8);
\draw (0,3/4)--(1/8,1);
\draw (1/8,1)--(3/8,1/2);
\draw [->] (1/8,1)--(1/4,3/4);
\draw (3/8,1/2)--(1/2,3/4);
\draw (1/2,3/4)--(1/2,3/2);
\draw [-<] (1/2,3/4)--(1/2,5/4);
\draw (1/8,1)--(1/8,5/4);
\draw [->] (1/8,1)--(1/8,9/8);
\draw (1/8,5/4)--(0,3/2);
\draw [-<] (1/8,5/4)--(1/16,11/8);
\draw (1/8,5/4)--(1/4,3/2);
\draw [-<] (1/8,5/4)--(3/16,11/8);
\node [below] at (0,0) {$-$};
\node [below] at (1/4,0) {$0$};
\node [below] at (1/2,0) {$+$};
\end{tikzpicture}\\+q^{-7} \begin{tikzpicture}[baseline=.5ex]
\draw [line width=1.5, ->] (-1/4,0)--(3/4,0);
\draw (0,0)--(0,3/2);
\draw [-<] (0,0)--(0,3/8);
\draw (1/4,0)--(3/8,1/4);
\draw [-<] (1/4,0)--(5/16,1/8);
\draw (1/2,0)--(3/8,1/4);
\draw [-<] (1/2,0)--(7/16,1/8);
\draw (3/8,1/4)--(3/8,1/2);
\draw [->] (3/8, 1/4)--(3/8,3/8);
\draw (3/8,1/2)--(1/4,3/4);
\draw [<-] (1/4,3/4)--(3/8,1);
\draw (3/8,1/2)--(1/2,3/4);
\draw [<-] (1/2,3/4)--(3/8,1);
\draw (3/8,1)--(3/8,5/4);
\draw [->] (3/8,1)--(3/8,9/8);
\draw (3/8,5/4)--(1/4,3/2);
\draw [-<] (3/8,5/4)--(5/16,11/8);
\draw (3/8,5/4)--(1/2,3/2);
\draw [-<] (3/8,5/4)--(7/16,11/8);
\node [below] at (0,0) {$-$};
\node [below] at (1/4,0) {$0$};
\node [below] at (1/2,0) {$+$};
\end{tikzpicture}+q^{-7} \begin{tikzpicture}[baseline=.5ex, xscale=-1]
\draw [line width=1.5, <-] (-1/4,0)--(3/4,0);
\draw (0,0)--(0,3/4);
\draw [-<] (0,0)--(0,3/8);
\draw (1/4,0)--(3/8,1/4);
\draw [-<] (1/4,0)--(5/16,1/8);
\draw (1/2,0)--(3/8,1/4);
\draw [-<] (1/2,0)--(7/16,1/8);
\draw (3/8,1/4)--(3/8,1/2);
\draw [->] (3/8, 1/4)--(3/8,3/8);
\draw (0,3/4)--(1/8,1);
\draw (1/8,1)--(3/8,1/2);
\draw [->] (1/8,1)--(1/4,3/4);
\draw (3/8,1/2)--(1/2,3/4);
\draw (1/2,3/4)--(1/2,3/2);
\draw [-<] (1/2,3/4)--(1/2,5/4);
\draw (1/8,1)--(1/8,5/4);
\draw [->] (1/8,1)--(1/8,9/8);
\draw (1/8,5/4)--(0,3/2);
\draw [-<] (1/8,5/4)--(1/16,11/8);
\draw (1/8,5/4)--(1/4,3/2);
\draw [-<] (1/8,5/4)--(3/16,11/8);
\node [below] at (0,0) {$+$};
\node [below] at (1/4,0) {$0$};
\node [below] at (1/2,0) {$-$};
\end{tikzpicture}+q^{-9}\begin{tikzpicture}[baseline=.5ex]
\draw [line width=1.5, ->] (-1/4,0)--(3/4,0);
\draw (0,0)--(0,3/4);
\draw [-<] (0,0)--(0,3/8);
\draw (1/4,0)--(3/8,1/4);
\draw [-<] (1/4,0)--(5/16,1/8);
\draw (1/2,0)--(3/8,1/4);
\draw [-<] (1/2,0)--(7/16,1/8);
\draw (3/8,1/4)--(3/8,1/2);
\draw [->] (3/8, 1/4)--(3/8,3/8);
\draw (0,3/4)--(1/8,1);
\draw (1/8,1)--(3/8,1/2);
\draw [->] (1/8,1)--(1/4,3/4);
\draw (3/8,1/2)--(1/2,3/4);
\draw (1/2,3/4)--(1/2,3/2);
\draw [-<] (1/2,3/4)--(1/2,5/4);
\draw (1/8,1)--(1/8,5/4);
\draw [->] (1/8,1)--(1/8,9/8);
\draw (1/8,5/4)--(0,3/2);
\draw [-<] (1/8,5/4)--(1/16,11/8);
\draw (1/8,5/4)--(1/4,3/2);
\draw [-<] (1/8,5/4)--(3/16,11/8);
\draw (1/2,3/2)--(3/8,7/4);
\draw (1/4,3/2)--(3/8,7/4);
\draw (3/8,7/4)--(3/8,2);
\draw (3/8,2)--(1/4,9/4);
\draw (3/8,2)--(1/2,9/4);
\draw (0,3/2)--(0,9/4);
\node [below] at (0,0) {$-$};
\node [below] at (1/4,0) {$0$};
\node [below] at (1/2,0) {$+$};
\end{tikzpicture}
\end{align*}

\begin{align*}
\stackrel{(I3),(I2),(B4)}{=} q^{-3} \begin{tikzpicture}[baseline=.5ex]
\draw [line width=1.5, ->] (-1/2,0)--(1/2,0);
\draw (0,0)--(0,3/4);
\draw [-<] (0,0)--(0,3/8);
\draw (1/4,0)--(1/4,3/4);
\draw [-<] (1/4,0)--(1/4,3/8);
\draw (-1/4,0)--(-1/4,3/4);
\draw [-<] (-1/4,0)--(-1/4,3/8);
\node [below] at (-1/4,0) {$-$};
\node [below] at (0,0) {$0$};
\node [below] at (1/4,0) {$+$};
\end{tikzpicture}+q^{-5}\begin{tikzpicture}[baseline=.5ex]
\draw [line width=1.5,->] (-1/2,0)--(1/2,0);
\draw (-1/4,0)--(-1/4,3/4);
\draw [-<] (-1/4,0)--(-1/4,3/8);
\draw (0,0)--(1/8,1/4);
\draw [-<] (0,0)--(1/16,1/8);
\draw (1/8,1/4)--(1/8,1/2);
\draw [->] (1/8,1/4)--(1/8,3/8);
\draw (1/4,0)--(1/8,1/4);
\draw [-<] (1/4,0)--(3/16,1/8);
\draw (1/8,1/2)--(0,3/4);
\draw [-<] (1/8,1/2)--(1/16,5/8);
\draw (1/8,1/2)--(1/4,3/4);
\draw [-<] (1/8,1/2)--(3/16,5/8);
\node [below] at (-1/4,0) {$-$};
\node [below] at (0,0) {$0$};
\node [below] at (1/4,0) {$+$};
\end{tikzpicture}+q^{-5}\begin{tikzpicture}[baseline=.5ex]
\draw [line width=1.5,->] (-1/4,0)--(3/4,0);
\draw (1/2,0)--(1/2,3/4);
\draw [-<] (1/2,0)--(1/2,3/8);
\draw (0,0)--(1/8,1/4);
\draw [-<] (0,0)--(1/16,1/8);
\draw (1/8,1/4)--(1/8,1/2);
\draw [->] (1/8,1/4)--(1/8,3/8);
\draw (1/4,0)--(1/8,1/4);
\draw [-<] (1/4,0)--(3/16,1/8);
\draw (1/8,1/2)--(0,3/4);
\draw [-<] (1/8,1/2)--(1/16,5/8);
\draw (1/8,1/2)--(1/4,3/4);
\draw [-<] (1/8,1/2)--(3/16,5/8);
\node [below] at (0,0) {$-$};
\node [below] at (1/4,0) {$0$};
\node [below] at (1/2,0) {$+$};
\end{tikzpicture}+q^{-7}\begin{tikzpicture}[baseline=.5ex]
\draw [line width=1.5, ->] (-1/4,0)--(3/4,0);
\draw (0,0)--(0,3/4);
\draw [-<] (0,0)--(0,3/8);
\draw (1/4,0)--(3/8,1/4);
\draw [-<] (1/4,0)--(5/16,1/8);
\draw (1/2,0)--(3/8,1/4);
\draw [-<] (1/2,0)--(7/16,1/8);
\draw (3/8,1/4)--(3/8,1/2);
\draw [->] (3/8, 1/4)--(3/8,3/8);
\draw (0,3/4)--(1/8,1);
\draw (1/8,1)--(3/8,1/2);
\draw [->] (1/8,1)--(1/4,3/4);
\draw (3/8,1/2)--(1/2,3/4);
\draw (1/2,3/4)--(1/2,3/2);
\draw [-<] (1/2,3/4)--(1/2,5/4);
\draw (1/8,1)--(1/8,5/4);
\draw [->] (1/8,1)--(1/8,9/8);
\draw (1/8,5/4)--(0,3/2);
\draw [-<] (1/8,5/4)--(1/16,11/8);
\draw (1/8,5/4)--(1/4,3/2);
\draw [-<] (1/8,5/4)--(3/16,11/8);
\node [below] at (0,0) {$-$};
\node [below] at (1/4,0) {$0$};
\node [below] at (1/2,0) {$+$};
\end{tikzpicture}\\ +q^{-7} \begin{tikzpicture}[baseline=.5ex, xscale=-1]
\draw [line width=1.5, <-] (-1/4,0)--(3/4,0);
\draw (0,0)--(0,3/4);
\draw [-<] (0,0)--(0,3/8);
\draw (1/4,0)--(3/8,1/4);
\draw [-<] (1/4,0)--(5/16,1/8);
\draw (1/2,0)--(3/8,1/4);
\draw [-<] (1/2,0)--(7/16,1/8);
\draw (3/8,1/4)--(3/8,1/2);
\draw [->] (3/8, 1/4)--(3/8,3/8);
\draw (0,3/4)--(1/8,1);
\draw (1/8,1)--(3/8,1/2);
\draw [->] (1/8,1)--(1/4,3/4);
\draw (3/8,1/2)--(1/2,3/4);
\draw (1/2,3/4)--(1/2,3/2);
\draw [-<] (1/2,3/4)--(1/2,5/4);
\draw (1/8,1)--(1/8,5/4);
\draw [->] (1/8,1)--(1/8,9/8);
\draw (1/8,5/4)--(0,3/2);
\draw [-<] (1/8,5/4)--(1/16,11/8);
\draw (1/8,5/4)--(1/4,3/2);
\draw [-<] (1/8,5/4)--(3/16,11/8);
\node [below] at (0,0) {$+$};
\node [below] at (1/4,0) {$0$};
\node [below] at (1/2,0) {$-$};
\end{tikzpicture}+q^{-5}\begin{tikzpicture}[baseline=.5ex, xscale=2/3]
\draw (0,1/2)--(0,1);
\draw [-<] (0,1/2)--(0,3/4);
\draw (0,1/2)--(-1/2,1);
\draw [-<] (0,1/2)--(-1/4,3/4);
\draw (0,1/2)--(1/2,1);
\draw [-<] (0,1/2)--(1/4,3/4);
\draw [line width=1.5, ->] (-3/4,0)--(3/4,0);
\end{tikzpicture}.
\end{align*}

If, instead, we first apply (B2) to the left two endpoints, and then we continue to apply (B2) until there are no longer any bad pairs, we obtain the same linear combination but with the diagrams reflected in a vertical line (but with the state locations and boundary orientation remaining the same). By noting the coefficients in our last equation are symmetric with respect to this reflection, we see that we obtain the same answer in both cases.

\underline{(B2) and (B3):}

If (B2) overlaps with (B3): the overlap must take one of the following forms.

\begin{equation*}
\begin{tikzpicture}[baseline=.5ex]
\draw [line width=1.5, ->] (-1/4,0)--(3/4,0);
\draw (0,0)--(0,1/2);
\draw [-<] (0,0)--(0,1/4);
\draw (1/4,0)--(3/8,1/4);
\draw [-<] (1/4,0)--(5/16,1/8);
\draw (1/2,0)--(3/8,1/4);
\draw [-<] (1/2,0)--(7/16,1/8);
\draw (3/8,1/4)--(3/8,1/2);
\node [below] at (0,0) {$b$};
\node [below] at (1/4,-.07) {$a$};
\node [below] at (1/2,-.07) {$a$};
\end{tikzpicture} \text{ or } \begin{tikzpicture}[baseline=.5ex, xscale=-1]
\draw [line width=1.5, <-] (-1/4,0)--(3/4,0);
\draw (0,0)--(0,1/2);
\draw [-<] (0,0)--(0,1/4);
\draw (1/4,0)--(3/8,1/4);
\draw [-<] (1/4,0)--(5/16,1/8);
\draw (1/2,0)--(3/8,1/4);
\draw [-<] (1/2,0)--(7/16,1/8);
\draw (3/8,1/4)--(3/8,1/2);
\node [below] at (0,-.07) {$a$};
\node [below] at (1/4,0) {$b$};
\node [below] at (1/2,0) {$b$};
\end{tikzpicture} \text{ ($b>a$)}
\end{equation*}

Both cases are handled symmetrically, so we will focus on the left case. If we apply (B3) first, we obtain zero. So we must show that if we instead apply (B2) first we eventually obtain zero. We do this by computing:

\begin{align*}
\begin{tikzpicture}[baseline=.5ex]
\draw [line width=1.5, ->] (-1/4,0)--(3/4,0);
\draw (0,0)--(0,1/2);
\draw [-<] (0,0)--(0,1/4);
\draw (1/4,0)--(3/8,1/4);
\draw [-<] (1/4,0)--(5/16,1/8);
\draw (1/2,0)--(3/8,1/4);
\draw [-<] (1/2,0)--(7/16,1/8);
\draw (3/8,1/4)--(3/8,1/2);
\node [below] at (0,0) {$b$};
\node [below] at (1/4,-.07) {$a$};
\node [below] at (1/2,-.07) {$a$};
\end{tikzpicture} &\stackrel{(B2)}{=} q^{-1} \begin{tikzpicture}[baseline=.5ex]
\draw [line width=1.5, ->] (-1/4,0)--(3/4,0);
\draw (0,0)--(0,1/2);
\draw [-<] (0,0)--(0,1/4);
\draw (1/4,0)--(3/8,1/4);
\draw [-<] (1/4,0)--(5/16,1/8);
\draw (1/2,0)--(3/8,1/4);
\draw [-<] (1/2,0)--(7/16,1/8);
\draw (3/8,1/4)--(3/8,1/2);
\node [below] at (0,-.07) {$a$};
\node [below] at (1/4,0) {$b$};
\node [below] at (1/2,-.07) {$a$};
\end{tikzpicture}+q^{-3} \begin{tikzpicture}[baseline=.5ex]
\draw [line width=1.5, ->] (-1/4,0)--(3/4,0);
\draw (0,0)--(1/8,1/4);
\draw [-<] (0,0)--(1/16,1/8);
\draw (1/4,0)--(1/8,1/4);
\draw [-<] (1/4,0)--(3/16,1/8);
\draw (1/8,1/4)--(1/8,1/2);
\draw (1/8,1/2)--(0,3/4);
\draw (1/8,1/2)--(3/8,1);
\draw (1/2,0)--(1/2,3/4);
\draw [-<] (1/2,0)--(1/2,3/8);
\draw (1/2,3/4)--(3/8,1);
\draw (3/8,1)--(3/8,5/4);
\draw (0,3/4)--(0,5/4);
\node [below] at (0,-.07) {$a$};
\node [below] at (1/4,0) {$b$};
\node [below] at (1/2,-.07) {$a$};
\end{tikzpicture}\\
&\stackrel{(B2)}{=} q^{-2}\begin{tikzpicture}[baseline=.5ex]
\draw [line width=1.5, ->] (-1/4,0)--(3/4,0);
\draw (0,0)--(0,1/2);
\draw [-<] (0,0)--(0,1/4);
\draw (1/4,0)--(3/8,1/4);
\draw [-<] (1/4,0)--(5/16,1/8);
\draw (1/2,0)--(3/8,1/4);
\draw [-<] (1/2,0)--(7/16,1/8);
\draw (3/8,1/4)--(3/8,1/2);
\node [below] at (0,-.07) {$a$};
\node [below] at (1/4,-.07) {$a$};
\node [below] at (1/2,0) {$b$};
\end{tikzpicture}+q^{-4}
\begin{tikzpicture}[baseline=.5ex]
\draw [line width=1.5, ->] (-1/4,0)--(3/4,0);
\draw (0,0)--(0,5/4);
\draw [-<] (0,0)--(0,1/4);
\draw (1/4,0)--(3/8,1/4);
\draw [-<] (1/4,0)--(5/16,1/8);
\draw (1/2,0)--(3/8,1/4);
\draw [-<] (1/2,0)--(7/16,1/8);
\draw (3/8,1/4)--(3/8,1/2);
\draw (3/8,1/2)--(1/4,3/4);
\draw (3/8,1/2)--(1/2,3/4);
\draw (1/4,3/4)--(3/8,1);
\draw (1/2,3/4)--(3/8,1);
\draw (3/8,1)--(3/8,5/4);
\node [below] at (0,-.07) {$a$};
\node [below] at (1/4,-.07) {$a$};
\node [below] at (1/2,0) {$b$};
\end{tikzpicture}+q^{-4}\begin{tikzpicture}[baseline=.5ex]
\draw [line width=1.5, ->] (-1/4,0)--(3/4,0);
\draw (0,0)--(1/8,1/4);
\draw [-<] (0,0)--(1/16,1/8);
\draw (1/4,0)--(1/8,1/4);
\draw [-<] (1/4,0)--(3/16,1/8);
\draw (1/8,1/4)--(1/8,1/2);
\draw (1/8,1/2)--(0,3/4);
\draw (1/8,1/2)--(3/8,1);
\draw (1/2,0)--(1/2,3/4);
\draw [-<] (1/2,0)--(1/2,3/8);
\draw (1/2,3/4)--(3/8,1);
\draw (3/8,1)--(3/8,5/4);
\draw (0,3/4)--(0,5/4);
\node [below] at (0,-.07) {$a$};
\node [below] at (1/4,-.07) {$a$};
\node [below] at (1/2,0) {$b$};
\end{tikzpicture}+q^{-6}\begin{tikzpicture}[baseline=.5ex]
\draw [line width=1.5, ->] (-1/4,0)--(3/4,0);
\draw (0,0)--(0,3/4);
\draw [-<] (0,0)--(0,1/4);
\draw (1/4,0)--(3/8,1/4);
\draw [-<] (1/4,0)--(5/16,1/8);
\draw (1/2,0)--(3/8,1/4);
\draw [-<] (1/2,0)--(7/16,1/8);
\draw (3/8,1/4)--(3/8,1/2);
\draw (3/8,1/2)--(1/4,3/4);
\draw (3/8,1/2)--(1/2,3/4);
\draw (0,3/4)--(1/8,1);
\draw (1/4,3/4)--(1/8,1);
\draw (1/2,3/4)--(1/2,6/4);
\draw (1/8,1)--(1/8,5/4);
\draw (1/8,5/4)--(0,6/4);
\draw (1/8,5/4)--(1/4,6/4);
\draw (0,6/4)--(0,2);
\draw (1/4,6/4)--(3/8,7/4);
\draw (1/2,6/4)--(3/8,7/4);
\draw (3/8,7/4)--(3/8,2);
\node [below] at (0,-.07) {$a$};
\node [below] at (1/4,-.07) {$a$};
\node [below] at (1/2,0) {$b$};
\end{tikzpicture}\\
&\stackrel{(I2),(I3),(B3)}{=} (q^{-2}-q^{-4}q^3[2]+1) \begin{tikzpicture}[baseline=.5ex]
\draw [line width=1.5, ->] (-1/4,0)--(3/4,0);
\draw (0,0)--(0,1/2);
\draw [-<] (0,0)--(0,1/4);
\draw (1/4,0)--(3/8,1/4);
\draw [-<] (1/4,0)--(5/16,1/8);
\draw (1/2,0)--(3/8,1/4);
\draw [-<] (1/2,0)--(7/16,1/8);
\draw (3/8,1/4)--(3/8,1/2);
\node [below] at (0,-.07) {$a$};
\node [below] at (1/4,-.07) {$a$};
\node [below] at (1/2,0) {$b$};
\end{tikzpicture}\\
&=0,
\end{align*} resolving this case.

Since (B4) is the same as ($C_0$) the last overlap we need to check is an overlap between (B2) and ($C_k$) for any $k\geq0.$

\underline{(B2) and ($C_k$):}

There are four cases for such an overlap. Consider first the following two cases:
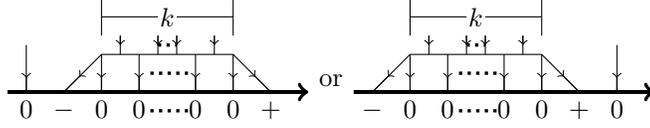
\begin{equation*}
\begin{tikzpicture}[baseline=.5ex]
\draw [line width=1.5, ->] (-1/2,0)--(3.5,0);
\draw (-1/4,0)--(-1/4,5/8);
\draw [-<] (-1/4,0)--(-1/4,1/4);
\draw (1/4,0)--(3/4,1/2);
\draw [-<] (1/4,0)--(1/2,1/4);
\draw (3/4,0)--(3/4,1/2);
\draw [-<] (3/4,0)--(3/4,1/4);
\draw (5/4,0)--(5/4,1/2);
\draw [-<] (5/4,0)--(5/4,1/4);
\draw (2,0)--(2,1/2);
\draw [-<] (2,0)--(2,1/4);
\draw (2+1/2,0)--(2+1/2,1/2);
\draw [-<] (2+1/2,0)--(2+1/2,1/4);
\draw (3,0)--(2+1/2,1/2);
\draw [-<] (3,0)--(2+3/4,1/4);
\draw (3/4,1/2)--(2+1/2,1/2);
\draw (1,1/2)--(1,3/4);
\draw [->] (1,3/4)--(1,5/8);
\draw (1+1/2,1/2)--(1+1/2,3/4);
\draw [->] (1+1/2,3/4)--(1+1/2,5/8);
\draw (1+3/4,1/2)--(1+3/4,3/4);
\draw [->] (1+3/4,3/4)--(1+3/4,5/8);
\draw (2+1/4,1/2)--(2+1/4,3/4);
\draw [->] (2+1/4,3/4)--(2+1/4,5/8);
\node [below] at (-1/4,0) {$0$};
\node [below] at (1/4,0) {$-$};
\node [below] at (3/4,0) {$0$};
\node [below] at (5/4,0) {$0$};
\node [below] at (2,0) {$0$};
\node [below] at (2+1/2,0) {$0$};
\node [below] at (3,0) {$+$};
\draw [line width=1.25, dotted] (1.4,-1/4)--(1.9,-1/4);
\draw [line width=1.25, dotted] (1.4,1/4)--(1.9,1/4);
\draw [line width=1.25, dotted] (1.5,5/8)--(1.75,5/8);
\draw (3/4,3/4)--(3/4,5/4);
\draw (3/4,1)--(6/4,1);
\draw (7/4,1)--(2+1/2,1);
\draw (5/2,3/4)--(5/2,5/4);
\node at (13/8,1) {$k$};
\end{tikzpicture} \text{ or } \begin{tikzpicture}[baseline=.5ex]
\draw [line width=1.5, ->] (0,0)--(4,0);
\draw (1/4,0)--(3/4,1/2);
\draw [-<] (1/4,0)--(1/2,1/4);
\draw (3/4,0)--(3/4,1/2);
\draw [-<] (3/4,0)--(3/4,1/4);
\draw (5/4,0)--(5/4,1/2);
\draw [-<] (5/4,0)--(5/4,1/4);
\draw (2,0)--(2,1/2);
\draw [-<] (2,0)--(2,1/4);
\draw (2+1/2,0)--(2+1/2,1/2);
\draw [-<] (2+1/2,0)--(2+1/2,1/4);
\draw (3,0)--(2+1/2,1/2);
\draw [-<] (3,0)--(2+3/4,1/4);
\draw (3/4,1/2)--(2+1/2,1/2);
\draw (1,1/2)--(1,3/4);
\draw [->] (1,3/4)--(1,5/8);
\draw (1+1/2,1/2)--(1+1/2,3/4);
\draw [->] (1+1/2,3/4)--(1+1/2,5/8);
\draw (1+3/4,1/2)--(1+3/4,3/4);
\draw [->] (1+3/4,3/4)--(1+3/4,5/8);
\draw (2+1/4,1/2)--(2+1/4,3/4);
\draw [->] (2+1/4,3/4)--(2+1/4,5/8);
\draw (3.5,0)--(3.5,5/8);
\draw [-<] (3.5,0)--(3.5,1/4);
\node [below] at (1/4,0) {$-$};
\node [below] at (3/4,0) {$0$};
\node [below] at (5/4,0) {$0$};
\node [below] at (2,0) {$0$};
\node [below] at (2+1/2,0) {$0$};
\node [below] at (3,0) {$+$};
\node [below] at (3.5,0) {$0$};
\draw [line width=1.25, dotted] (1.4,-1/4)--(1.9,-1/4);
\draw [line width=1.25, dotted] (1.4,1/4)--(1.9,1/4);
\draw [line width=1.25, dotted] (1.5,5/8)--(1.75,5/8);
\draw (3/4,3/4)--(3/4,5/4);
\draw (3/4,1)--(6/4,1);
\draw (7/4,1)--(2+1/2,1);
\draw (5/2,3/4)--(5/2,5/4);
\node at (13/8,1) {$k$};
\end{tikzpicture}
\end{equation*}

These two cases are handled symmetrically, so we will focus on the left case. If we apply ($C_k$) first, we obtain
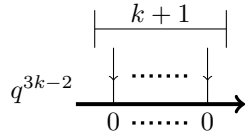
\begin{equation*}
q^{3k-2}\begin{tikzpicture}[baseline=.5ex]
\draw [line width=1.5, ->] (1/2,0)--(2+3/4,0);
\draw (1,0)--(1,3/4);
\draw [-<] (1,0)--(1,3/8);
\draw (2+1/4,0)--(2+1/4,3/4);
\draw [-<] (2+1/4,0)--(2+1/4,3/8);
\node [below] at (1,0) {$0$};
\node [below] at (2+1/4,0) {$0$};
\draw [line width=1.25, dotted] (1+1/4,3/8)--(2,3/8);
\draw [line width=1.25, dotted] (1+1/4,-1/4)--(2,-1/4);
\draw (3/4,3/4)--(3/4,5/4);
\draw (3/4,1)--(2+1/2,1);
\draw (5/2,3/4)--(5/2,5/4);
\node at (13/8,1.2) {$k+1$};
\end{tikzpicture}
\end{equation*}

If we apply (B2) first, we obtain

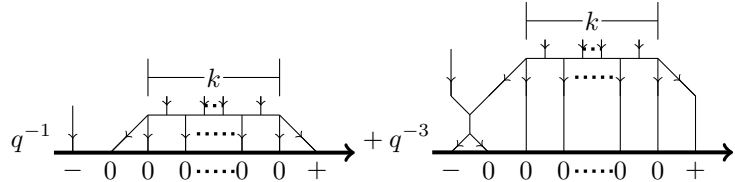
\begin{equation*}
q^{-1}\begin{tikzpicture}[baseline=.5ex]
\draw [line width=1.5, ->] (-1/2,0)--(3.5,0);
\draw (-1/4,0)--(-1/4,5/8);
\draw [-<] (-1/4,0)--(-1/4,1/4);
\draw (1/4,0)--(3/4,1/2);
\draw [-<] (1/4,0)--(1/2,1/4);
\draw (3/4,0)--(3/4,1/2);
\draw [-<] (3/4,0)--(3/4,1/4);
\draw (5/4,0)--(5/4,1/2);
\draw [-<] (5/4,0)--(5/4,1/4);
\draw (2,0)--(2,1/2);
\draw [-<] (2,0)--(2,1/4);
\draw (2+1/2,0)--(2+1/2,1/2);
\draw [-<] (2+1/2,0)--(2+1/2,1/4);
\draw (3,0)--(2+1/2,1/2);
\draw [-<] (3,0)--(2+3/4,1/4);
\draw (3/4,1/2)--(2+1/2,1/2);
\draw (1,1/2)--(1,3/4);
\draw [->] (1,3/4)--(1,5/8);
\draw (1+1/2,1/2)--(1+1/2,3/4);
\draw [->] (1+1/2,3/4)--(1+1/2,5/8);
\draw (1+3/4,1/2)--(1+3/4,3/4);
\draw [->] (1+3/4,3/4)--(1+3/4,5/8);
\draw (2+1/4,1/2)--(2+1/4,3/4);
\draw [->] (2+1/4,3/4)--(2+1/4,5/8);
\node [below] at (-1/4,0) {$-$};
\node [below] at (1/4,0) {$0$};
\node [below] at (3/4,0) {$0$};
\node [below] at (5/4,0) {$0$};
\node [below] at (2,0) {$0$};
\node [below] at (2+1/2,0) {$0$};
\node [below] at (3,0) {$+$};
\draw [line width=1.25, dotted] (1.4,-1/4)--(1.9,-1/4);
\draw [line width=1.25, dotted] (1.4,1/4)--(1.9,1/4);
\draw [line width=1.25, dotted] (1.5,5/8)--(1.75,5/8);
\draw (3/4,3/4)--(3/4,5/4);
\draw (3/4,1)--(6/4,1);
\draw (7/4,1)--(2+1/2,1);
\draw (5/2,3/4)--(5/2,5/4);
\node at (13/8,1) {$k$};
\end{tikzpicture}+q^{-3}\begin{tikzpicture}[baseline=-4.5ex]
\draw [line width=1.5, ->] (-1/2,-3/4)--(3.5,-3/4);
\draw (-1/4,0)--(-1/4,5/8);
\draw [-<] (-1/4,0)--(-1/4,1/4);
\draw (1/4,0)--(3/4,1/2);
\draw [-<] (1/4,0)--(1/2,1/4);
\draw (3/4,0)--(3/4,1/2);
\draw [-<] (3/4,0)--(3/4,1/4);
\draw (5/4,0)--(5/4,1/2);
\draw [-<] (5/4,0)--(5/4,1/4);
\draw (2,0)--(2,1/2);
\draw [-<] (2,0)--(2,1/4);
\draw (2+1/2,0)--(2+1/2,1/2);
\draw [-<] (2+1/2,0)--(2+1/2,1/4);
\draw (3,0)--(2+1/2,1/2);
\draw [-<] (3,0)--(2+3/4,1/4);
\draw (3/4,1/2)--(2+1/2,1/2);
\draw (1,1/2)--(1,3/4);
\draw [->] (1,3/4)--(1,5/8);
\draw (1+1/2,1/2)--(1+1/2,3/4);
\draw [->] (1+1/2,3/4)--(1+1/2,5/8);
\draw (1+3/4,1/2)--(1+3/4,3/4);
\draw [->] (1+3/4,3/4)--(1+3/4,5/8);
\draw (2+1/4,1/2)--(2+1/4,3/4);
\draw [->] (2+1/4,3/4)--(2+1/4,5/8);
\node [below] at (-1/4,-3/4) {$-$};
\node [below] at (1/4,-3/4) {$0$};
\node [below] at (3/4,-3/4) {$0$};
\node [below] at (5/4,-3/4) {$0$};
\node [below] at (2,-3/4) {$0$};
\node [below] at (2+1/2,-3/4) {$0$};
\node [below] at (3,-3/4) {$+$};
\draw [line width=1.25, dotted] (1.4,-1)--(1.9,-1);
\draw [line width=1.25, dotted] (1.4,1/4)--(1.9,1/4);
\draw [line width=1.25, dotted] (1.5,5/8)--(1.75,5/8);
\draw (3/4,3/4)--(3/4,5/4);
\draw (3/4,1)--(6/4,1);
\draw (7/4,1)--(2+1/2,1);
\draw (5/2,3/4)--(5/2,5/4);
\node at (13/8,1) {$k$};
\draw (-1/4,0)--(0,-1/4);
\draw (1/4,0)--(0,-1/4);
\draw (0,-1/4)--(0,-1/2);
\draw (0,-1/2)--(-1/4,-3/4);
\draw [-<] (-1/4,-3/4)--(-1/8,-3/4+1/8);
\draw (0,-1/2)--(1/4,-3/4);
\draw [-<] (1/4,-3/4)--(1/8,-3/4+1/8);
\draw (3/4,0)--(3/4,-3/4);
\draw (5/4,0)--(5/4,-3/4);
\draw (2,0)--(2,-3/4);
\draw (2+1/2,0)--(2+1/2,-3/4);
\draw (3,0)--(3,-3/4);
\end{tikzpicture}
\end{equation*}

The first term in this linear combination becomes zero after applying (B3). The diagram in the second term is isotopic to the diagram appearing on the left side of ($C_{k+1}$). After application of ($C_{k+1}$) we obtain confluence in this case.

The other two possible overlaps between (B2) and ($C_k$) are of the following forms:

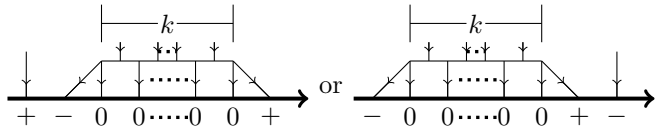
\begin{equation*}
\begin{tikzpicture}[baseline=.5ex]
\draw [line width=1.5, ->] (-1/2,0)--(3.5,0);
\draw (-1/4,0)--(-1/4,5/8);
\draw [-<] (-1/4,0)--(-1/4,1/4);
\draw (1/4,0)--(3/4,1/2);
\draw [-<] (1/4,0)--(1/2,1/4);
\draw (3/4,0)--(3/4,1/2);
\draw [-<] (3/4,0)--(3/4,1/4);
\draw (5/4,0)--(5/4,1/2);
\draw [-<] (5/4,0)--(5/4,1/4);
\draw (2,0)--(2,1/2);
\draw [-<] (2,0)--(2,1/4);
\draw (2+1/2,0)--(2+1/2,1/2);
\draw [-<] (2+1/2,0)--(2+1/2,1/4);
\draw (3,0)--(2+1/2,1/2);
\draw [-<] (3,0)--(2+3/4,1/4);
\draw (3/4,1/2)--(2+1/2,1/2);
\draw (1,1/2)--(1,3/4);
\draw [->] (1,3/4)--(1,5/8);
\draw (1+1/2,1/2)--(1+1/2,3/4);
\draw [->] (1+1/2,3/4)--(1+1/2,5/8);
\draw (1+3/4,1/2)--(1+3/4,3/4);
\draw [->] (1+3/4,3/4)--(1+3/4,5/8);
\draw (2+1/4,1/2)--(2+1/4,3/4);
\draw [->] (2+1/4,3/4)--(2+1/4,5/8);
\node [below] at (-1/4,0) {$+$};
\node [below] at (1/4,0) {$-$};
\node [below] at (3/4,0) {$0$};
\node [below] at (5/4,0) {$0$};
\node [below] at (2,0) {$0$};
\node [below] at (2+1/2,0) {$0$};
\node [below] at (3,0) {$+$};
\draw [line width=1.25, dotted] (1.4,-1/4)--(1.9,-1/4);
\draw [line width=1.25, dotted] (1.4,1/4)--(1.9,1/4);
\draw [line width=1.25, dotted] (1.5,5/8)--(1.75,5/8);
\draw (3/4,3/4)--(3/4,5/4);
\draw (3/4,1)--(6/4,1);
\draw (7/4,1)--(2+1/2,1);
\draw (5/2,3/4)--(5/2,5/4);
\node at (13/8,1) {$k$};
\end{tikzpicture} \text{ or } \begin{tikzpicture}[baseline=.5ex]
\draw [line width=1.5, ->] (0,0)--(4,0);
\draw (1/4,0)--(3/4,1/2);
\draw [-<] (1/4,0)--(1/2,1/4);
\draw (3/4,0)--(3/4,1/2);
\draw [-<] (3/4,0)--(3/4,1/4);
\draw (5/4,0)--(5/4,1/2);
\draw [-<] (5/4,0)--(5/4,1/4);
\draw (2,0)--(2,1/2);
\draw [-<] (2,0)--(2,1/4);
\draw (2+1/2,0)--(2+1/2,1/2);
\draw [-<] (2+1/2,0)--(2+1/2,1/4);
\draw (3,0)--(2+1/2,1/2);
\draw [-<] (3,0)--(2+3/4,1/4);
\draw (3/4,1/2)--(2+1/2,1/2);
\draw (1,1/2)--(1,3/4);
\draw [->] (1,3/4)--(1,5/8);
\draw (1+1/2,1/2)--(1+1/2,3/4);
\draw [->] (1+1/2,3/4)--(1+1/2,5/8);
\draw (1+3/4,1/2)--(1+3/4,3/4);
\draw [->] (1+3/4,3/4)--(1+3/4,5/8);
\draw (2+1/4,1/2)--(2+1/4,3/4);
\draw [->] (2+1/4,3/4)--(2+1/4,5/8);
\draw (3.5,0)--(3.5,5/8);
\draw [-<] (3.5,0)--(3.5,1/4);
\node [below] at (1/4,0) {$-$};
\node [below] at (3/4,0) {$0$};
\node [below] at (5/4,0) {$0$};
\node [below] at (2,0) {$0$};
\node [below] at (2+1/2,0) {$0$};
\node [below] at (3,0) {$+$};
\node [below] at (3.5,0) {$-$};
\draw [line width=1.25, dotted] (1.4,-1/4)--(1.9,-1/4);
\draw [line width=1.25, dotted] (1.4,1/4)--(1.9,1/4);
\draw [line width=1.25, dotted] (1.5,5/8)--(1.75,5/8);
\draw (3/4,3/4)--(3/4,5/4);
\draw (3/4,1)--(6/4,1);
\draw (7/4,1)--(2+1/2,1);
\draw (5/2,3/4)--(5/2,5/4);
\node at (13/8,1) {$k$};
\end{tikzpicture}
\end{equation*}

Since these two cases are handled symmetrically, we will focus on the left case.

We introduce some notation to simplify this computation. We will use symbols placed next to each other to represent certain diagrams appearing next to each other. We represent the diagrams in the left case above by $\downarrow_{+} \cdot C_k.$ We denote by $0_i$ the diagram involving $i$ parallel strands that terminate in good endpoints with states labeled $0.$ We also denote by $X_i$ the diagram 

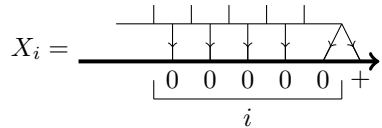
\begin{equation*}
X_i=\begin{tikzpicture}[baseline=.5ex]
\draw [line width=1.5,->] (-1/4,0)--(7/2+1/4,0);
\draw (1,0)--(1,1/2);
\draw [-<] (1,0)--(1,1/4);
\draw (3/2,0)--(3/2,1/2);
\draw [-<] (3/2,0)--(3/2,1/4);
\draw (2,0)--(2,1/2);
\draw [-<] (2,0)--(2,1/4);
\draw (5/2,0)--(5/2,1/2);
\draw [-<] (5/2,0)--(5/2,1/4);
\draw (3,0)--(3+1/4,1/2);
\draw [-<] (3,0)--(3+1/8,1/4);
\draw (7/2,0)--(7/2-1/4,1/2);
\draw [-<] (7/2,0)--(7/2-1/8,1/4);
\draw (1/4,1/2)--(3+1/4,1/2);
\draw (3/4,1/2)--(3/4,3/4);
\draw (5/4,1/2)--(5/4,3/4);
\draw (7/4,1/2)--(7/4,3/4);
\draw (9/4,1/2)--(9/4,3/4);
\draw (11/4,1/2)--(11/4,3/4);
\node [below] at (1,0) {$0$};
\node [below] at (3/2,0) {$0$};
\node [below] at (2,0) {$0$};
\node [below] at (5/2,0) {$0$};
\node [below] at (3,0) {$0$};
\node [below] at (7/2,0) {$+$};
\draw (3/4,-1/2)--(3+1/4,-1/2);
\draw (3/4,-1/2)--(3/4,-1/4);
\draw (3+1/4,-1/2)--(3+1/4,-1/4);
\node [below] at (2,-1/2) {$i$};
\end{tikzpicture}
\end{equation*}

By applying the relation ($C_k$) to $\downarrow _{+} \cdot C_k$ we obtain $q^{3k-2} \downarrow _{+} \cdot 0_k.$ Consider the effect of using relation (B2) on $q^{3k-2} \downarrow_{+} \cdot 0_k$ to get rid of bad pairs, using (B3) at each opportunity. The reduced result is of the form $$q^{3k-2}\sum_{l=0}^k q^{-l}q^{-3(k-l)} 0_l\cdot X_{k-l}=\sum_{l=0}^k q^{2l-2} 0_l \cdot X_{k-l}.$$

We now check that we reach the same reduced result if we instead apply (B2) first to $\downarrow _{+} \cdot C_k$. We introduce another piece of notation. The diagram $A_{i,j}$ has $i$ $0\text{-states}$ on the left of the $+\text{-state}$ and $j$ $ 0 \text{-states}$ on the right.

\begin{equation*}
A_{i,j}=\begin{tikzpicture}[baseline=.5ex]
\draw [line width=1.5,->] (-1/4,0)--(7/2+1/4,0);
\draw (0,0)--(1/4,1/2);
\draw [-<] (0,0)--(1/8,1/4);
\draw (1/2,0)--(1/4,1/2);
\draw [-<] (1/2,0)--(3/8,1/4);
\draw (1,0)--(1,1/2);
\draw [-<] (1,0)--(1,1/4);
\draw (3/2,0)--(3/2,1/2);
\draw [-<] (3/2,0)--(3/2,1/4);
\draw (2,0)--(2,1/2);
\draw [-<] (2,0)--(2,1/4);
\draw (5/2,0)--(5/2,1/2);
\draw [-<] (5/2,0)--(5/2,1/4);
\draw (3,0)--(3+1/4,1/2);
\draw [-<] (3,0)--(3+1/8,1/4);
\draw (7/2,0)--(7/2-1/4,1/2);
\draw [-<] (7/2,0)--(7/2-1/8,1/4);
\draw (1/4,1/2)--(3+1/4,1/2);
\draw (3/4,1/2)--(3/4,3/4);
\draw (5/4,1/2)--(5/4,3/4);
\draw (7/4,1/2)--(7/4,3/4);
\draw (9/4,1/2)--(9/4,3/4);
\draw (11/4,1/2)--(11/4,3/4);
\node [below] at (0,0) {$-$};
\node [below] at (1/2,0) {$0$};
\node [below] at (1,0) {$0$};
\node [below] at (3/2,0) {$0$};
\node [below] at (2,0) {$+$};
\node [below] at (5/2,0) {$0$};
\node [below] at (3,0) {$0$};
\node [below] at (7/2,0) {$+$};
\draw (1/4,-1/2)--(3/2+1/4,-1/2);
\draw (1/4,-1/2)--(1/4,-1/4);
\draw (3/2+1/4,-1/2)--(3/2+1/4,-1/4);
\node [below] at (1,-1/2) {$i$};
\draw (5/2-1/4,-1/2)--(3+1/4,-1/2);
\draw (5/2-1/4,-1/2)--(5/2-1/4,-1/4);
\draw (3+1/4,-1/2)--(3+1/4,-1/4);
\node [below] at (11/4,-1/2) {$j$};
\end{tikzpicture}
\end{equation*}

We also note that diagrams of the following form

\begin{equation*}
\begin{tikzpicture}[baseline=.5ex]
\draw [line width=1.5,->] (-1/4,0)--(7/2+1/4,0);
\draw (0,0)--(1/4,1/2);
\draw [-<] (0,0)--(1/8,1/4);
\draw (1/2,0)--(1/4,1/2);
\draw [-<] (1/2,0)--(3/8,1/4);
\draw (1,0)--(1,1/2);
\draw [-<] (1,0)--(1,1/4);
\draw (3/2,0)--(3/2,1/2);
\draw [-<] (3/2,0)--(3/2,1/4);
\draw (2,0)--(2,1/2);
\draw [-<] (2,0)--(2,1/4);
\draw (5/2,0)--(5/2,1/2);
\draw [-<] (5/2,0)--(5/2,1/4);
\draw (3,0)--(3+1/4,1/2);
\draw [-<] (3,0)--(3+1/8,1/4);
\draw (7/2,0)--(7/2-1/4,1/2);
\draw [-<] (7/2,0)--(7/2-1/8,1/4);
\draw (1/4,1/2)--(3+1/4,1/2);
\draw (3/4,1/2)--(3/4,3/4);
\draw (5/4,1/2)--(5/4,3/4);
\draw (7/4,1/2)--(7/4,3/4);
\draw (9/4,1/2)--(9/4,3/4);
\draw (11/4,1/2)--(11/4,3/4);
\node [below] at (0,0) {$0$};
\node [below] at (1/2,0) {$+$};
\node [below] at (1,0) {$0$};
\node [below] at (3/2,0) {$0$};
\node [below] at (2,0) {$0$};
\node [below] at (5/2,0) {$0$};
\node [below] at (3,0) {$0$};
\node [below] at (7/2,0) {$+$};
\end{tikzpicture}=0
\end{equation*}

are zero, as can be shown by induction on the number of zero states appearing between the two + states. The inductive hypothesis can be applied after applying (B2) once to improve the order of the states and then applying (I2) to remove the square that forms.

If we apply relation (B2) to $\downarrow _{+} \cdot C_k$ then one of the resulting terms will become zero as it is of the form above. We are then left with $$\downarrow _{+} \cdot C_k=q^{-3}A_{0,k+1}.$$

Now consider the diagram $A_{l,m}$ for some $l,m \geq 0.$ We have $A_{l,0}=0$ by relation (B3). For $m>0$ we can apply relation (B2) followed by (I2) and, ignoring the term with the zero diagram as above, we see that $$A_{l,m}=q^{-1}A_{l+1,m-1}+q^3C_l\cdot X_{m-1}.$$

A repeated application of this equation yields \begin{align*}
q^{-3}A_{0,k+1}&= q^{-3}q^3\sum_{i=0}^k q^{-i}C_{i} \cdot X_{k-i}\\
& \stackrel{(C_i)}{=} \sum_{i=0}^k q^{-i}q^{3i-2}0_i \cdot X_{k-i}\\
&= \sum_{i=0}^k q^{2i-2}0_i \cdot X_{k-i}.
\end{align*}

Thus, we have reached local confluence in this last case. The Diamond Lemma now gives us the result.
\end{proof}

We define the $\textit{interior skein algebra}$ $\mathring{\mathcal{S}}_q^{SL_3}(\Sigma)$ as the module freely spanned by closed webs contained in the interior of $\Sigma$ modulo the interior relations (I1a)-(I4b) only.

\begin{corollary}
There is an algebra embedding $$\mathring{\mathcal{S}}_q^{SL_3}(\Sigma) \rightarrow \mathcal{S}_q^{SL_3}(\Sigma)$$ induced by the inclusion map on diagrams.
\end{corollary}

\begin{proof}
Using the reduction rules (I1a)-(I4b) and (S), the Diamond Lemma applies to give a basis for $\mathring{\mathcal{S}}_q^{SL_3}(\Sigma)$. This set of basis diagrams is a subset of basis diagrams of $\mathcal{S}_q^{SL_3}(\Sigma)$, thus the inclusion induces an injective map.
\end{proof}

\section{Bialgebra and comodule structure associated to the bigon}

The surface made by removing one point from the boundary of a closed disk is called the monogon and will be denoted $\mathfrak{M}.$ The surface obtained by removing two points from the boundary of a closed disk is called the bigon and will be denoted $\mathfrak{B}$.
\begin{center}
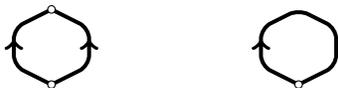

\begin{tikzpicture}
\draw [line width=1.5, ->, rounded corners] (0,0)--(-1/2,1/4)--(-1/2,1/2+.1);
\draw [line width=1.5, rounded corners] (-1/2,1/2)--(-1/2,3/4)--(0,1);
\draw [line width=1.5,->, rounded corners]
(0,0)--(1/2,1/4)--(1/2,1/2+.1);
\draw [line width=1.5, rounded corners]
(1/2,1/2)--(1/2,3/4)--(0,1);
\filldraw [fill=white, draw=black] (0,0) circle [radius=.05];
\filldraw [fill=white, draw=black] (0,1) circle [radius=.05]; 
\end{tikzpicture}
\hspace{2cm}
\begin{tikzpicture}
\draw [line width=1.5, ->, rounded corners] (0,0)--(-1/2,1/4)--(-1/2,1/2+.1);
\draw [line width=1.5, rounded corners] (-1/2,1/2)--(-1/2,3/4)--(0,1)--(1/2,3/4)--(1/2,1/4)--(0,0);
\filldraw [fill=white, draw=black] (0,0) circle [radius=.05];
\end{tikzpicture}
\captionof{figure}{Bigon $\mathfrak{B}$ on left and monogon $\mathfrak{M}$ on right.}
\end{center}

\begin{proposition}
We have that $$\mathcal{S}_q^{SL_3}(\mathfrak{M}) \cong \mathcal{R}.$$
\end{proposition}

\begin{proof}
We show that $\mathcal{S}_q^{SL_3}(\mathfrak{M})$ is spanned by the empty diagram. The fact that the empty diagram is nonzero follows from the fact that it is irreducible and is thus a basis element.

Consider a web diagram $W$ in $\mathcal{S}_q^{SL_3}(\mathfrak{M})$. We can use relations (I1a) and (I1b) to inductively write $W$ as a linear combination of crossingless diagrams. We can use relations (l) or (m) to get rid of vertices near the boundary. If there are strands between a vertex and the boundary we can apply relations (d) or (f) to create room for the vertex to slide over to the boundary without introducing crossings.

So by induction we can write $W$ as a linear combination of diagrams with no crossings and no vertices. After applying relations (I4a) and (I4b) to get rid of circles, these diagrams only have arcs connected to the single boundary arc. By applying relations (g) and (i), these diagrams become scalar multiples of the empty diagram.
\end{proof}

We recall that in \cite{Kup94}, Kuperberg used an Euler characteristic argument to show that the module spanned by closed webs in the plane is 1-dimensional. We remark that by Proposition 3 along with the corollary to the construction of the basis, we obtain an alternate proof that Kuperberg's relations are enough to reduce any closed web in the plane to a scalar multiple of the empty web, and that this reduction can be performed algorithmically  by iteratively applying the left sides of the interior relations. We also observe that Proposition 3 and the algorithm produced by the Diamond Lemma imply that any stated web in $\mathfrak{M}$ can be reduced to a scalar multiple of the empty diagram by iteratively applying just the left sides of the defining relations and ($C_k$).

We next describe the bialgebra structure of $\mathcal{S}_q^{SL_3}(\mathfrak{B}).$ For a counit, we will construct an algebra morphism $\varepsilon: \mathcal{S}_q^{SL_3}(\mathfrak{B}) \rightarrow \mathcal{S}_q^{SL_3}(\mathfrak{M})\cong \mathcal{R}$. As in \cite{Le18} we will use an edge inversion map.

\begin{definition}
If $b$ is a boundary arc of $\Sigma$ with the orientation given in the defining relations of $\mathcal{S}_q^{SL_3}(\Sigma)$ we define the inversion along $b$, $\text{inv}_b: \mathcal{S}_q^{SL_3}(\Sigma) \rightarrow \mathcal{S}_q^{SL_3}(\Sigma)$ to be the $\mathcal{R} \text{-module}$ homomorphism defined on web diagrams by reversing the height order of $b,$ switching the states to their negatives, and multiplying by scalars $C_s^{\uparrow}$ and $C_s^{\downarrow}$ for each endpoint on $b.$ Here, we use $C_s^{\downarrow}=-q^{-4/3}$ for each good endpoint on $b$ with any state $s$ and we use $C_t^{\uparrow}=-q^{-4/3}q^{-2t}$ for each bad endpoint on $b$ with a state $t \in \{-,0,+\}.$
\end{definition}

\begin{proposition}
The map $\text{inv}_b$ defined above is a well-defined $\mathcal{R}\text{-module}$ automorphism.
\end{proposition}

\begin{proof}
We must check that the map respects the defining boundary relations. To do so, we apply the map to both sides of a boundary relation and then reduce the results using the Diamond Lemma algorithm to see that we obtain the same answers in each case. Thus, the map is well-defined. Alternatively, it is easier to use the relations in Section 3 to check that $\text{inv}_b$ respects the relations (c),(e),(h),(j),(k), and (n). We then observe that these relations imply relations (B1)-(B4).  To check that it is an automorphism, one needs to check that the obvious candidate for its inverse is well-defined in the same way. 
\end{proof}

We define $\varepsilon: \mathcal{S}_q^{SL_3}(\mathfrak{B}) \rightarrow \mathcal{S}_q^{SL_3}(\mathfrak{M})$ to be the map given by the result of inverting the the right boundary arc $e_r$ of the bigon with $\text{inv}_{e_r}$ and then filling in the puncture. The map is well-defined since it is a composition of well-defined maps. The fact that it is an algebra morphism is an easy diagrammatic observation, and can be seen in the same way as in \cite{CL19}.

The comultiplication $\Delta: \mathcal{S}_q^{SL_3}(\mathfrak{B}) \rightarrow \mathcal{S}_q^{SL_3}(\mathfrak{B}) \otimes \mathcal{S}_q^{SL_3}(\mathfrak{B})$ is given by the splitting morphism $\Delta_c$ for an ideal arc $c$ that travels from the bottom puncture to the top puncture. By Theorem 1, $\Delta$ is an algebra morphism and satisfies the coassociativity property.

To check that $\varepsilon$ satisfies the counit property, we only need to check on generators. To find a nice set of generators, we use the method in the proof of Proposition 3 to see that any web in the bigon can be written as a linear combination of webs which have no crossings, no vertices, and no circles. Any trivial arcs that start and end on the same boundary arc can be replaced by scalars, and we are left with a linear combination of webs containing only parallel and antiparallel strands with one endpoint on each boundary arc. Thus, $\mathcal{S}_q^{SL_3}(\mathfrak{B})$ has a generating set consisting of diagrams, each of which contain a single strand traveling from one boundary arc of the diagram to the other. We denote such diagrams $\alpha_{st}$ and $\beta_{st}$ depending on the strand orientation and states.

\begin{center}
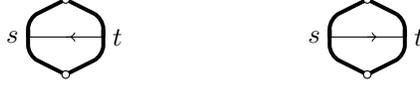

\begin{tikzpicture}
\draw [line width=1.5, rounded corners] (0,0)--(-1/2,1/4)--(-1/2,1/2+.1);
\draw [line width=1.5, rounded corners] (-1/2,1/2)--(-1/2,3/4)--(0,1);
\draw [line width=1.5, rounded corners]
(0,0)--(1/2,1/4)--(1/2,1/2+.1);
\draw [line width=1.5, rounded corners]
(1/2,1/2)--(1/2,3/4)--(0,1);
\filldraw [fill=white, draw=black] (0,0) circle [radius=.05];
\filldraw [fill=white, draw=black] (0,1) circle [radius=.05];
\draw [-<] (-1/2,1/2)--(1/8,1/2);
\draw (0,1/2)--(1/2,1/2);
\node [left] at (-1/2,1/2) {$s$};
\node [right] at (1/2,1/2) {$t$};
\end{tikzpicture}
\hspace{2cm}
\begin{tikzpicture}
\draw [line width=1.5, rounded corners] (0,0)--(-1/2,1/4)--(-1/2,1/2+.1);
\draw [line width=1.5, rounded corners] (-1/2,1/2)--(-1/2,3/4)--(0,1);
\draw [line width=1.5, rounded corners]
(0,0)--(1/2,1/4)--(1/2,1/2+.1);
\draw [line width=1.5, rounded corners]
(1/2,1/2)--(1/2,3/4)--(0,1);
\filldraw [fill=white, draw=black] (0,0) circle [radius=.05];
\filldraw [fill=white, draw=black] (0,1) circle [radius=.05];
\draw [->] (-1/2,1/2)--(1/8,1/2);
\draw (0,1/2)--(1/2,1/2);
\node [left] at (-1/2,1/2) {$s$};
\node [right] at (1/2,1/2) {$t$};
\end{tikzpicture}
\captionof{figure}{Generator $\alpha_{st}$ on left and generator $\beta_{st}$ on right.}
\end{center}

We use our diagrammatic definition of $\varepsilon$ to  compute that \begin{align*}
\varepsilon(\alpha_{st})&= \varepsilon(\begin{tikzpicture}[baseline=3ex]
\draw [line width=1.5, rounded corners] (0,0)--(-1/2,1/4)--(-1/2,1/2+.1);
\draw [line width=1.5, rounded corners] (-1/2,1/2)--(-1/2,3/4)--(0,1);
\draw [line width=1.5, rounded corners]
(0,0)--(1/2,1/4)--(1/2,1/2+.1);
\draw [line width=1.5, rounded corners]
(1/2,1/2)--(1/2,3/4)--(0,1);
\filldraw [fill=white, draw=black] (0,0) circle [radius=.05];
\filldraw [fill=white, draw=black] (0,1) circle [radius=.05];
\draw [-<] (-1/2,1/2)--(1/8,1/2);
\draw (0,1/2)--(1/2,1/2);
\node [left] at (-1/2,1/2) {$s$};
\node [right] at (1/2,1/2) {$t$};
\end{tikzpicture}) \\
&=-q^{-4/3}q^{-2t} ( \begin{tikzpicture}[baseline=3ex]
\draw [line width=1.5, ->, rounded corners] (0,0)--(-1/2,1/4)--(-1/2,3/4)--(0,1);
\draw [line width=1.5, rounded corners] (-1/2,1/2)--(-1/2,3/4)--(0,1)--(1/2,3/4)--(1/2,1/4)--(0,0);
\filldraw [fill=white, draw=black] (0,0) circle [radius=.05];
\draw [-<] (-1/2,1/2)--(1/8,1/2);
\draw (0,1/2)--(1/2,1/2);
\node [left] at (-1/2,1/2) {$s$};
\node [right] at (1/2,1/2) {$-t$};
\end{tikzpicture})\\
&\stackrel{(i)}{=} -q^{-4/3}q^{-2t}(-q^{4/3}q^{2t}\delta_{s-t,0})\\
&= \delta_{st}
\end{align*}

We similarly compute that $\varepsilon(\beta_{st})=\delta_{st}$.

By the definition of $\Delta$, we compute that $$\Delta(\alpha_{st})=\sum_{l \in \{-,0,+\}} \alpha_{sl} \otimes \alpha_{lt}.$$ Similarly, $$\Delta(\beta_{st})=\sum_{l \in \{-,0,+\}} \beta_{sl} \otimes \beta_{lt}.$$

These equations allow us to verify that $$(\varepsilon \otimes \text{id}) \circ \Delta(\alpha_{st})=\alpha_{st}=(\text{id} \otimes \varepsilon) \circ \Delta(\alpha_{st}).$$ The same equations hold for $\beta_{st}$ and we have proven the following proposition.

\begin{proposition}
The algebra $\mathcal{S}_q^{SL_3}(\mathfrak{B})$ has a natural biaglebra structure given by the maps $\Delta$ and $\varepsilon$ defined above.
\end{proposition}

The ingredients here are now the same as in \cite{CL19} and so we obtain an analogue of their Proposition 4.1

\begin{proposition}
Suppose $b$ is a boundary arc of $\Sigma.$ The map defined by splitting $\Sigma$  along an ideal arc isotopic to $b$ so as to split off a bigon $\mathfrak{B}$ whose right edge is $b$ gives an $\mathcal{R}\text{-algebra}$ homomorphism $$\Delta_b: \mathcal{S}_q^{SL_3}(\Sigma) \rightarrow \mathcal{S}_q^{SL_3}(\Sigma)\otimes \mathcal{S}_q^{SL_3}(\mathfrak{B}).$$ This endows $\mathcal{S}_q^{SL_3}(\Sigma)$ with a right comodule-algebra structure over $\mathcal{S}_q^{SL_3}(\mathfrak{B}).$ Similarly, the map $_b\Delta$ defined by splitting off from $\Sigma$ a bigon $\mathfrak{B}$ whose left edge is $b$ gives an $\mathcal{R}\text{-algebra}$ homomorphism $$_b\Delta: \mathcal{S}_q^{SL_3}(\Sigma) \rightarrow  \mathcal{S}_q^{SL_3}(\mathfrak{B})\otimes \mathcal{S}_q^{SL_3}(\Sigma).$$ This endows $\mathcal{S}_q^{SL_3}(\Sigma)$ with a left comodule-algebra structure over $\mathcal{S}_q^{SL_3}(\mathfrak{B}).$
\end{proposition}

\section{Gluing or cutting along a triangle}

Consider a punctured bordered surface $\Sigma$ with two distinct boundary arcs $a$ and $b.$ Also consider an ideal triangle $\mathfrak{T}$, which is a disk with three points removed from its boundary. We will denote the punctured bordered surface $\Sigma \# \mathfrak{T}$ obtained by gluing $\Sigma$ to $\mathfrak{T}$ along $a$ and $b.$ We label the edges of $\mathfrak{T}$ as in the following diagram.

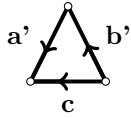
\begin{center}
\begin{tikzpicture}
%bottom triangle edge
\draw[line width=1.5] (-1/2,0)--(1/2,0);
\draw[line width =1.5, -<] (-1/2,0)--(0,0);
%left edge
\draw[line width=1.5] (-1/2,0)--(0,1);
\draw[line width=1.5, -<] (-1/2,0)--(-1/4,1/2);
%right edge
\draw[line width=1.5] (1/2,0)--(0,1);
\draw[line width=1.5,->] (1/2,0)--(1/4,1/2);
%punctures
\filldraw[fill=white, draw=black] (1/2,0) circle [radius=.05];
\filldraw[fill=white, draw=black] (-1/2,0) circle [radius=.05];
\filldraw[fill=white, draw=black] (0,1) circle [radius=.05];
%labels
\node[left] at (-1/4-1/8,1/2+1/8) {\textbf{a'}};
\node[right] at (1/4+1/8,1/2+1/8) {\textbf{b'}};
\node[below] at (0,-1/8) {\textbf{c}};
\end{tikzpicture}
\end{center}

There is a well-defined $\mathcal{R}\text{-module}$ homomorphism: $$\text{glue}_\mathfrak{T}: \mathcal{S}_q^{SL_3}(\Sigma) \rightarrow \mathcal{S}_q^{SL_3}(\Sigma \# \mathfrak{T})$$ defined on diagrams by continuing the strands with endpoints on $a$ or $b$ until they reach $c$. The map is depicted in the following diagram.

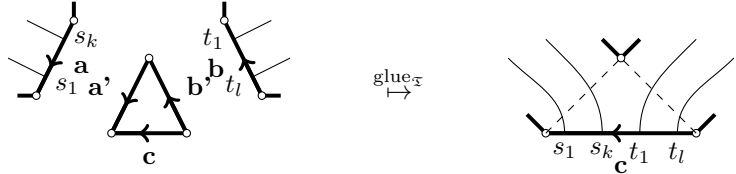
\begin{equation*}
\begin{tikzpicture}[baseline=3ex]
%bottom triangle edge
\draw[line width=1.5] (-1/2,0)--(1/2,0);
\draw[line width =1.5, -<] (-1/2,0)--(0,0);
%left edge
\draw[line width=1.5] (-1/2,0)--(0,1);
\draw[line width=1.5, -<] (-1/2,0)--(-1/4,1/2);
%right edge
\draw[line width=1.5] (1/2,0)--(0,1);
\draw[line width=1.5,->] (1/2,0)--(1/4,1/2);
%punctures
\filldraw[fill=white, draw=black] (1/2,0) circle [radius=.05];
\filldraw[fill=white, draw=black] (-1/2,0) circle [radius=.05];
\filldraw[fill=white, draw=black] (0,1) circle [radius=.05];
%labels
\node[left] at (-1/4-1/8,1/2+1/8) {\textbf{a'}};
\node[right] at (1/4+1/8,1/2+1/8) {\textbf{b'}};
\node[below] at (0,-1/8) {\textbf{c}};
%a
\draw[line width=1.5] (-1.5,.5)--(-1,1.5);
\draw[line width=1.5, -<] (-1.5,.5)--(-1.25,1);
\node[right] at (-1.25+1/8,1-1/8) {\textbf{a}};
\draw[line width=1.5] (-1,1.5)--(-1,1.5+1/4);
\draw[line width=1.5] (-1.5,.5)--(-1.5-1/4,.5);
\filldraw[fill=white, draw=black] (-1.5,.5) circle [radius=.05];
\filldraw[fill=white, draw=black] (-1,1.5) circle [radius=.05];
\draw (-1.5+1/8,.5+1/4)--(-1.5+1/8-1/2,.5+1/4+1/4);
\draw (-1-1/8,1+1/4)--(-1-1/8-1/2,1+1/4+1/4);
\node[right] at (-1.5+1/8,.5+1/8) {$s_1$};
\node[right] at (-1-1/8,1+1/4) {$s_k$};

%b
\draw[line width=1.5] (1.5,.5)--(1,1.5);
\draw[line width=1.5,->] (1.5,.5)--(1.25,1);
\node [left] at (1.25-1/8,1-1/8) {\textbf{b}};
\draw[line width=1.5] (1,1.5)--(1,1.75);
\draw[line width=1.5] (1.5,.5)--(1.75,.5);
\filldraw[fill=white, draw=black] (1.5,.5) circle [radius=.05];
\filldraw[fill=white, draw=black] (1,1.5) circle [radius=.05];
\draw (1.5-1/8,.5+1/4)--(1.5-1/8+1/2,.5+1/4+1/4);
\draw (1+1/8,1+1/4)--(1+1/8+1/2,1+1/4+1/4);
\node[left] at (1.5-1/8,.5+1/8) {$t_l$};
\node[left] at (1+1/8,1+1/4) {$t_1$};
\end{tikzpicture} \hspace{1cm} \stackrel{\text{glue}_\mathfrak{T}}{\mapsto} \hspace{1cm} \begin{tikzpicture}[baseline=3ex]
%bottom edge
\draw[line width=1.5](-1,0)--(1,0);
\draw[line width=1.5,-<](-1,0)--(0,0);
%side edges
\draw[line width=1.5] (-1,0)--(-1-1/4,1/4);
\draw[line width=1.5] (1,0)--(1+1/4,1/4);
\draw[line width=1.5] (0,1)--(-1/4,1+1/4);
\draw[line width=1.5] (0,1)--(1/4,1+1/4);
%punctures
\filldraw[fill=white,draw=black] (-1,0) circle [radius=.05];
\filldraw[fill=white,draw=black] (1,0) circle [radius=.05];
\filldraw[fill=white,draw=black] (0,1) circle [radius=.05];
%strands
\draw (-1.5,1) to[out=-45,in=90] (-3/4,0);
\draw (-1,1.25) to[out=-45,in=90] (-1/4,0);
\draw (1.5,1) to[out=225,in=90] (3/4,0);
\draw (1,1.25) to[out=225,in=90] (1/4,0);
\node [below] at (-3/4,0) {$s_1$};
\node [below] at (-1/4,0) {$s_k$};
\node [below] at (1/4,0) {$t_1$};
\node [below] at (3/4,0) {$t_l$};
\node[below] at (0,-1/4) {\textbf{c}};
%dashed lines
\draw [dashed] (-1,0)--(0,1);
\draw [dashed] (1,0)--(0,1);
\end{tikzpicture}
\end{equation*}

The map $\text{glue}_\mathfrak{T}$ was introduced in \cite{CL19} for the $SL_2$ case. In general, $\text{glue}_\mathfrak{T}$ does not respect the algebra structure, but it gives rise to an algebra structure that is called a \textit{self braided tensor product} in \cite{CL19}. In Section 10 of this paper, we describe a special case of this structure, called the \textit{braided tensor product}. In this section, we are interested in $\text{glue}_\mathfrak{T}$ because it is an $\mathcal{R}\text{-linear}$ isomorphism. We will show this by constructing a natural inverse.

The triangle $\mathfrak{T}$ admits an analogue of the bigon's counit. We define $$\varepsilon_\mathfrak{T}: \mathcal{S}_q^{SL_3}(\mathfrak{T})\rightarrow \mathcal{S}_q^{SL_3}(\mathfrak{M})$$ as the map obtained by applying $\textit{inv}_{b'} \circ \textit{inv}_{a'}$ and then filling in the punctures between $c$ and $a'$ and between $a'$ and $b'$ as in the following figure.

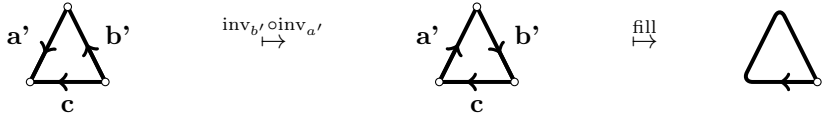
\begin{equation*}
\begin{tikzpicture}[baseline=3ex]
%bottom triangle edge
\draw[line width=1.5] (-1/2,0)--(1/2,0);
\draw[line width =1.5, -<] (-1/2,0)--(0,0);
%left edge
\draw[line width=1.5] (-1/2,0)--(0,1);
\draw[line width=1.5, -<] (-1/2,0)--(-1/4,1/2);
%right edge
\draw[line width=1.5] (1/2,0)--(0,1);
\draw[line width=1.5,->] (1/2,0)--(1/4,1/2);
%punctures
\filldraw[fill=white, draw=black] (1/2,0) circle [radius=.05];
\filldraw[fill=white, draw=black] (-1/2,0) circle [radius=.05];
\filldraw[fill=white, draw=black] (0,1) circle [radius=.05];
%labels
\node[left] at (-1/4-1/8,1/2+1/8) {\textbf{a'}};
\node[right] at (1/4+1/8,1/2+1/8) {\textbf{b'}};
\node[below] at (0,-1/8) {\textbf{c}};
\end{tikzpicture} \hspace{1cm} \stackrel{\text{inv}_{b'} \circ \text{inv}_{a'}}{\mapsto} \hspace{1cm} \begin{tikzpicture}[baseline=3ex]
%bottom triangle edge
\draw[line width=1.5] (-1/2,0)--(1/2,0);
\draw[line width =1.5, -<] (-1/2,0)--(0,0);
%left edge
\draw[line width=1.5] (-1/2,0)--(0,1);
\draw[line width=1.5, ->] (-1/2,0)--(-1/4,1/2);
%right edge
\draw[line width=1.5] (1/2,0)--(0,1);
\draw[line width=1.5,-<] (1/2,0)--(1/4,1/2);
%punctures
\filldraw[fill=white, draw=black] (1/2,0) circle [radius=.05];
\filldraw[fill=white, draw=black] (-1/2,0) circle [radius=.05];
\filldraw[fill=white, draw=black] (0,1) circle [radius=.05];
%labels
\node[left] at (-1/4-1/8,1/2+1/8) {\textbf{a'}};
\node[right] at (1/4+1/8,1/2+1/8) {\textbf{b'}};
\node[below] at (0,-1/8) {\textbf{c}};
\end{tikzpicture} \hspace{1cm} \stackrel{\text{fill}}{\mapsto} \hspace{1cm} \begin{tikzpicture}[baseline=3ex]
\draw [line width=1.5, rounded corners] (1/2,0)--(-1/2,0)--(0,1)--(1/2,0);
\draw[line width=1.5, ->] (1/2,0)--(0,0);
\filldraw[fill=white, draw=black] (1/2,0) circle [radius=.05];
\end{tikzpicture}
\end{equation*}

Since $\varepsilon_{\mathfrak{T}}$ is defined as a composition of well-defined $\mathcal{R}\text{-linear}$ maps it is an $\mathcal{R}\text{-linear}$ map. What makes $\varepsilon_{\mathfrak{T}}$ an analogue of $\varepsilon$ is that if $\varepsilon_{\mathfrak{T}}$ is applied to a diagram $W$ of the following form (with any choice of strand orientations):

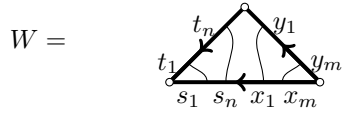
\begin{equation*}
W= \hspace{1cm} \begin{tikzpicture}[baseline=3ex]
%bottom edge
\draw [line width=1.5] (1,0)--(-1,0);
\draw [line width=1.5,->] (1,0)--(-1/8,0);
%left edge
\draw[line width=1.5] (-1,0)--(0,1);
\draw[line width=1.5, -<] (-1,0)--(-1/2,1/2);
%right edge
\draw[line width=1.5] (1,0)--(0,1);
\draw[line width=1.5, ->] (1,0)--(1/2,1/2);
%punctures
\filldraw[fill=white, draw=black] (-1,0) circle [radius=.05];
\filldraw[fill=white, draw=black] (1,0) circle [radius=.05];
\filldraw[fill=white, draw=black] (0,1) circle [radius=.05];
%strands
\draw (-3/4,1/4) to [out=-45,in=90] (-1/2,0);
\draw (-1/4,3/4) to [out=-45,in=90] (-1/4,0);
\draw (1/4,3/4) to [out=-135, in=90] (1/4,0);
\draw (3/4,1/4) to [out=-135, in=90] (1/2,0);
%labels
\node [below] at (-3/4,0) {$s_1$};
\node [below] at (-1/4,0) {$s_n$};
\node [below] at (1/4,0) {$x_1$};
\node [below] at (3/4,0) {$x_m$};
\node [left] at (-3/4,1/4) {$t_1$};
\node [left] at (-1/4,3/4) {$t_n$};
\node [right] at (1/4,3/4) {$y_1$};
\node [right] at (3/4,1/4) {$y_m$};
\end{tikzpicture}
\end{equation*}

the result is $$\varepsilon_{\mathfrak{T}}(W)=(\prod_{i=1}^{n}\delta_{s_i,t_i})(\prod_{j=1}^{m}\delta_{x_j,y_j}).$$

We next define an $\mathcal{R}\text{-linear}$ map $$\text{cut}_{\mathfrak{T}}: \mathcal{S}_q^{SL_3}(\Sigma \# \mathfrak{T}) \rightarrow \mathcal{S}_q^{SL_3}(\Sigma).$$ Recall the notation of the projection $p: \Sigma \sqcup \mathfrak{T} \rightarrow \Sigma \# \mathfrak{T}$ associated to gluing $\Sigma$ to the triangle along $a$ and $b$. If $a''=p(a')=p(a)$ and $b''=p(b')=p(b),$ we define $\text{cut}_\mathfrak{T}$ by $$\text{cut}_\mathfrak{T}=(\varepsilon_{\mathfrak{T}}\otimes \text{id})\circ (\Delta_{b''}\circ \Delta_{a''}).$$

Since $(\Delta_{b''}\circ \Delta_{a''})$ cuts out a triangle, we view it as a linear map $\mathcal{S}_q^{SL_3}(\Sigma \# \mathfrak{T}) \rightarrow S_q^{SL_3}(\mathfrak{T}) \otimes S_q^{SL_3}(\Sigma),$ so the composition above makes sense.

\begin{proposition}
The $\mathcal{R}\text{-linear}$ maps $\text{glue}_\mathfrak{T}$ and $\text{cut}_\mathfrak{T}$ satisfy $$\text{cut}_\mathfrak{T}\circ\text{glue}_\mathfrak{T}=\text{id}_{\mathcal{S}_q^{SL_3}(\Sigma)}$$ and $$\text{glue}_\mathfrak{T} \circ \text{cut}_\mathfrak{T}=\text{id}_{\mathcal{S}_q^{SL_3}(\Sigma \# \mathfrak{T})}.$$
\end{proposition}

\begin{proof}
We will check each equality on a spanning set for the skein algebra involved. For the case of $S_q^{SL_3}(\Sigma)$ we consider the spanning set consisting of all stated web diagrams. Suppose $D$ is a stated web diagram on $\Sigma$. If we examine the diagrams that appear in the triangle cut out by $(\Delta_{b''} \circ \Delta_{a''})\circ \text{glue}_\mathfrak{T}(D)$, we see that they are all of the form $W$ above. Thus, the computation for $\varepsilon_\mathfrak{T}(W)$ above shows that $$(\varepsilon_{\mathfrak{T}}\otimes \text{id})(\Delta_{b''} \circ \Delta_{a''})\circ \text{glue}_\mathfrak{T}(D)=D.$$ This proves the first equality of Proposition 7.

For the second equality, we wish to use a smaller spanning set of $\mathcal{S}_q^{SL_3}(\Sigma \# \mathfrak{T}).$ Consider a stated web diagram $D$ on $(\Sigma \# \mathfrak{T})$ and examine it in a neighborhood of $p(\mathfrak{T})$. By applying an isotopy we can guarantee that $p(\mathfrak{T})$ contains only arcs, and that any arc that enters the triangle through one of the sides either leaves through the other side or terminates at an endpoint on $c.$ After such an isotopy, we obtain a diagram of the following form (for some choice of strand orientations): 

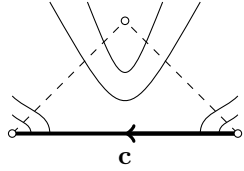
\begin{center}
\begin{tikzpicture}
%c
\draw[line width=1.5] (1.5,0)--(-1.5,0);
\draw[line width=1.5, ->] (1.5,0)--(0,0);
\node [below] at (0,-1/8) {\textbf{c}};
%strands
\draw (-1.5,1/4) to [out=-45, in=90] (-1.25,0);
\draw (-1.5,1/2) to [out=-45, in=90] (-1,0);
\draw (1.5,1/4) to [out=-135, in=90] (1.25,0);
\draw (1.5,1/2) to [out=-135, in=90] (1,0);
\draw (-1,1.75) .. controls (0,0) .. (1,1.75);
\draw (-1/2,1.75) .. controls (0,.5) .. (1/2,1.75);
%dashed lines
\draw [dashed] (-1.5,0)--(0,1.5);
\draw [dashed] (1.5,0)--(0,1.5);
%punctures
\filldraw[fill=white, draw=black] (1.5,0) circle [radius=.05];
\filldraw[fill=white, draw=black] (-1.5,0) circle [radius=.05];
\filldraw[fill=white, draw=black] (0,1.5) circle [radius=.05];
\end{tikzpicture}
\end{center}

Using relations (f) and (j) we can break up the strands that pass through both $a''$ and $b''$ and thus write our diagram $D$ as a linear combination of diagrams of the following form:

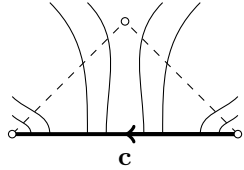
\begin{center}
\begin{tikzpicture}
%c
\draw[line width=1.5] (1.5,0)--(-1.5,0);
\draw[line width=1.5, ->] (1.5,0)--(0,0);
\node [below] at (0,-1/8) {\textbf{c}};
%strands
\draw (-1.5,1/4) to [out=-45, in=90] (-1.25,0);
\draw (-1.5,1/2) to [out=-45, in=90] (-1,0);
\draw (1.5,1/4) to [out=-135, in=90] (1.25,0);
\draw (1.5,1/2) to [out=-135, in=90] (1,0);
\draw (-1,1.75) to [out=-45, in=90] (-1/2,0);
\draw (-1/2,1.75) to [out=-45, in=90] (-1/4,0) ;
\draw (1,1.75) to [out=-135, in=90] (1/2,0);
\draw (1/2,1.75) to [out=-135, in=90] (1/4,0);
%dashed lines
\draw [dashed] (-1.5,0)--(0,1.5);
\draw [dashed] (1.5,0)--(0,1.5);
%punctures
\filldraw[fill=white, draw=black] (1.5,0) circle [radius=.05];
\filldraw[fill=white, draw=black] (-1.5,0) circle [radius=.05];
\filldraw[fill=white, draw=black] (0,1.5) circle [radius=.05];
\end{tikzpicture}
\end{center}

So a spanning set consists of diagrams on $\Sigma \# \mathfrak{T}$ that are of the above form in a neighborhood of $p(\mathfrak{T})$. Let $E$ be such a diagram. We see that the triangles that appear in the terms of $(\Delta_{b''} \circ \Delta_{a''})(E)$ are all of the form $W$ above. Again the computation of $\varepsilon_{\mathfrak{T}}(W)$ above allows us to see that $$\text{glue}_\mathfrak{T}\circ (\varepsilon_{\mathfrak{T}}\otimes \text{id})\circ (\Delta_{b''} \circ \Delta_{a''})(E)=E.$$ This proves the second equality of Proposition 7.

\end{proof}

\begin{corollary}
Suppose $c$ is a boundary arc of a punctured bordered surface $\bar{\Sigma}$ and that $a''$ and $b''$ are ideal arcs with disjoint interiors such that $a''\cup b'' \cup c$ bound an ideal triangle. Then both $\Delta_{a''}$ and $\Delta_{b''}$ are injective.
\end{corollary}

\begin{proof}
Let $\mathfrak{T}$ be the ideal triangle that is split off from $\bar{\Sigma}$ if $\Delta_{b''}\circ \Delta_{a''}$ is applied. Then $\bar{\Sigma}=\Sigma \# \mathfrak{T}$ for the punctured bordered surface $\Sigma$ containing two distinct boundary arcs $a$ and $b$ resulting from the splitting maps. Proposition 7 tells us that $\text{cut}_\mathfrak{T}$ is injective. By the definition of $\text{cut}_\mathfrak{T}$ we see that $\Delta_{b''}\circ \Delta_{a''}$ is injective. Thus, $\Delta_{a''}$ is injective. By Theorem 1, we see that $\Delta_{b''}\circ \Delta_{a''}=\Delta_{a''} \circ \Delta_{b''}$. Thus, $\Delta_{b''}$ is injective as well.
\end{proof}

\section{The triangular decomposition}

We are now able to prove the following addendum to Theorem 1.

\begin{theorem}
Suppose $\bar{\Sigma}$ is a punctured bordered surface and $a''$ is an ideal arc on $\bar{\Sigma}.$ Then the map $\Delta_{a''}$ is injective.
\end{theorem}

\begin{proof}
Let $b''$ be an ideal arc isotopic to $a''$ so that the ideal arcs have disjoint interiors and bound a bigon. Let $c''$ be an ideal arc that bounds a monogon whose ideal vertex is an endpoint of $a''$, and such that $a'',b'',c''$ have disjoint interiors and $a''\cup b'' \cup c''$ bounds an ideal triangle. The following diagram depicts the map $\Delta_{c''}$.

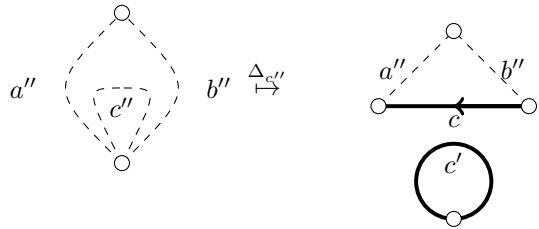
\begin{equation*}
\begin{tikzpicture}[scale=2,baseline=6ex]
%dashed lines
\draw [dashed] (0,0) .. controls (1/4,1/2) .. (0,1/2);
\draw [dashed] (0,0) .. controls (-1/4,1/2) .. (0,1/2);
\draw [dashed] (0,0).. controls (-1/2,1/2)..(0,1);
\draw [dashed] (0,0).. controls (1/2,1/2)..(0,1);
%punctures
\filldraw[fill=white, draw=black] (0,0) circle [radius=.05];
\filldraw[fill=white, draw=black] (0,1) circle [radius=.05];
%labels
\node [below] at (0,1/2) {$c''$};
\node [left] at (-1/2,1/2) {$a''$};
\node [right] at (1/2,1/2) {$b''$};
\end{tikzpicture}  \stackrel{\Delta_{c''}}{\mapsto} \hspace {1cm} \begin{tikzpicture}[scale=2, baseline=1ex]
\draw[line width=1.5] (.5,0)--(-.5,0);
\draw [line width=1.5,->] (.5,0)--(0,0);
% dashed lines
\draw [dashed] (-.5,0)--(0,.5);
\draw [dashed] (.5,0)--(0,.5);
%monogon
\draw [line width=1.5] (0,-.5) circle [radius=.25];
%punctures
\filldraw[fill=white, draw=black] (.5,0) circle [radius=.05];
\filldraw[fill=white, draw=black] (0,.5) circle [radius=.05];
\filldraw[fill=white, draw=black] (-.5,0) circle [radius=.05];
\filldraw[fill=white, draw=black] (0,-.75) circle [radius=.05];
%labels
\node[left] at (-.25,.25) {$a''$};
\node[right] at (.25,.25) {$b''$};
\node[below] at (0,0) {$c$};
\node [below] at (0,-.25) {$c'$};
\end{tikzpicture}
\end{equation*}

Consider the application of $\Delta_{c''}$ to the set of basis diagrams described in Theorem 2. Each irreducible diagram $D$ can be isotoped so that it does not intersect the monogon bounded by $c''.$ This allows us to observe that $\Delta_{c''}(D)$ is an irreducible diagram on its surface as well, and that the isotopy class of $D$ can be completely determined by the isotopy class of this irreducible representative of $\Delta_{c''}(D).$ Thus, $\Delta_{c''}$ maps a basis to a linearly independent set and we conclude that $\Delta_{c''}$ is injective.

After splitting off the monogon bounded by $c''$ we are left with a surface $\Sigma$ that contains a boundary arc $c$ such that $p(c)=c''.$ Now the ideal arcs $a'',b''$ and the boundary arc $c$ satisfy the hypothesis of the previous corollary. By the corollary, $\Delta_{a''}$ is injective on the image of $\Delta_{c''}$ and thus $\Delta_{a''} \circ \Delta_{c''}$ is an injective map. The fact that these maps commute implies that $\Delta_{a''}$ is injective on $\mathcal{S}_q^{SL_3}(\bar{\Sigma})$ as well.
\end{proof}

Now that we have determined the splitting morphisms have trivial kernels, we discuss their images.

Suppose $\Sigma$ is a punctured bordered surface with distinct boundary arcs $a$ and $b.$ Let $\bar{\Sigma}=\Sigma/(a=b)$ and denote by $c$ the common image of $a$ and $b$ under the gluing map. Recall the comodule structure maps associated to the boundary arcs $a,$ and $b.$ We will be interested in $\Delta_a: \mathcal{S}_q^{SL_3}(\Sigma) \rightarrow \mathcal{S}_q^{SL_3}(\Sigma) \otimes \mathcal{S}_q^{SL_3}(\mathfrak{B})$ and $\tau \circ _b\Delta: \mathcal{S}_q^{SL_3}(\Sigma) \rightarrow \mathcal{S}_q^{SL_3}(\Sigma) \otimes \mathcal{S}_q^{SL_3}(\mathfrak{B}),$ where $\tau$ only transposes the tensor factors. We are interested in the following result.

\begin{theorem}
Let $\bar{\Sigma}=\Sigma/(a=b)$ and denote by $c$ the common image of $a$ and $b$ under the gluing map. Then we have $$\text{im}(\Delta_c)=\text{ker}(\Delta_a-\tau\circ_b\Delta).$$
\end{theorem}

\begin{proof}
The inclusion $\text{im}(\Delta_c)\subseteq\text{ker}(\Delta_a-\tau\circ_b\Delta)$ follows by coassociativity of splitting $\bar{\Sigma}$ along $c$ and an ideal arc isotopic to $c.$

To prove the other inclusion, we assume that $y \in \mathcal{S}_q^{SL_3}(\Sigma)$ satisfies $\Delta_a(y)=\tau \circ _b\Delta(y).$ Our goal is to find some $x \in \mathcal{S}_q^{SL_3}(\bar{\Sigma})$ such that $y=\Delta_c(x).$ The element $y$ is represented by a linear combination of stated web diagrams on $\Sigma.$  We will find a candidate for $x$ by trying to weld the strands with endpoints on $a$ or on $b$ to each other. This process uses a map similar to the edge inversion maps $\text{inv}$ before, but this time with a different choice of scalars associated to the endpoints.

For a boundary arc $e$ with positive orientation, we define the edge reversal map $\text{rev}_e$ to be the $\mathcal{R}\text{-linear}$ automorphism of the stated skein module that reverses the height order on $e$, flips the states to their negatives and multiplies by the following scalars for each endpoint on $e$: $^\downarrow _s C= -q^{-4/3}q^{2s}$ for good endpoints with a state $s$ and $^\uparrow_s C=-q^{-4/3}$ for bad endpoints with a state $s$. We can check that this map is well-defined and an automorphism in the same way that we checked this for $\text{inv}_e.$

Let $z=\Delta_a(y)=\tau \circ _b\Delta(y).$ Denote the left boundary arc of the bigon of $\Sigma \sqcup \mathfrak{B}$ by $e_l$ and the right arc by $e_r.$ Let $\mathfrak{T}_1$ and $\mathfrak{T}_2$ be two triangles. We will use the gluing maps $\text{glue}_\mathfrak{T}$ defined in Section 7. Denote the left, right, and bottom edges of the triangles $t_{1l}$, $t_{2l}$, $t_{1r}$, $t_{2r},$ and $t_{1b}$, $t_{2b}$, respectively. We will consider the result of reversing the arc $a$, reversing the arc $e_r,$ then gluing to the triangles. To glue to $\mathfrak{T}_2$ we glue $b$ to $t_{2r}$ and glue $e_r$ to $t_{2l}.$ To glue to $\mathfrak{T}_1,$ we glue $e_l$ to $t_{1r}$ and glue $a$ to $t_{1l}.$

We can write the new element as $\text{glue}_{\mathfrak{T}_1} \circ \text{glue}_{\mathfrak{T}_2} \circ \text{rev}_{e_r} \circ \text {rev}_a(z).$ This gluing is depicted in the following diagram. 

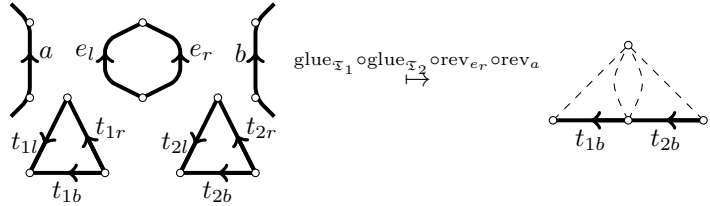
\begin{equation*}
\begin{tikzpicture}[baseline=1ex]
%bigon
\draw [line width=1.5, rounded corners,->] (0,0)--(-1/2,1/4)--(-1/2,1/2+.1);
\draw [line width=1.5, rounded corners] (-1/2,1/2)--(-1/2,3/4)--(0,1);
\draw [line width=1.5, rounded corners,->] (0,0)--(1/2,1/4)--(1/2,1/2+.1);
\draw [line width=1.5, rounded corners] (1/2,1/2)--(1/2,3/4)--(0,1);
\node [left] at (-1/2,1/2+.1) {$e_l$};
\node [right] at (1/2,1/2+.1) {$e_r$};
%bigon punctures
\filldraw [fill=white, draw=black] (0,0) circle [radius=.05];
\filldraw [fill=white, draw=black] (0,1) circle [radius=.05];
%a
\draw [line width=1.5,rounded corners, ->] (-1.75,-1/4)--(-1.5,0)--(-1.5,1/2+.1);
\draw [line width=1.5, rounded corners] (-1.5, 1/2)--(-1.5,1)--(-1.75,5/4);
\filldraw [fill=white, draw=black] (-1.5,0) circle [radius=.05];
\filldraw [fill=white, draw=black] (-1.5,1) circle [radius=.05];
\node [right] at (-1.5,1/2+.1) {$a$};
%b
\draw [line width=1.5,rounded corners, ->] (1.75,-1/4)--(1.5,0)--(1.5,1/2+.1);
\draw [line width=1.5, rounded corners] (1.5, 1/2)--(1.5,1)--(1.75,5/4);
\filldraw [fill=white, draw=black] (1.5,0) circle [radius=.05];
\filldraw [fill=white, draw=black] (1.5,1) circle [radius=.05];
\node [left] at (1.5,1/2+.1) {$b$};
%T1
\draw [line width=1.5] (-.5,-1)--(-1.5,-1);
\draw [line width=1.5,->] (-.5,-1)--(-1,-1);
\draw [line width=1.5] (-1.5,-1)--(-1,0);
\draw [line width=1.5, -<] (-1.5,-1)--(-1.25,-.5);
\draw [line width=1.5] (-.5,-1)--(-1,0);
\draw [line width=1.5, ->] (-.5,-1)--(-.75,-.5);
\filldraw [fill=white, draw=black] (-1.5,-1) circle [radius=.05];
\filldraw [fill=white, draw=black] (-.5,-1) circle [radius=.05];
\filldraw [fill=white, draw=black] (-1,0) circle [radius=.05];
\node [below] at (-1,-1) {$t_{1b}$};
\node [left] at (-1.25,-.6) {$t_{1l}$};
\node [right] at (-.75,-.4) {$t_{1r}$};
%T2
\draw [line width=1.5] (1.5,-1)--(.5,-1);
\draw [line width=1.5,->] (1.5,-1)--(1,-1);
\draw [line width=1.5] (.5,-1)--(1,0);
\draw [line width=1.5, -<] (.5,-1)--(.75,-.5);
\draw [line width=1.5] (1.5,-1)--(1,0);
\draw [line width=1.5, ->] (1.5,-1)--(1.25,-.5);
\filldraw [fill=white, draw=black] (.5,-1) circle [radius=.05];
\filldraw [fill=white, draw=black] (1.5,-1) circle [radius=.05];
\filldraw [fill=white, draw=black] (1,0) circle [radius=.05];
\node [below] at (1,-1) {$t_{2b}$};
\node [left] at (.75,-.6) {$t_{2l}$};
\node [right] at (1.25,-.4) {$t_{2r}$};
\end{tikzpicture} \stackrel{\text{glue}_{\mathfrak{T}_1}\circ \text{glue}_{\mathfrak{T}_2}\circ \text{rev}_{e_r} \circ \text{rev}_{a}}{\mapsto}\begin{tikzpicture}[baseline=3ex]
\draw [line width=1.5] (1,0)--(-1,0);
\draw [line width=1.5, ->] (1,0)--(1/2,0);
\draw [line width=1.5, ->] (0,0)--(-1/2,0);
\node [below] at (-1/2,0) {$t_{1b}$};
\node [below] at (1/2,0) {$t_{2b}$};
%dashed lines
\draw [dashed] (-1,0)--(0,1);
\draw [dashed] (1,0)--(0,1);
\draw [dashed] (0,0).. controls (-1/4,1/2)..(0,1);
\draw [dashed] (0,0).. controls (1/4,1/2)..(0,1);
%punctures
\filldraw[fill=white, draw=black] (-1,0) circle [radius=.05];
\filldraw[fill=white, draw=black] (1,0) circle [radius=.05];
\filldraw[fill=white, draw=black] (0,1) circle [radius=.05];
\filldraw[fill=white, draw=black] (0,0) circle [radius=.05];
\end{tikzpicture}
\end{equation*}

First, we view $z$ as $z=\tau \circ _b\Delta(y)$. Write $y$ as a linear combination of diagrams $D_i$. For each $i,$ $\tau \circ _b\Delta(D_i)$ is a linear combination of diagrams $D_{ij}.$ Each $D_{ij}$ has $k_i$ endpoints on $e_r$, $k_i$ endpoints on $b$, and the states of corresponding endpoints match. After applying $\text{rev}_{e_r}$ to $D_{ij}$ and then gluing to $\mathfrak{T_2},$ we see that there are $2k_i$ endpoints on $t_{2b},$ and that the endpoints which are $k_i\text{-th}$ and $k_i+1\text{-st}$ in the height order have opposite states and opposite orientations. The scalars associated with the application of $\text{rev}_{e_r}$ guarantee that relations (d) or (f) are applicable and allow us to reduce the number of endpoints on $t_{2b}.$ After applying these relations $k_i$ times for each $D_i$, we see that we can write $\text{glue}_{\mathfrak{T}_1} \circ \text{glue}_{\mathfrak{T}_2} \circ \text{rev}_{e_r} \circ \text {rev}_a(z)$ as a linear combination of diagrams, where no diagram has an endpoint on $t_{2b}.$ As no reduction rule from our Diamond Lemma algorithm can result in an endpoint appearing on a boundary arc that previously contained no endpoints, we see that when we write $\text{glue}_{\mathfrak{T}_1} \circ \text{glue}_{\mathfrak{T}_2} \circ \text{rev}_{e_r} \circ \text {rev}_a(z)$ as a linear combination of our basis diagrams, each basis diagram that appears in the linear combination has no endpoints on $t_{2b}.$

Next, we view $z$ as $z=\Delta_a(y).$ In a similar way as the last paragraph, we see that after applying $\text{rev}_a$ and gluing to $\mathfrak{T_1},$ we can apply relations (d) and (f) to write $\text{glue}_{\mathfrak{T}_1} \circ \text{glue}_{\mathfrak{T}_2} \circ \text{rev}_{e_r} \circ \text {rev}_a(z)$ as a linear combination of basis diagrams such that no diagram has an endpoint on $t_{1b}.$ By the uniqueness of this linear combination we see that we can write it as a linear combination of basis diagrams so that no diagram appearing in the linear combination has an endpoint on $t_{1b}$ or on $t_{2b}.$ 

Now, $\text{glue}_{\mathfrak{T}_1} \circ \text{glue}_{\mathfrak{T}_2} \circ \text{rev}_{e_r} \circ \text {rev}_a(z)$ is a linear combination of basis diagrams on $(\Sigma \sqcup \mathfrak{B})\#\mathfrak{T}_1\#\mathfrak{T}_2.$ Consider the surface $(\Sigma \sqcup \mathfrak{B})\#\mathfrak{T}_1\#\mathfrak{T}_2\setminus(t_{1b}\cup t_{2b})$. This is not a punctured bordered surface, but depending on whether the appropriate endpoint of $c$ was a boundary puncture or was an interior puncture, this surface is either a punctured bordered surface missing an interval on its boundary or it is a punctured bordered surface missing a boundary circle. In either case, it is naturally diffeomorphic to the original punctured bordered surface $\bar{\Sigma}$ by replacing this missing boundary interval or boundary circle with a single puncture. There is a linear map defined on the submodule of $\mathcal{S}_q^{SL_3}((\Sigma \sqcup \mathfrak{B})\#\mathfrak{T}_1\#\mathfrak{T}_2)$ spanned by basis diagrams that have no endpoints on $t_{1b}$ or $t_{2b}$ that takes such a basis diagram and embeds it in $(\Sigma \sqcup \mathfrak{B})\#\mathfrak{T}_1\#\mathfrak{T}_2\setminus(t_{1b}\cup t_{2b})$. After applying this map to $\text{glue}_{\mathfrak{T}_1} \circ \text{glue}_{\mathfrak{T}_2} \circ \text{rev}_{e_r} \circ \text {rev}_a(z)$ and composing with our diffeomorphism, we obtain our candidate $x \in \mathcal{S}_q^{SL_3}(\bar{\Sigma}).$

To see that $x$ is the correct choice, we consider $\Delta_c(x)$ and then apply the same process to it as we did to $y$ and observe that $$\text{glue}_{\mathfrak{T}_1} \circ \text{glue}_{\mathfrak{T}_2} \circ \text{rev}_{e_r} \circ \text {rev}_a\circ \Delta_a(y)=\text{glue}_{\mathfrak{T}_1} \circ \text{glue}_{\mathfrak{T}_2} \circ \text{rev}_{e_r} \circ \text {rev}_a \circ\Delta_a (\Delta_c (x)).$$ The injectivity of the maps involved allow us to conclude that $\Delta_c(x)=y.$
\end{proof}

We say that a punctured bordered surface is \textit{ideal triangulable} if it can be obtained from a finite collection of disjoint triangles by gluing some pairs of edges together. It is known that a punctured bordered surface is ideal triangulable if it has no connected component that is one of the following: a closed surface, a sphere with fewer than three punctures, a bigon, or a monogon.

If $\Sigma$ is an ideal triangulable punctured bordered surface, then the images of the glued edges are ideal arcs on $\Sigma$ with disjoint interiors. These form the set of interior edges $\mathcal{E}$ for the $\textit{ideal triangulation}$ of $\Sigma$. Let $p: \displaystyle\sqcup_{i=1}^n \mathfrak{T}_i \rightarrow \Sigma$ be the gluing map. If $e\in \mathcal{E}$, then its preimage $p^{-1}(e)=\{e',e''\}$ consists of two triangle edges. The composition $\Delta$ of the splitting maps $\Delta_e$ for $e \in \mathcal{E}$ gives an algebra embedding $$\Delta: \mathcal{S}_q^{SL_3}(\Sigma) \rightarrow \bigotimes_{i=1}^n \mathcal{S}_q^{SL_3}(\mathfrak{T}_i).$$ The composition $^L\Delta$ of all left comodule maps $_{e''}\Delta$ gives a map $$^L\Delta: \bigotimes_{i=1}^n \mathcal{S}_q^{SL_3}(\mathfrak{T}_i)\rightarrow (\bigotimes_{e \in \mathcal{E}} \mathcal{S}_q^{SL_3}(\mathfrak{B}))\otimes (\bigotimes_{i=1}^n \mathcal{S}_q^{SL_3}(\mathfrak{T}_i)).$$ The composition $\Delta^R$ of all right comodule maps $\Delta_{e'}$ gives a map $$\Delta^R:\bigotimes_{i=1}^n \mathcal{S}_q^{SL_3}(\mathfrak{T}_i)\rightarrow (\bigotimes_{i=1}^n \mathcal{S}_q^{SL_3}(\mathfrak{T}_i))\otimes(\bigotimes_{e \in \mathcal{E}} \mathcal{S}_q^{SL_3}(\mathfrak{B})).$$ Then Theorem 3 and Theorem 4 allow us to observe the following corollary.

\begin{corollary}
	If $\Sigma$ admits an ideal triangulation with a set of interior edges $\mathcal{E},$ then the following sequence is exact:
	
	$$0\rightarrow \mathcal{S}_q^{SL_3}(\Sigma) \stackrel{\Delta}{\rightarrow} \bigotimes_{i=1}^n \mathcal{S}_q^{SL_3}(\mathfrak{T}_i) \stackrel{\Delta^R-\tau \circ ^L\Delta}{\displaystyle\rightarrow} (\bigotimes_{i=1}^n \mathcal{S}_q^{SL_3}(\mathfrak{T}_i))\otimes(\bigotimes_{e \in \mathcal{E}} \mathcal{S}_q^{SL_3}(\mathfrak{B})).$$
\end{corollary}

\section{The stated skein algebra of the bigon}

In \cite{CL19}, it was shown that the Kauffman bracket stated skein algebra of the bigon is isomorphic to $\mathcal{O}_{q}(SL_2)$ as a Hopf algebra (with a suitable renormalization of $q$). They showed this by defining a bialgebra map between $\mathcal{O}_{q}(SL_2)$ and the Kauffman bracket stated skein algebra of the bigon. The fact that this map is an isomorphism follows because it maps the canonical basis of the stated skein algebra to a well known basis of $\mathcal{O}_{q}(SL_2).$ There is an analogous isomorphism between our $SL_3$ stated skein algebra of the bigon and $\mathcal{O}_q(SL_3).$ However, the proof here will require us to define maps in both directions since it is not otherwise clear that the canonical basis of the $SL_3$ stated skein algebra of the bigon matches up with a basis of $\mathcal{O}_q(SL_3).$

We first recall the $R\text{-matrix}$ definition of $\mathcal{O}_q(SL_3).$ Consider the free $\mathcal{R}\text{-module}$ $V$ with basis $\{x_1, x_2, x_3\}$. The standard $R\text{-matrix}$ for $SL_3$ is a linear map $$R: V \otimes V \rightarrow V \otimes V$$ defined by $$R(x_i \otimes x_j) =q^{-1/3} \begin{cases}
q x_i \otimes x_j & \text{ (if $i=j$)}\\
x_j \otimes x_i & \text{ (if $i<j$)}\\
x_j \otimes x_i + (q-q^{-1})x_i \otimes x_j & \text{ (if $i>j$)}.\\
\end{cases}$$

We develop some notation for the matrix entries $R_{ij}^{kl}$ of $R.$ We have that $R(x_i \otimes x_j)$ is uniquely written as $$R(x_i \otimes x_j) = \sum_{1 \leq k,l \leq 3} R_{ij}^{kl} x_k \otimes x_l.$$

We define $\mathcal{O}_q(SL_3)$ as the free $\mathcal{R}\text{-algebra}$ generated by elements $\{X_{ij}\}_{1\leq i,j \leq 3}$ modulo the following relations $$\begin{cases}
\displaystyle\sum_{1 \leq k,l \leq 3} R_{ij}^{kl}X_{km}X_{ln}=\displaystyle\sum_{1 \leq k,l \leq 3} R_{kl}^{mn} X_{ik}X_{jl}& \text{(for $1 \leq i,j,m,n \leq 3$)}\\
\displaystyle\sum_{\sigma \in S_3} (-q)^{l(\sigma)} X_{\sigma_11}X_{\sigma_22}X_{\sigma_33}=1.\\
\end{cases}$$

Here, we consider $(\sigma_1,\sigma_2,\sigma_3)=(1,2,3)$ the identity permutation.

The left side of the second equation is called the \textit{quantum determinant}, $\det_q,$ of the matrix of generators $(A)_{ij}=X_{ij}$. We will also make use of notation $A[i\vert j]$ to mean the quantum minor of $A$ after deleting row $i$ and column $j.$  

$\mathcal{O}_q(SL_3)$ has a Hopf algebra structure with structure maps given by $$\varepsilon(X_{ij})=\delta_{ij}$$ and $$\Delta(X_{ij})=\sum_{k=1}^3 X_{ik} \otimes X_{kj}.$$

The antipode $S: \mathcal{O}_q(SL_3)\rightarrow \mathcal{O}_q(SL_3)$ is defined by $$S(X_{ij})=(-q)^{i-j} A[j|i].$$

For the purpose of notation to match up our stated skein algebra with the standard definition of $\mathcal{O}_q(SL_3),$ we define a bijection $t: \{1,2,3\} \rightarrow \{-,0,+\}$ given by $t(1)=+$, $t(2)=0$, $t(3)=-.$ Since $t$ reverses the order we've placed on the sets $\{1,2,3\}$ and $\{-,0,+\}$ we will have to take care when we apply relations (k)-(n) to diagrams.

\begin{proposition}
There is a unique bialgebra morphism $\phi: \mathcal{O}_q(SL_3) \rightarrow \mathcal{S}_q^{SL_3}(\mathfrak{B})$ defined by $$\phi(X_{ij})=%\includegraphics[width=2cm]{bigon generator}
\begin{tikzpicture}[baseline=3ex]
\draw [line width=1.5, rounded corners] (0,0)--(-1/2,1/4)--(-1/2,1/2+.1);
\draw [line width=1.5, rounded corners] (-1/2,1/2)--(-1/2,3/4)--(0,1);
\draw [line width=1.5, rounded corners]
(0,0)--(1/2,1/4)--(1/2,1/2+.1);
\draw [line width=1.5, rounded corners]
(1/2,1/2)--(1/2,3/4)--(0,1);
\filldraw [fill=white, draw=black] (0,0) circle [radius=.05];
\filldraw [fill=white, draw=black] (0,1) circle [radius=.05];
\draw [->] (-1/2,1/2)--(1/8,1/2);
\draw (0,1/2)--(1/2,1/2);
\node [left] at (-1/2,1/2) {$t(i)$};
\node [right] at (1/2,1/2) {$t(j)$};
\end{tikzpicture}
$$
\end{proposition}

\begin{proof}
Since the elements $X_{ij}$ generate $\mathcal{O}_q(SL_3)$, the morphism will be unique if it exists. By construction, such a morphism will preserve the bialgebra structure. To prove that $\phi$ gives a well-defined algebra morphism we must check that it respects the defining relations of $\mathcal{O}_q(SL_3)$. We must show that the relations $$\displaystyle\sum_{1 \leq k,l \leq 3} R_{ij}^{kl}\phi(X_{km})\phi(X_{ln})=\displaystyle\sum_{1 \leq k,l \leq 3} R_{kl}^{mn} \phi(X_{ik})\phi(X_{jl})$$ and $$\displaystyle\sum_{\sigma \in S_3} (-q)^{l(\sigma)} \phi(X_{\sigma_1 1})\phi(X_{\sigma_22})\phi(X_{\sigma_33})=1$$ hold in $\mathcal{S}_q^{SL_3}(\mathfrak{B}).$ For this, we recall the bialgebra structure of the bigon given in Section 5. We consider the result of applying $(\varepsilon \otimes \text{id})\circ \Delta$ to the following diagram in two different ways.

\begin{center}
\begin{tikzpicture}[use Hobby shortcut, baseline=3ex]
%bigon
\draw [line width=1.5, rounded corners,->] (0,0)--(-1/2,1/4)--(-1/2,1/2+.1);
\draw [line width=1.5, rounded corners] (-1/2,1/2)--(-1/2,3/4)--(0,1);
\draw [line width=1.5, rounded corners,->] (0,0)--(1/2,1/4)--(1/2,1/2+.1);
\draw [line width=1.5, rounded corners] (1/2,1/2)--(1/2,3/4)--(0,1);
%punctures
\filldraw [fill=white, draw=black] (0,0) circle [radius=.05];
\filldraw [fill=white, draw=black] (0,1) circle [radius=.05];
%strands
\begin{knot}[
consider self intersections=true,
%  draft mode=crossings,
%flip crossing=1,
ignore endpoint intersections=false]
\strand  (-1/2,3/4)..(1/2,1/4);
\strand  (-1/2,1/4)..(1/2,3/4);
\end{knot}
\draw [-<] (1/2,1/4)--(1/4,1/4+1/8);
\draw [-<] (1/2,3/4)--(1/4,3/4-1/8);
%labels
\node [left] at (-1/2,3/4) {$t(i)$};
\node [left] at (-1/2,1/4) {$t(j)$};
\node [right] at (1/2,3/4) {$t(m)$};
\node [right] at (1/2,1/4) {$t(n)$};
\end{tikzpicture}
\end{center}

For the first way, we split the bigon along an ideal arc that stays to the right of the crossing and obtain

\begin{equation*}
\displaystyle\sum_{1 \leq k,l \leq 3} \varepsilon\bigg(\begin{tikzpicture}[use Hobby shortcut, baseline=3ex]
%bigon
\draw [line width=1.5, rounded corners,->] (0,0)--(-1/2,1/4)--(-1/2,1/2+.1);
\draw [line width=1.5, rounded corners] (-1/2,1/2)--(-1/2,3/4)--(0,1);
\draw [line width=1.5, rounded corners,->] (0,0)--(1/2,1/4)--(1/2,1/2+.1);
\draw [line width=1.5, rounded corners] (1/2,1/2)--(1/2,3/4)--(0,1);
%punctures
\filldraw [fill=white, draw=black] (0,0) circle [radius=.05];
\filldraw [fill=white, draw=black] (0,1) circle [radius=.05];
%strands
\begin{knot}[
consider self intersections=true,
%  draft mode=crossings,
%flip crossing=1,
ignore endpoint intersections=false]
\strand  (-1/2,3/4)..(1/2,1/4);
\strand  (-1/2,1/4)..(1/2,3/4);
\end{knot}
\draw [-<] (1/2,1/4)--(1/4,1/4+1/8);
\draw [-<] (1/2,3/4)--(1/4,3/4-1/8);
%labels
\node [left] at (-1/2,3/4) {$t(i)$};
\node [left] at (-1/2,1/4) {$t(j)$};
\node [right] at (1/2,3/4) {$t(k)$};
\node [right] at (1/2,1/4) {$t(l)$};
\end{tikzpicture}\bigg) \begin{tikzpicture}[baseline=3ex]
%bigon
\draw [line width=1.5, rounded corners,->] (0,0)--(-1/2,1/4)--(-1/2,1/2+.1);
\draw [line width=1.5, rounded corners] (-1/2,1/2)--(-1/2,3/4)--(0,1);
\draw [line width=1.5, rounded corners,->] (0,0)--(1/2,1/4)--(1/2,1/2+.1);
\draw [line width=1.5, rounded corners] (1/2,1/2)--(1/2,3/4)--(0,1);
%punctures
\filldraw [fill=white, draw=black] (0,0) circle [radius=.05];
\filldraw [fill=white, draw=black] (0,1) circle [radius=.05];
%strands
\draw (-1/2,3/4)--(1/2,3/4);
\draw [->] (-1/2,3/4)--(0,3/4);
\draw (-1/2,1/4)--(1/2,1/4);
\draw [->] (-1/2,1/4)--(0,1/4);
%labels
\node [left] at (-1/2,3/4) {$t(k)$};
\node [left] at (-1/2,1/4) {$t(l)$};
\node [right] at (1/2,3/4) {$t(m)$};
\node [right] at (1/2,1/4) {$t(n)$};
\end{tikzpicture}
\end{equation*}

For the second way, we split the bigon along an ideal arc that stays to the left of the crossing and then apply $\text{id} \otimes \varepsilon.$

\begin{equation*}
\displaystyle\sum_{1 \leq k,l \leq 3}  \begin{tikzpicture}[baseline=3ex]
%bigon
\draw [line width=1.5, rounded corners,->] (0,0)--(-1/2,1/4)--(-1/2,1/2+.1);
\draw [line width=1.5, rounded corners] (-1/2,1/2)--(-1/2,3/4)--(0,1);
\draw [line width=1.5, rounded corners,->] (0,0)--(1/2,1/4)--(1/2,1/2+.1);
\draw [line width=1.5, rounded corners] (1/2,1/2)--(1/2,3/4)--(0,1);
%punctures
\filldraw [fill=white, draw=black] (0,0) circle [radius=.05];
\filldraw [fill=white, draw=black] (0,1) circle [radius=.05];
%strands
\draw (-1/2,3/4)--(1/2,3/4);
\draw [->] (-1/2,3/4)--(0,3/4);
\draw (-1/2,1/4)--(1/2,1/4);
\draw [->] (-1/2,1/4)--(0,1/4);
%labels
\node [left] at (-1/2,3/4) {$t(i)$};
\node [left] at (-1/2,1/4) {$t(j)$};
\node [right] at (1/2,3/4) {$t(k)$};
\node [right] at (1/2,1/4) {$t(l)$};
\end{tikzpicture} \varepsilon\bigg(\begin{tikzpicture}[use Hobby shortcut, baseline=3ex]
%bigon
\draw [line width=1.5, rounded corners,->] (0,0)--(-1/2,1/4)--(-1/2,1/2+.1);
\draw [line width=1.5, rounded corners] (-1/2,1/2)--(-1/2,3/4)--(0,1);
\draw [line width=1.5, rounded corners,->] (0,0)--(1/2,1/4)--(1/2,1/2+.1);
\draw [line width=1.5, rounded corners] (1/2,1/2)--(1/2,3/4)--(0,1);
%punctures
\filldraw [fill=white, draw=black] (0,0) circle [radius=.05];
\filldraw [fill=white, draw=black] (0,1) circle [radius=.05];
%strands
\begin{knot}[
consider self intersections=true,
%  draft mode=crossings,
%flip crossing=1,
ignore endpoint intersections=false]
\strand  (-1/2,3/4)..(1/2,1/4);
\strand  (-1/2,1/4)..(1/2,3/4);
\end{knot}
\draw [-<] (1/2,1/4)--(1/4,1/4+1/8);
\draw [-<] (1/2,3/4)--(1/4,3/4-1/8);
%labels
\node [left] at (-1/2,3/4) {$t(k)$};
\node [left] at (-1/2,1/4) {$t(l)$};
\node [right] at (1/2,3/4) {$t(m)$};
\node [right] at (1/2,1/4) {$t(n)$};
\end{tikzpicture}\bigg)
\end{equation*}

The bialgebra axiom $(\varepsilon\otimes \text{id})\Delta=(\text{id} \otimes \varepsilon)\Delta$ along with the isotopy invariance of the splitting map guarantees that both answers must be the same.

We can use the defining relations to compute that

\begin{equation*}
\varepsilon\bigg(\begin{tikzpicture}[use Hobby shortcut, baseline=3ex]
%bigon
\draw [line width=1.5, rounded corners,->] (0,0)--(-1/2,1/4)--(-1/2,1/2+.1);
\draw [line width=1.5, rounded corners] (-1/2,1/2)--(-1/2,3/4)--(0,1);
\draw [line width=1.5, rounded corners,->] (0,0)--(1/2,1/4)--(1/2,1/2+.1);
\draw [line width=1.5, rounded corners] (1/2,1/2)--(1/2,3/4)--(0,1);
%punctures
\filldraw [fill=white, draw=black] (0,0) circle [radius=.05];
\filldraw [fill=white, draw=black] (0,1) circle [radius=.05];
%strands
\begin{knot}[
consider self intersections=true,
%  draft mode=crossings,
%flip crossing=1,
ignore endpoint intersections=false]
\strand  (-1/2,3/4)..(1/2,1/4);
\strand  (-1/2,1/4)..(1/2,3/4);
\end{knot}
\draw [-<] (1/2,1/4)--(1/4,1/4+1/8);
\draw [-<] (1/2,3/4)--(1/4,3/4-1/8);
%labels
\node [left] at (-1/2,3/4) {$t(a)$};
\node [left] at (-1/2,1/4) {$t(b)$};
\node [right] at (1/2,3/4) {$t(c)$};
\node [right] at (1/2,1/4) {$t(d)$};
\end{tikzpicture}\bigg)=R_{ab}^{cd}
\end{equation*}

Equating our two answers shows that the relations $$\displaystyle\sum_{1 \leq k,l \leq 3} R_{ij}^{kl}\phi(X_{km})\phi(X_{ln})=\displaystyle\sum_{1 \leq k,l \leq 3} R_{kl}^{mn} \phi(X_{ik})\phi(X_{jl})$$ hold in $\mathcal{S}_q^{SL_3}(\mathfrak{B}).$

Next, we consider the following diagram

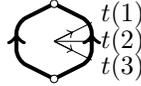
\begin{center}
\begin{tikzpicture}
%bigon
\draw [line width=1.5, rounded corners,->] (0,0)--(-1/2,1/4)--(-1/2,1/2+.1);
\draw [line width=1.5, rounded corners] (-1/2,1/2)--(-1/2,3/4)--(0,1);
\draw [line width=1.5, rounded corners,->] (0,0)--(1/2,1/4)--(1/2,1/2+.1);
\draw [line width=1.5, rounded corners] (1/2,1/2)--(1/2,3/4)--(0,1);
%strands
\draw (0,1/2)--(1/2,3/4);
\draw [->] (0,1/2)--(1/4,5/8);
\draw (0,1/2)--(1/2,1/2);
\draw (0,1/2)--(1/2,1/4);
\draw [->] (0,1/2)--(1/4,3/8);
%punctures
\filldraw [fill=white, draw=black] (0,0) circle [radius=.05];
\filldraw [fill=white, draw=black] (0,1) circle [radius=.05];
%labels
\node [right] at (1/2,3/4+.1) {$t(1)$};
\node [right] at (1/2,1/2) {$t(2)$};
\node [right] at (1/2,1/4-.1) {$t(3)$};
\end{tikzpicture}
\end{center}

On one hand, we can evaluate this diagram using relation (k) from Section 2 along the right edge of the bigon. On the other hand, we could use relation (l) along the left edge of the bigon.

This gives us the relation

\begin{equation*}
q^{-2}=q^{-2} \displaystyle \sum_{\sigma \in S_3} (-q)^{l(\sigma)} \begin{tikzpicture}[baseline=3ex]
%bigon
\draw [line width=1.5, rounded corners,->] (0,0)--(-1/2,1/4)--(-1/2,1/2+.1);
\draw [line width=1.5, rounded corners] (-1/2,1/2)--(-1/2,3/4)--(0,1);
\draw [line width=1.5, rounded corners,->] (0,0)--(1/2,1/4)--(1/2,1/2+.1);
\draw [line width=1.5, rounded corners] (1/2,1/2)--(1/2,3/4)--(0,1);
%strands
\draw (-1/2,3/4)--(1/2,3/4);
\draw [->] (-1/2,3/4)--(0,3/4);
\draw (-1/2,1/2)--(1/2,1/2);
\draw [->] (-1/2,1/2)--(0,1/2);
\draw (-1/2,1/4)--(1/2,1/4);
\draw [->] (-1/2,1/4)--(0,1/4);
%punctures
\filldraw [fill=white, draw=black] (0,0) circle [radius=.05];
\filldraw [fill=white, draw=black] (0,1) circle [radius=.05];
%labels
\node [right] at (1/2,3/4+.1) {$t(1)$};
\node [right] at (1/2,1/2) {$t(2)$};
\node [right] at (1/2,1/4-.1) {$t(3)$};
\node [left] at (-1/2,3/4+.1) {$t(\sigma_1)$};
\node [left] at (-1/2,1/2) {$t(\sigma_2)$};
\node [left] at (-1/2,1/4-.1) {$t(\sigma_3)$};
\end{tikzpicture}
\end{equation*}

Thus, the relation $$\displaystyle\sum_{\sigma \in S_3} (-q)^{l(\sigma)} \phi(X_{\sigma_1 1})\phi(X_{\sigma_22})\phi(X_{\sigma_33})=1$$ holds in $\mathcal{S}_q^{SL_3}(\mathfrak{B})$. Thus, $\phi$ is well-defined.

\end{proof}

To prove that $\phi$ is an isomorphism, we will construct an inverse function. We will define an algebra morphism $\psi: S_q^{SL_3}(\mathfrak{B}) \rightarrow \mathcal{O}_q(SL_3)$ by defining it on diagrams and then checking that it is well-defined.

In order for $\psi$ to be the inverse of $\phi$ we are forced to define it on the diagrams $\alpha_{t(i)t(j)}$ and $\beta_{t(i)t(j)}$ from Section 6 as $$\psi(\beta_{t(i)t(j)})=X_{ij}$$ and $$\psi(\alpha_{t(i)t(j)})= (-q)^{j-i} A[4-i \vert 4-j]. $$

As was noted in Section 6, the diagrams $\alpha_{t(i)t(j)}$ and $\beta_{t(i)t(j)}$ generate $\mathcal{S}_q^{SL_3}(\mathfrak{B})$. So the values of $\psi$ on these diagrams would determine $\psi$ on $\mathcal{S}_q^{SL_3}(\mathfrak{B})$. However, as we do not a priori have a definition of $\mathcal{S}_q^{SL_3}(\mathfrak{B})$ as a quotient of a free algebra by relators, it will be tricky to check that the map is well-defined. Instead, we have a definition of $\mathcal{S}_q^{SL_3}(\mathfrak{B})$ as a quotient of a free module and so we will define $\psi$ on any diagram by giving specific directions on how to write the diagram in terms of the diagrams $\alpha_{t(i)t(j)}$ and $\beta_{t(i)t(j)}$ and then check that this process leads to a well-defined map.

Given a diagram $D$, we obtain $\psi(D)$ by performing the following algorithm:

\begin{itemize}
	\item Apply $\Delta$ by splitting $D$ near the right boundary arc of $\mathfrak{B}$ so that $\Delta(D)$ is written as $$\Delta(D)=\sum  D_i \otimes E_i,$$ where the diagrams $E_i$ each contain only parallel and antiparallel strands.
	\item Apply $(\varepsilon \otimes \text{id})$ to $\Delta(D)$ to write $$(\varepsilon \otimes \text{id})\Delta(D)=\sum \varepsilon(D_i) E_i.$$
	\item Obtain $$\psi(D)=\sum \varepsilon(D_i)\psi(E_i) \in \mathcal{O}_q(SL_3),$$ where $\psi(E_i)$ is determined by the values of $\psi(\alpha_{t(i)t(j)})$ and $\psi(\beta_{t(i)t(j)})$ given above.
\end{itemize}

\begin{proposition}
The map $\psi: \mathcal{S}_q^{SL_3}(\mathfrak{B}) \rightarrow \mathcal{O}_q(SL_3)$ described above is a well-defined algebra homomorphsim.
\end{proposition}

\begin{proof}
We observe that if $\psi$ is well-defined, then it does respect the natural multiplication of diagrams in $\mathcal{S}_q^{SL_3}(\mathfrak{B}).$

We must check that the process outlined in the bulletpoints above respects the defining relations of the stated skein algebra. We split the relations into three cases: interior relations, boundary relations along the left boundary arc of $\mathfrak{B}$, boundary relations along the right boundary arc of $\mathfrak{B}.$

Consider a relation falling under the first two cases. Such a relation only affects the diagrams $D_i$ during the process. Since $\varepsilon$ is well-defined, application of such relations will result in identical representatives in $\mathcal{O}_q(SL_3)$, and so the process respects these relations.

The case of a relation along the right boundary arc of $\mathfrak{B}$ is more difficult since it will change the diagrams $E_i$ and will thus ultimately produce different representatives in $\mathcal{O}_q(SL_3)$. It is our task to show that these representatives are equivalent. We handle each relation separately.

% Why can we reduce to simple cases? horizontally and vertically

\underline{Relation (B1):}

To prove that $\psi$ respects relation (B1) it will suffice to check that 

\begin{equation*}
\psi\bigg(\begin{tikzpicture}[baseline=3ex]
%bigon
\draw [line width=1.5, rounded corners,->] (0,0)--(-1/2,1/4)--(-1/2,1/2+.1);
\draw [line width=1.5, rounded corners] (-1/2,1/2)--(-1/2,3/4)--(0,1);
\draw [line width=1.5, rounded corners,->] (0,0)--(1/2,1/4)--(1/2,1/2+.1);
\draw [line width=1.5, rounded corners] (1/2,1/2)--(1/2,3/4)--(0,1);
%strands
\draw (-1/2,1/2)--(1/2,1/2);
\draw [-<] (-1/2,1/2)--(0,1/2);
%punctures
\filldraw [fill=white, draw=black] (0,0) circle [radius=.05];
\filldraw [fill=white, draw=black] (0,1) circle [radius=.05];
%labels
\node [left] at (-1/2-.1,1/2) {$e$};
\node [right] at (1/2+.1,1/2) {$a+b$};
\end{tikzpicture}
\bigg)=(-1)^{a+b}q^{-1/3-(a+b)}\psi \bigg(\begin{tikzpicture}[baseline=3ex]
%bigon
\draw [line width=1.5, rounded corners,->] (0,0)--(-1/2,1/4)--(-1/2,1/2+.1);
\draw [line width=1.5, rounded corners] (-1/2,1/2)--(-1/2,3/4)--(0,1);
\draw [line width=1.5, rounded corners,->] (0,0)--(1/2,1/4)--(1/2,1/2+.1);
\draw [line width=1.5, rounded corners] (1/2,1/2)--(1/2,3/4)--(0,1);
%strands
\draw (-1/2,1/2)--(0,1/2);
\draw [-<] (-1/2,1/2)--(-1/4,1/2);
\draw (0,1/2)--(1/2,3/4);
\draw [->] (0,1/2)--(1/4,5/8);
\draw(0,1/2)--(1/2,1/4);
\draw[->] (0,1/2)--(1/4,3/8);
%punctures
\filldraw [fill=white, draw=black] (0,0) circle [radius=.05];
\filldraw [fill=white, draw=black] (0,1) circle [radius=.05];
%labels
\node [left] at (-1/2-.1,1/2) {$e$};
\node [right] at (1/2,3/4) {$b$};
\node [right] at (1/2,1/4) {$a$};
\end{tikzpicture}\bigg)
\end{equation*}

for any states $e,a,b \in \{-,0,+\}$ with $a<b.$

Fix such $e,a,b$ and let $i=t^{-1}(e)$ and $j=t^{-1}(a+b)$ be the corresponding integers in $\{1,2,3\}$. Then by the definition of $\psi,$ the left side of our relation is $(-q)^{j-i} A[4-i \vert 4-j]$.

We now compute the right side of the equation. It will be convenient to let $c,d$ be the unique states in $\{-,0,+\}$ such that $c<d$ and $c+d=e.$

By the definition of $\psi,$ we compute that

\begin{equation*}
\psi \bigg(\begin{tikzpicture}[baseline=3ex]
%bigon
\draw [line width=1.5, rounded corners,->] (0,0)--(-1/2,1/4)--(-1/2,1/2+.1);
\draw [line width=1.5, rounded corners] (-1/2,1/2)--(-1/2,3/4)--(0,1);
\draw [line width=1.5, rounded corners,->] (0,0)--(1/2,1/4)--(1/2,1/2+.1);
\draw [line width=1.5, rounded corners] (1/2,1/2)--(1/2,3/4)--(0,1);
%strands
\draw (-1/2,1/2)--(0,1/2);
\draw [-<] (-1/2,1/2)--(-1/4,1/2);
\draw (0,1/2)--(1/2,3/4);
\draw [->] (0,1/2)--(1/4,5/8);
\draw(0,1/2)--(1/2,1/4);
\draw[->] (0,1/2)--(1/4,3/8);
%punctures
\filldraw [fill=white, draw=black] (0,0) circle [radius=.05];
\filldraw [fill=white, draw=black] (0,1) circle [radius=.05];
%labels
\node [left] at (-1/2-.1,1/2) {$c+d$};
\node [right] at (1/2,3/4) {$b$};
\node [right] at (1/2,1/4) {$a$};
\end{tikzpicture}\bigg)=\displaystyle\sum_{x,y} \varepsilon\bigg(\begin{tikzpicture}[baseline=3ex]
%bigon
\draw [line width=1.5, rounded corners,->] (0,0)--(-1/2,1/4)--(-1/2,1/2+.1);
\draw [line width=1.5, rounded corners] (-1/2,1/2)--(-1/2,3/4)--(0,1);
\draw [line width=1.5, rounded corners,->] (0,0)--(1/2,1/4)--(1/2,1/2+.1);
\draw [line width=1.5, rounded corners] (1/2,1/2)--(1/2,3/4)--(0,1);
%strands
\draw (-1/2,1/2)--(0,1/2);
\draw [-<] (-1/2,1/2)--(-1/4,1/2);
\draw (0,1/2)--(1/2,3/4);
\draw [->] (0,1/2)--(1/4,5/8);
\draw(0,1/2)--(1/2,1/4);
\draw[->] (0,1/2)--(1/4,3/8);
%punctures
\filldraw [fill=white, draw=black] (0,0) circle [radius=.05];
\filldraw [fill=white, draw=black] (0,1) circle [radius=.05];
%labels
\node [left] at (-1/2-.1,1/2) {$c+d$};
\node [right] at (1/2,3/4) {$x$};
\node [right] at (1/2,1/4) {$y$};
\end{tikzpicture}\bigg)\psi \bigg(\begin{tikzpicture}[baseline=3ex]
%bigon
\draw [line width=1.5, rounded corners,->] (0,0)--(-1/2,1/4)--(-1/2,1/2+.1);
\draw [line width=1.5, rounded corners] (-1/2,1/2)--(-1/2,3/4)--(0,1);
\draw [line width=1.5, rounded corners,->] (0,0)--(1/2,1/4)--(1/2,1/2+.1);
\draw [line width=1.5, rounded corners] (1/2,1/2)--(1/2,3/4)--(0,1);
%strands
\draw (-1/2,3/4)--(1/2,3/4);
\draw [->] (-1/2,3/4)--(0,3/4);
\draw (-1/2,1/4)--(1/2,1/4);
\draw [->] (-1/2,1/4)--(0,1/4);
%punctures
\filldraw [fill=white, draw=black] (0,0) circle [radius=.05];
\filldraw [fill=white, draw=black] (0,1) circle [radius=.05];
%labels
\node [left] at (-1/2,3/4) {$x$};
\node [left] at (-1/2,1/4) {$y$};
\node [right] at (1/2,3/4) {$b$};
\node [right] at (1/2, 1/4) {$a$};
\end{tikzpicture} \bigg)
\end{equation*}

We will denote the values of the counit appearing in the above equation as $\varepsilon_{c+d,x,y}.$ We use (B3) and (B1) to compute that $\varepsilon_{c+d,x,y}=0$ unless $\{x,y\}=\{c,d\}$ and we use (B2) to see that $\varepsilon_{c+d,c,d}=-q\varepsilon_{c+d,d,c}.$ We also use (B1) to compute that $\varepsilon_{c+d,d,c}=(-1)^{c+d}q^{1/3+(c+d)}.$

The right side of our relation becomes

\begin{align*}
&=(-1)^{a+b}q^{-1/3-(a+b)}(-1)^{c+d}q^{1/3+(c+d)}\bigg( \psi \bigg( \begin{tikzpicture}[baseline=3ex]
%bigon
\draw [line width=1.5, rounded corners,->] (0,0)--(-1/2,1/4)--(-1/2,1/2+.1);
\draw [line width=1.5, rounded corners] (-1/2,1/2)--(-1/2,3/4)--(0,1);
\draw [line width=1.5, rounded corners,->] (0,0)--(1/2,1/4)--(1/2,1/2+.1);
\draw [line width=1.5, rounded corners] (1/2,1/2)--(1/2,3/4)--(0,1);
%strands
\draw (-1/2,3/4)--(1/2,3/4);
\draw [->] (-1/2,3/4)--(0,3/4);
\draw (-1/2,1/4)--(1/2,1/4);
\draw [->] (-1/2,1/4)--(0,1/4);
%punctures
\filldraw [fill=white, draw=black] (0,0) circle [radius=.05];
\filldraw [fill=white, draw=black] (0,1) circle [radius=.05];
%labels
\node [left] at (-1/2,3/4) {$d$};
\node [left] at (-1/2,1/4) {$c$};
\node [right] at (1/2,3/4) {$b$};
\node [right] at (1/2, 1/4) {$a$};
\end{tikzpicture} \bigg)-q \psi \bigg( \begin{tikzpicture}[baseline=3ex]
%bigon
\draw [line width=1.5, rounded corners,->] (0,0)--(-1/2,1/4)--(-1/2,1/2+.1);
\draw [line width=1.5, rounded corners] (-1/2,1/2)--(-1/2,3/4)--(0,1);
\draw [line width=1.5, rounded corners,->] (0,0)--(1/2,1/4)--(1/2,1/2+.1);
\draw [line width=1.5, rounded corners] (1/2,1/2)--(1/2,3/4)--(0,1);
%strands
\draw (-1/2,3/4)--(1/2,3/4);
\draw [->] (-1/2,3/4)--(0,3/4);
\draw (-1/2,1/4)--(1/2,1/4);
\draw [->] (-1/2,1/4)--(0,1/4);
%punctures
\filldraw [fill=white, draw=black] (0,0) circle [radius=.05];
\filldraw [fill=white, draw=black] (0,1) circle [radius=.05];
%labels
\node [left] at (-1/2,3/4) {$c$};
\node [left] at (-1/2,1/4) {$d$};
\node [right] at (1/2,3/4) {$b$};
\node [right] at (1/2, 1/4) {$a$};
\end{tikzpicture}\bigg) \bigg)\\
&=(-q)^{(c+d)-(a+b)}\bigg( \psi \bigg( \begin{tikzpicture}[baseline=3ex]
%bigon
\draw [line width=1.5, rounded corners,->] (0,0)--(-1/2,1/4)--(-1/2,1/2+.1);
\draw [line width=1.5, rounded corners] (-1/2,1/2)--(-1/2,3/4)--(0,1);
\draw [line width=1.5, rounded corners,->] (0,0)--(1/2,1/4)--(1/2,1/2+.1);
\draw [line width=1.5, rounded corners] (1/2,1/2)--(1/2,3/4)--(0,1);
%strands
\draw (-1/2,3/4)--(1/2,3/4);
\draw [->] (-1/2,3/4)--(0,3/4);
\draw (-1/2,1/4)--(1/2,1/4);
\draw [->] (-1/2,1/4)--(0,1/4);
%punctures
\filldraw [fill=white, draw=black] (0,0) circle [radius=.05];
\filldraw [fill=white, draw=black] (0,1) circle [radius=.05];
%labels
\node [left] at (-1/2,3/4) {$d$};
\node [left] at (-1/2,1/4) {$c$};
\node [right] at (1/2,3/4) {$b$};
\node [right] at (1/2, 1/4) {$a$};
\end{tikzpicture} \bigg)-q \psi \bigg( \begin{tikzpicture}[baseline=3ex]
%bigon
\draw [line width=1.5, rounded corners,->] (0,0)--(-1/2,1/4)--(-1/2,1/2+.1);
\draw [line width=1.5, rounded corners] (-1/2,1/2)--(-1/2,3/4)--(0,1);
\draw [line width=1.5, rounded corners,->] (0,0)--(1/2,1/4)--(1/2,1/2+.1);
\draw [line width=1.5, rounded corners] (1/2,1/2)--(1/2,3/4)--(0,1);
%strands
\draw (-1/2,3/4)--(1/2,3/4);
\draw [->] (-1/2,3/4)--(0,3/4);
\draw (-1/2,1/4)--(1/2,1/4);
\draw [->] (-1/2,1/4)--(0,1/4);
%punctures
\filldraw [fill=white, draw=black] (0,0) circle [radius=.05];
\filldraw [fill=white, draw=black] (0,1) circle [radius=.05];
%labels
\node [left] at (-1/2,3/4) {$c$};
\node [left] at (-1/2,1/4) {$d$};
\node [right] at (1/2,3/4) {$b$};
\node [right] at (1/2, 1/4) {$a$};
\end{tikzpicture}\bigg) \bigg)
\end{align*}

We check that this formula agrees with $$(-q)^{t^{-1}(a+b)-t^{-1}(c+d)}(X_{t^{-1}(d)t^{-1}(b)}X_{t^{-1}(c)t^{-1}(a)}-qX_{t^{-1}(c)t^{-1}(b)}X_{t^{-1}(d)t^{-1}(a)})$$ which is $$(-q)^{(j-i)}A[4-i\vert 4-j],$$ as required.

\underline{Relation (B2):}

To show that $\psi$ respects relation (B2) it suffices to check that the following relation holds in $\mathcal{O}_q(SL_3).$

\begin{equation*}
\psi \bigg(\begin{tikzpicture}[baseline=3ex]
%bigon
\draw [line width=1.5, rounded corners,->] (0,0)--(-1/2,1/4)--(-1/2,1/2+.1);
\draw [line width=1.5, rounded corners] (-1/2,1/2)--(-1/2,3/4)--(0,1);
\draw [line width=1.5, rounded corners,->] (0,0)--(1/2,1/4)--(1/2,1/2+.1);
\draw [line width=1.5, rounded corners] (1/2,1/2)--(1/2,3/4)--(0,1);
%strands
\draw (-1/2,3/4)--(1/2,3/4);
\draw [->] (-1/2,3/4)--(0,3/4);
\draw (-1/2,1/4)--(1/2,1/4);
\draw [->] (-1/2,1/4)--(0,1/4);
%punctures
\filldraw [fill=white, draw=black] (0,0) circle [radius=.05];
\filldraw [fill=white, draw=black] (0,1) circle [radius=.05];
%labels
\node [left] at (-1/2,3/4) {$t(i)$};
\node [left] at (-1/2,1/4) {$t(j)$};
\node [right] at (1/2,3/4) {$t(m)$};
\node [right] at (1/2, 1/4) {$t(n)$};
\end{tikzpicture} \bigg)=q^{-1} \psi \bigg( \begin{tikzpicture}[baseline=3ex]
%bigon
\draw [line width=1.5, rounded corners,->] (0,0)--(-1/2,1/4)--(-1/2,1/2+.1);
\draw [line width=1.5, rounded corners] (-1/2,1/2)--(-1/2,3/4)--(0,1);
\draw [line width=1.5, rounded corners,->] (0,0)--(1/2,1/4)--(1/2,1/2+.1);
\draw [line width=1.5, rounded corners] (1/2,1/2)--(1/2,3/4)--(0,1);
%strands
\draw (-1/2,3/4)--(1/2,3/4);
\draw [->] (-1/2,3/4)--(0,3/4);
\draw (-1/2,1/4)--(1/2,1/4);
\draw [->] (-1/2,1/4)--(0,1/4);
%punctures
\filldraw [fill=white, draw=black] (0,0) circle [radius=.05];
\filldraw [fill=white, draw=black] (0,1) circle [radius=.05];
%labels
\node [left] at (-1/2,3/4) {$t(i)$};
\node [left] at (-1/2,1/4) {$t(j)$};
\node [right] at (1/2,3/4) {$t(n)$};
\node [right] at (1/2, 1/4) {$t(m)$};
\end{tikzpicture}\bigg) +q^{-3} \psi \bigg( \begin{tikzpicture}[baseline=3ex]
%bigon
\draw [line width=1.5, rounded corners,->] (0,0)--(-1/2,1/4)--(-1/2,1/2+.1);
\draw [line width=1.5, rounded corners] (-1/2,1/2)--(-1/2,3/4)--(0,1);
\draw [line width=1.5, rounded corners,->] (0,0)--(1/2,1/4)--(1/2,1/2+.1);
\draw [line width=1.5, rounded corners] (1/2,1/2)--(1/2,3/4)--(0,1);
%strands
\draw (-1/2,3/4)--(-1/4,1/2);
\draw (-1/2,1/4)--(-1/4,1/2);
\draw (1/2,3/4)--(1/4,1/2);
\draw (1/2,1/4)--(1/4,1/2);
\draw (-1/4,1/2)--(1/4,1/2);
\draw [-<](-1/4,1/2)--(0,1/2);
%punctures
\filldraw [fill=white, draw=black] (0,0) circle [radius=.05];
\filldraw [fill=white, draw=black] (0,1) circle [radius=.05];
%labels
\node [left] at (-1/2,3/4) {$t(i)$};
\node [left] at (-1/2,1/4) {$t(j)$};
\node [right] at (1/2,3/4) {$t(n)$};
\node [right] at (1/2, 1/4) {$t(m)$};
\end{tikzpicture}\bigg)
\end{equation*}

for $i,j,m,n \in \{1,2,3\}$ such that $n<m.$ So we must show that

\begin{equation*}
\psi \bigg( \begin{tikzpicture}[baseline=3ex]
%bigon
\draw [line width=1.5, rounded corners,->] (0,0)--(-1/2,1/4)--(-1/2,1/2+.1);
\draw [line width=1.5, rounded corners] (-1/2,1/2)--(-1/2,3/4)--(0,1);
\draw [line width=1.5, rounded corners,->] (0,0)--(1/2,1/4)--(1/2,1/2+.1);
\draw [line width=1.5, rounded corners] (1/2,1/2)--(1/2,3/4)--(0,1);
%strands
\draw (-1/2,3/4)--(-1/4,1/2);
\draw (-1/2,1/4)--(-1/4,1/2);
\draw (1/2,3/4)--(1/4,1/2);
\draw (1/2,1/4)--(1/4,1/2);
\draw (-1/4,1/2)--(1/4,1/2);
\draw [-<](-1/4,1/2)--(0,1/2);
%punctures
\filldraw [fill=white, draw=black] (0,0) circle [radius=.05];
\filldraw [fill=white, draw=black] (0,1) circle [radius=.05];
%labels
\node [left] at (-1/2,3/4) {$t(i)$};
\node [left] at (-1/2,1/4) {$t(j)$};
\node [right] at (1/2,3/4) {$t(n)$};
\node [right] at (1/2, 1/4) {$t(m)$};
\end{tikzpicture}\bigg)=q^3X_{im}X_{jn}-q^2X_{in}X_{jm}.
\end{equation*}

From relation (I1a) and the computations of $\varepsilon(\beta_{st})$ from Section 5, we compute that

\begin{equation*}
\varepsilon \bigg( \begin{tikzpicture}[baseline=3ex]
%bigon
\draw [line width=1.5, rounded corners,->] (0,0)--(-1/2,1/4)--(-1/2,1/2+.1);
\draw [line width=1.5, rounded corners] (-1/2,1/2)--(-1/2,3/4)--(0,1);
\draw [line width=1.5, rounded corners,->] (0,0)--(1/2,1/4)--(1/2,1/2+.1);
\draw [line width=1.5, rounded corners] (1/2,1/2)--(1/2,3/4)--(0,1);
%strands
\draw (-1/2,3/4)--(-1/4,1/2);
\draw (-1/2,1/4)--(-1/4,1/2);
\draw (1/2,3/4)--(1/4,1/2);
\draw (1/2,1/4)--(1/4,1/2);
\draw (-1/4,1/2)--(1/4,1/2);
\draw [-<](-1/4,1/2)--(0,1/2);
%punctures
\filldraw [fill=white, draw=black] (0,0) circle [radius=.05];
\filldraw [fill=white, draw=black] (0,1) circle [radius=.05];
%labels
\node [left] at (-1/2,3/4) {$t(i)$};
\node [left] at (-1/2,1/4) {$t(j)$};
\node [right] at (1/2,3/4) {$t(n)$};
\node [right] at (1/2, 1/4) {$t(m)$};
\end{tikzpicture}\bigg)=q^{3+1/3}R_{ij}^{kl}-q^4\delta_{ik}\delta_{jl}.
\end{equation*}

Thus, we must show that $$(\displaystyle\sum_{k,l} q^{3+1/3}R_{ij}^{kl}X_{kn}X_{lm})-q^4X_{in}X_{jm}=q^3X_{im}X_{jn}-q^2X_{in}X_{jm}.$$

We apply the identity $$\displaystyle\sum_{k,l}R_{ij}^{kl}X_{kn}X_{lm}=\displaystyle\sum_{k,l} R_{kl}^{nm}X_{ik}X_{jl}.$$ Since $n<m$, we have that $R_{nm}^{nm}=q^{-1/3}(q-q^{-1})$ and $R_{mn}^{nm}=q^{-1/3}$ are the only nonzero values of $R_{kl}^{nm}$ as $k$ and $l$ vary.

The left side of our equation now becomes \begin{align*} & (\displaystyle\sum_{k,l} q^{3+1/3}R_{ij}^{kl}X_{kn}X_{lm})-q^4X_{in}X_{jm}\\&=(\displaystyle\sum_{k,l} q^{3+1/3}R_{kl}^{nm}X_{ik}X_{jl})-q^4X_{in}X_{jm} \\
&=q^{3}(q-q^{-1})X_{in}X_{jm}+q^{3}X_{im}X_{jn}-q^{4}X_{in}X_{jm}\\ &=q^{3}X_{im}X_{jn}-q^{2}X_{in}X_{jm},
\end{align*} as required. So $\psi$ respects (B2).

\underline{Relation (B3):}

To show that $\psi$ respects (B3) we need to show that

\begin{equation*}
\psi \bigg( \begin{tikzpicture}[baseline=3ex]
%bigon
\draw [line width=1.5, rounded corners,->] (0,0)--(-1/2,1/4)--(-1/2,1/2+.1);
\draw [line width=1.5, rounded corners] (-1/2,1/2)--(-1/2,3/4)--(0,1);
\draw [line width=1.5, rounded corners,->] (0,0)--(1/2,1/4)--(1/2,1/2+.1);
\draw [line width=1.5, rounded corners] (1/2,1/2)--(1/2,3/4)--(0,1);
%strands
\draw (-1/2,1/2)--(0,1/2);
\draw [-<] (-1/2,1/2)--(-1/4,1/2);
\draw (0,1/2)--(1/2,3/4);
\draw [->] (0,1/2)--(1/4,5/8);
\draw(0,1/2)--(1/2,1/4);
\draw[->] (0,1/2)--(1/4,3/8);
%punctures
\filldraw [fill=white, draw=black] (0,0) circle [radius=.05];
\filldraw [fill=white, draw=black] (0,1) circle [radius=.05];
%labels
\node [left] at (-1/2-.1,1/2) {$t(i)$};
\node [right] at (1/2,3/4) {$t(j)$};
\node [right] at (1/2,1/4) {$t(j)$};
\end{tikzpicture}\bigg)=0
\end{equation*}

for any $i,j \in \{1,2,3\}.$

By the definition of $\psi,$ we have

\begin{equation*}
\psi \bigg( \begin{tikzpicture}[baseline=3ex]
%bigon
\draw [line width=1.5, rounded corners,->] (0,0)--(-1/2,1/4)--(-1/2,1/2+.1);
\draw [line width=1.5, rounded corners] (-1/2,1/2)--(-1/2,3/4)--(0,1);
\draw [line width=1.5, rounded corners,->] (0,0)--(1/2,1/4)--(1/2,1/2+.1);
\draw [line width=1.5, rounded corners] (1/2,1/2)--(1/2,3/4)--(0,1);
%strands
\draw (-1/2,1/2)--(0,1/2);
\draw [-<] (-1/2,1/2)--(-1/4,1/2);
\draw (0,1/2)--(1/2,3/4);
\draw [->] (0,1/2)--(1/4,5/8);
\draw(0,1/2)--(1/2,1/4);
\draw[->] (0,1/2)--(1/4,3/8);
%punctures
\filldraw [fill=white, draw=black] (0,0) circle [radius=.05];
\filldraw [fill=white, draw=black] (0,1) circle [radius=.05];
%labels
\node [left] at (-1/2-.1,1/2) {$t(i)$};
\node [right] at (1/2,3/4) {$t(j)$};
\node [right] at (1/2,1/4) {$t(j)$};
\end{tikzpicture}\bigg)= \displaystyle\sum_{k,l} \varepsilon \bigg( \begin{tikzpicture}[baseline=3ex]
%bigon
\draw [line width=1.5, rounded corners,->] (0,0)--(-1/2,1/4)--(-1/2,1/2+.1);
\draw [line width=1.5, rounded corners] (-1/2,1/2)--(-1/2,3/4)--(0,1);
\draw [line width=1.5, rounded corners,->] (0,0)--(1/2,1/4)--(1/2,1/2+.1);
\draw [line width=1.5, rounded corners] (1/2,1/2)--(1/2,3/4)--(0,1);
%strands
\draw (-1/2,1/2)--(0,1/2);
\draw [-<] (-1/2,1/2)--(-1/4,1/2);
\draw (0,1/2)--(1/2,3/4);
\draw [->] (0,1/2)--(1/4,5/8);
\draw(0,1/2)--(1/2,1/4);
\draw[->] (0,1/2)--(1/4,3/8);
%punctures
\filldraw [fill=white, draw=black] (0,0) circle [radius=.05];
\filldraw [fill=white, draw=black] (0,1) circle [radius=.05];
%labels
\node [left] at (-1/2-.1,1/2) {$t(i)$};
\node [right] at (1/2,3/4) {$t(k)$};
\node [right] at (1/2,1/4) {$t(l)$};
\end{tikzpicture}\bigg) X_{kj}X_{lj}.
\end{equation*}

We compute that

\begin{equation*}
\varepsilon \bigg( \begin{tikzpicture}[baseline=3ex]
%bigon
\draw [line width=1.5, rounded corners,->] (0,0)--(-1/2,1/4)--(-1/2,1/2+.1);
\draw [line width=1.5, rounded corners] (-1/2,1/2)--(-1/2,3/4)--(0,1);
\draw [line width=1.5, rounded corners,->] (0,0)--(1/2,1/4)--(1/2,1/2+.1);
\draw [line width=1.5, rounded corners] (1/2,1/2)--(1/2,3/4)--(0,1);
%strands
\draw (-1/2,1/2)--(0,1/2);
\draw [-<] (-1/2,1/2)--(-1/4,1/2);
\draw (0,1/2)--(1/2,3/4);
\draw [->] (0,1/2)--(1/4,5/8);
\draw(0,1/2)--(1/2,1/4);
\draw[->] (0,1/2)--(1/4,3/8);
%punctures
\filldraw [fill=white, draw=black] (0,0) circle [radius=.05];
\filldraw [fill=white, draw=black] (0,1) circle [radius=.05];
%labels
\node [left] at (-1/2-.1,1/2) {$t(i)$};
\node [right] at (1/2,3/4) {$t(k)$};
\node [right] at (1/2,1/4) {$t(l)$};
\end{tikzpicture}\bigg)=0
\end{equation*}
if $4-i$ is in $\{k,l\}$ or if $k=l.$

If $l<k$ we have $\varepsilon_{ikl}=-q\varepsilon_{ilk}.$ This can be computed by using relations (B2) and (I3).

Thus, $$\psi \bigg( \begin{tikzpicture}[baseline=3ex]
%bigon
\draw [line width=1.5, rounded corners,->] (0,0)--(-1/2,1/4)--(-1/2,1/2+.1);
\draw [line width=1.5, rounded corners] (-1/2,1/2)--(-1/2,3/4)--(0,1);
\draw [line width=1.5, rounded corners,->] (0,0)--(1/2,1/4)--(1/2,1/2+.1);
\draw [line width=1.5, rounded corners] (1/2,1/2)--(1/2,3/4)--(0,1);
%strands
\draw (-1/2,1/2)--(0,1/2);
\draw [-<] (-1/2,1/2)--(-1/4,1/2);
\draw (0,1/2)--(1/2,3/4);
\draw [->] (0,1/2)--(1/4,5/8);
\draw(0,1/2)--(1/2,1/4);
\draw[->] (0,1/2)--(1/4,3/8);
%punctures
\filldraw [fill=white, draw=black] (0,0) circle [radius=.05];
\filldraw [fill=white, draw=black] (0,1) circle [radius=.05];
%labels
\node [left] at (-1/2-.1,1/2) {$t(i)$};
\node [right] at (1/2,3/4) {$t(j)$};
\node [right] at (1/2,1/4) {$t(j)$};
\end{tikzpicture}\bigg)=\varepsilon_{ilk}(X_{lj}X_{kj}-qX_{kj}X_{lj})$$ for the unique suitable pair $l,k$ for which $\varepsilon_{ilk}$ is nonzero. The result follows from the identity $$X_{lj}X_{kj}=qX_{kj}X_{lj}$$ which holds in $\mathcal{O}_q(SL_3)$ for $l<k.$

\underline{Relation (B4):}

To check that $\psi$ respects relation (B4) it suffices to check 

\begin{equation*}
\psi \bigg( \begin{tikzpicture}[baseline=3ex]
%bigon
\draw [line width=1.5, rounded corners,->] (0,0)--(-1/2,1/4)--(-1/2,1/2+.1);
\draw [line width=1.5, rounded corners] (-1/2,1/2)--(-1/2,3/4)--(0,1);
\draw [line width=1.5, rounded corners,->] (0,0)--(1/2,1/4)--(1/2,1/2+.1);
\draw [line width=1.5, rounded corners] (1/2,1/2)--(1/2,3/4)--(0,1);
%strands
\draw (0,1/2)--(1/2,3/4);
\draw [->] (0,1/2)--(1/4,5/8);
\draw (0,1/2)--(1/2,1/2);
\draw (0,1/2)--(1/2,1/4);
\draw [->] (0,1/2)--(1/4,3/8);
%punctures
\filldraw [fill=white, draw=black] (0,0) circle [radius=.05];
\filldraw [fill=white, draw=black] (0,1) circle [radius=.05];
%labels
\node [right] at (1/2,3/4+.1) {$t(1)$};
\node [right] at (1/2,1/2) {$t(2)$};
\node [right] at (1/2,1/4-.1) {$t(3)$};
\end{tikzpicture}\bigg)=q^{-2}.
\end{equation*}

By the definition of $\psi,$ we compute

\begin{equation*}
\psi \bigg( \begin{tikzpicture}[baseline=3ex]
%bigon
\draw [line width=1.5, rounded corners,->] (0,0)--(-1/2,1/4)--(-1/2,1/2+.1);
\draw [line width=1.5, rounded corners] (-1/2,1/2)--(-1/2,3/4)--(0,1);
\draw [line width=1.5, rounded corners,->] (0,0)--(1/2,1/4)--(1/2,1/2+.1);
\draw [line width=1.5, rounded corners] (1/2,1/2)--(1/2,3/4)--(0,1);
%strands
\draw (0,1/2)--(1/2,3/4);
\draw [->] (0,1/2)--(1/4,5/8);
\draw (0,1/2)--(1/2,1/2);
\draw (0,1/2)--(1/2,1/4);
\draw [->] (0,1/2)--(1/4,3/8);
%punctures
\filldraw [fill=white, draw=black] (0,0) circle [radius=.05];
\filldraw [fill=white, draw=black] (0,1) circle [radius=.05];
%labels
\node [right] at (1/2,3/4+.1) {$t(1)$};
\node [right] at (1/2,1/2) {$t(2)$};
\node [right] at (1/2,1/4-.1) {$t(3)$};
\end{tikzpicture}\bigg)=\displaystyle \sum_{\sigma \in S_3} \varepsilon \bigg(\begin{tikzpicture}[baseline=3ex]
%bigon
\draw [line width=1.5, rounded corners,->] (0,0)--(-1/2,1/4)--(-1/2,1/2+.1);
\draw [line width=1.5, rounded corners] (-1/2,1/2)--(-1/2,3/4)--(0,1);
\draw [line width=1.5, rounded corners,->] (0,0)--(1/2,1/4)--(1/2,1/2+.1);
\draw [line width=1.5, rounded corners] (1/2,1/2)--(1/2,3/4)--(0,1);
%strands
\draw (0,1/2)--(1/2,3/4);
\draw [->] (0,1/2)--(1/4,5/8);
\draw (0,1/2)--(1/2,1/2);
\draw (0,1/2)--(1/2,1/4);
\draw [->] (0,1/2)--(1/4,3/8);
%punctures
\filldraw [fill=white, draw=black] (0,0) circle [radius=.05];
\filldraw [fill=white, draw=black] (0,1) circle [radius=.05];
%labels
\node [right] at (1/2,3/4+.1) {$t(\sigma_1)$};
\node [right] at (1/2,1/2) {$t(\sigma_2)$};
\node [right] at (1/2,1/4-.1) {$t(\sigma_3)$};
\end{tikzpicture} \bigg)X_{\sigma_1 1}X_{\sigma_2 2}X_{\sigma_3 3}.
\end{equation*}

We see that this is equal to \begin{align*}q^{-2} \sum_{\sigma \in S_3} (-q)^{\l(\sigma)} X_{\sigma_1 1}X_{\sigma_2 2}X_{\sigma_3 3} &=q^{-2}\text{det}_q\\
&=q^{-2}.
\end{align*}

So we see that $\psi$ respects (B4) and, thus, $\psi$ is well-defined.
\end{proof}

Our previous two propositions allow us to state the following theorem.

\begin{theorem}
We have that $$\mathcal{S}_q^{SL_3}(\mathfrak{B}) \cong \mathcal{O}_q(SL_3)$$ as Hopf algebras.
\end{theorem}

\begin{proof}
In Proposition 8 we showed that $\phi$ is a well-defined map of bialgebras. To show that $\phi$ is an isomorphism, it suffices to show that $\phi$ is invertible as a map of $\mathcal{R}\text{-modules}.$ We claim that $\psi$ is its inverse.

We observe that $\psi \circ \phi(X_{ij})=X_{ij}$ for all generators $X_{ij}$ of $\mathcal{O}_q(SL_3).$ Since $\psi$ and $\phi$ are both algebra maps, this implies that $$\psi \circ \phi = \text{id}_{\mathcal{O}_q(SL_3)}.$$

Similarly, $\phi \circ \psi$ agrees with $\text{id}_{\mathcal{S}_q^{SL_3}(\mathfrak{B})}$ for all generating diagrams $\alpha_{st}$ and $\beta_{st}.$ Thus, $$\phi \circ \psi = \text{id}_{\mathcal{S}_q^{SL_3}(\mathfrak{B})}.$$
	
Thus, $\mathcal{O}_q(SL_3)$ and $\mathcal{S}_q^{SL_3}(\mathfrak{B})$ are isomorphic as bialgebras. Since $\mathcal{O}_q(SL_3)$ is a Hopf algebra, then $\mathcal{O}_q(SL_3)$ and $\mathcal{S}_q^{SL_3}(\mathfrak{B})$ are isomorphic as Hopf algebras.
\end{proof}

\section{The stated skein algebra of the triangle}

The Hopf algebra $\mathcal{O}_q(SL_3)$ is equipped with a cobraiding $\rho: \mathcal{O}_q(SL_3) \otimes \mathcal{O}_q(SL_3) \rightarrow \mathcal{R}$. In \cite{CL19} the cobraiding for the $SL_2$ case was shown to have a simple diagrammatic definition, and an analogous definition will work here as well. This cobraiding will allow us to describe the $SL_3$ stated skein algebra of the triangle, $\mathfrak{T}.$

We define the cobraiding $\rho: \mathcal{S}_q^{SL_3} (\mathfrak{B})\otimes \mathcal{S}_q^{SL_3}(\mathfrak{B}) \rightarrow \mathcal{R}$ on diagrams by 

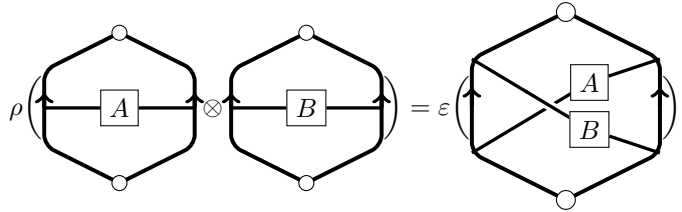
\begin{equation*}
\rho \bigg(\begin{tikzpicture}[scale=2, baseline=6ex]
%bigon
\draw [line width=1.5, rounded corners,->] (0,0)--(-1/2,1/4)--(-1/2,1/2+.1);
\draw [line width=1.5, rounded corners] (-1/2,1/2)--(-1/2,3/4)--(0,1);
\draw [line width=1.5, rounded corners,->] (0,0)--(1/2,1/4)--(1/2,1/2+.1);
\draw [line width=1.5, rounded corners] (1/2,1/2)--(1/2,3/4)--(0,1);
%punctures
\filldraw [fill=white, draw=black] (0,0) circle [radius=.05];
\filldraw [fill=white, draw=black] (0,1) circle [radius=.05];
\node at (0,1/2) [draw] (A) {$A$};
\draw [line width=1.25] (-1/2,1/2)--(A)--(1/2,1/2);
\end{tikzpicture} \otimes \begin{tikzpicture}[scale=2, baseline=6ex]
%bigon
\draw [line width=1.5, rounded corners,->] (0,0)--(-1/2,1/4)--(-1/2,1/2+.1);
\draw [line width=1.5, rounded corners] (-1/2,1/2)--(-1/2,3/4)--(0,1);
\draw [line width=1.5, rounded corners,->] (0,0)--(1/2,1/4)--(1/2,1/2+.1);
\draw [line width=1.5, rounded corners] (1/2,1/2)--(1/2,3/4)--(0,1);
%punctures
\filldraw [fill=white, draw=black] (0,0) circle [radius=.05];
\filldraw [fill=white, draw=black] (0,1) circle [radius=.05];
\node at (0,1/2) [draw] (B) {$B$};
\draw [line width=1.25] (-1/2,1/2)--(B)--(1/2,1/2);
\end{tikzpicture}  \bigg)= \varepsilon \bigg(
\begin{tikzpicture}[use Hobby shortcut, scale=2.5, baseline=2.5*3ex]
%bigon
\draw [line width=1.5, rounded corners,->] (0,0)--(-1/2,1/4)--(-1/2,1/2+.1);
\draw [line width=1.5, rounded corners] (-1/2,1/2)--(-1/2,3/4)--(0,1);
\draw [line width=1.5, rounded corners,->] (0,0)--(1/2,1/4)--(1/2,1/2+.1);
\draw [line width=1.5, rounded corners] (1/2,1/2)--(1/2,3/4)--(0,1);
%punctures
\filldraw [fill=white, draw=black] (0,0) circle [radius=.05];
\filldraw [fill=white, draw=black] (0,1) circle [radius=.05];
\begin{knot}[
consider self intersections=true,
%  draft mode=crossings,
flip crossing=1,
ignore endpoint intersections=false]
\strand [line width=1.25] (1/8,5/8)..(-1/2,1/4);
\strand [line width=1.25] (1/8,3/8)..(-1/2,3/4);
\end{knot}
\node at (1/8,5/8) [fill=white, draw] (A) {$A$};
\node at (1/8,3/8) [fill=white, draw] (B) {$B$};
%strands
\draw [line width=1.25]  (A)--(1/2,3/4);
\draw [line width=1.25]  (B)--(1/2,1/4);
\end{tikzpicture}\bigg).
\end{equation*}

In the diagrams above, the strands depict a bundle of parallel or antiparallel strands. The diagrammatic definition of the map makes it easy to see that it respects the defining relations of the stated skein algebra, so it is well-defined. The argument that this satisfies the cobraiding axioms is identical to the one in Section 3.7 of \cite{CL19}, but we do not need to use it in this paper.

We recall that a cobraiding is determined by its values on a set of generators and so we see that the map $\rho$ that we have defined diagrammatically satisfies $$\rho (X_{ij} \otimes X_{kl})=R_{ik}^{jl},$$ and thus matches up with the standard co-R-matrix.

In the situation that we have two algebras $M$ and $N$ which are both left comodule-algebras over $\mathcal{O}_q(SL_3)$ we can endow the vector space $M \otimes N$ with a left comodule-algebra structure using the cobraiding $\rho.$ We will denote this algebra by $M \underset{-}{\otimes} N$ and call it the \textit{braided tensor product} of the algebras $M$ and $N$. Using Sweedler's notation, its multiplication is defined as follows:

\begin{equation*}
(x\otimes y)\star (z \otimes t)  = (x \otimes 1)(\sum_{(z)(y)} \rho(z'\otimes y') (z'' \otimes y''))(1 \otimes t)
\end{equation*}

Equivalently, if we identify $M$ with $M \otimes \{1\}$ and $N$ with $\{1 \} \otimes N,$ then our product structure is given by

$$xy= \begin{cases}
xy & \text{if $x,y$ both in $M$ or both in $N$}\\
x \otimes y & \text{ if $x$ in $M$ and $y$ in $M$}\\
\displaystyle \sum_{(x)(y)} \rho (y' \otimes x') (y'' \otimes x'') &  \text{ if $x$ in $N$ and $y$ in $M$}\\
\end{cases}$$

Le and Costantino showed in \cite{CL19} that gluing disjoint surfaces along a triangle yields a braided tensor product of stated skein algebras for the $SL_2$ case. The same is true for the $SL_3$ case and Proposition 7 from this paper takes care of most of the work we need to do to show it.

\begin{theorem}
Let $\Sigma_1$ and $\Sigma_2$ be disjoint punctured bordered surfaces. If $a$ is a boundary arc of $\Sigma_1$ and $b$ is a boundary arc of $\Sigma_2$, then we have an algebra isomorphism
\begin{equation*}
\mathcal{S}_q^{SL_3}(\Sigma_1) \underset{-}{\otimes}\mathcal{S}_q^{SL_3}(\Sigma_2) \cong \mathcal{S}_q^{SL_3}((\Sigma_1 \sqcup \Sigma_2)\# \mathfrak{T})
\end{equation*} given by the map $\text{glue}_\mathfrak{T}$ defined in Section 6.
\end{theorem}

\begin{proof}
By Proposition 7, the map $$\text{glue}_\mathfrak{T}: \mathcal{S}_q^{SL_3}(\Sigma_1 \sqcup \Sigma_2)\rightarrow \mathcal{S}_q^{SL_3}((\Sigma_1 \sqcup \Sigma_2)\#\mathfrak{T})$$ is an $\mathcal{R}\text{-module}$ isomorphism. Since $S_q^{SL_3}(\Sigma_1 \sqcup \Sigma_2)$ is naturally isomorphic to $S_q^{SL_3}(\Sigma_1) \otimes S_q^{SL_3}(\Sigma_2),$ we see that the isomorphism claimed in Theorem 6 holds on the level of $\mathcal{R}\text{-modules.}$ To see that it holds on the level of $\mathcal{R}\text{-algebras}$ we must show that $\text{glue}_\mathfrak{T}$ respects the algebra structure.

For this fact, the same diagrammatic proof in \cite{CL19} works here. In each of the following cases: 
\begin{itemize}
	\item $x,y$ are both in $S_q^{SL_3}(\Sigma_1),$
	\item $x,y$ are both in $S_q^{SL_3}(\Sigma_2),$
	\item or $x$ is in $S_q^{SL_3}(\Sigma_1)$ while $y$ is in $S_q^{SL_3}(\Sigma_2),$
\end{itemize} it is clear that $\text{glue}_\mathfrak{T}(x)\text{glue}_\mathfrak{T}(y)=\text{glue}_\mathfrak{T}(xy).$

In the remaining case, we have that $x$ is in $\mathcal{S}_q^{SL_3}(\Sigma_2)$ and $y$ is in $\mathcal{S}_q^{SL_3}(\Sigma_1).$ We diagrammatically compute that

\begin{align*}
\text{glue}_{\mathfrak{T}}(x)\text{glue}_{\mathfrak{T}}(y) & = \begin{tikzpicture}[use Hobby shortcut, scale=2, baseline=6ex]
%bottom edge
\draw [line width=1.5] (1,0)--(-1,0);
\draw [line width=1.5, ->] (1,0)--(0,0);
%side edges
\draw [line width=1.5] (-1,0)--(-5/4,0);
\draw [line width=1.5] (1,0)--(5/4,0);
\draw [line width=1.5] (0,1)--(-1/4,5/4);
\draw [line width=1.5] (0,1)--(1/4,5/4);
%dashed lines
\draw [dashed] (-1,0)--(0,1);
\draw [dashed] (1,0)--(0,1);
%punctures
\filldraw[fill=white, draw=black] (1,0) circle [radius=.05];
\filldraw[fill=white, draw=black] (0,1) circle [radius=.05];
\filldraw[fill=white, draw=black] (-1,0) circle [radius=.05];
%crossing
{\begin{knot}[
	consider self intersections=true,
	%  draft mode=crossings,
	flip crossing=1,
	ignore endpoint intersections=false]
	\strand [line width=1.25] (-1/2,1/2)..(1/2,0);
	\strand [line width=1.25] (1/2,1/2)..(-1/2,0);
	\end{knot}}
%boxes
\node at (1/2,1/2) [fill=white, draw] (x) {$x$};
\node at (-1/2,1/2) [fill=white, draw] (y) {$y$};
%strands
\draw [line width=1.25] (y)--(-3/4,1);
\draw [line width=1.25] (x)--(3/4,1);
\end{tikzpicture} \\
&= \displaystyle\sum_{(x)(y)} \varepsilon \bigg( \begin{tikzpicture}[use Hobby shortcut, scale=2, baseline=6ex]
%sideways bigon
\draw [line width=1.5, rounded corners] (1/2,1/2)--(1/4,1)--(-1/4,1)--(-1/2,1/2)--(-1/4,0)--(1/4,0)--(1/2,1/2);
\draw [line width=1.5, rounded corners, ->] (1/2,1/2)--(1/4,1)--(0,1);
\draw [line width=1.5, rounded corners, ->] (1/2,1/2)--(1/4,0)--(0,0);
%punctures
\filldraw[fill=white, draw=black] (1/2,1/2) circle [radius=.05];
\filldraw[fill=white, draw=black] (-1/2,1/2) circle [radius=.05];
%crossing
{\begin{knot}[
	consider self intersections=true,
	%  draft mode=crossings,
	flip crossing=1,
	ignore endpoint intersections=false]
	\strand [line width=1.25] (-1/8,3/4)..(1/8,0);
	\strand [line width=1.25] (1/8,3/4)..(-1/8,0);
	\end{knot}}
%boxes
\node at (1/8+.05,3/4) [fill=white, draw] (x) {$x'$};
\node at (-1/8-.05,3/4) [fill=white, draw] (y) {$y'$};
%strands
\draw [line width=1.25] (x)--(1/4,1);
\draw [line width=1.25] (y)--(-1/4,1);
\end{tikzpicture}\bigg) \begin{tikzpicture}[scale=2, baseline=6ex]
%bottom edge
\draw [line width=1.5] (1,0)--(-1,0);
\draw [line width=1.5, ->] (1,0)--(0,0);
%side edges
\draw [line width=1.5] (-1,0)--(-5/4,0);
\draw [line width=1.5] (1,0)--(5/4,0);
\draw [line width=1.5] (0,1)--(-1/4,5/4);
\draw [line width=1.5] (0,1)--(1/4,5/4);
%dashed lines
\draw [dashed] (-1,0)--(0,1);
\draw [dashed] (1,0)--(0,1);
%punctures
\filldraw[fill=white, draw=black] (1,0) circle [radius=.05];
\filldraw[fill=white, draw=black] (0,1) circle [radius=.05];
\filldraw[fill=white, draw=black] (-1,0) circle [radius=.05];
%boxes
\node at (1/2,1/2) [fill=white, draw] (x) {$x''$};
\node at (-1/2,1/2) [fill=white, draw] (y) {$y''$};
%strands
\draw [line width=1.25] (y)--(-3/4,1);
\draw [line width=1.25] (x)--(3/4,1);
\draw [line width=1.25] (y)--(-1/2,0);
\draw [line width=1.25] (x)--(1/2,0);
\end{tikzpicture}\\
&= \displaystyle\sum_{(x)(y)} \rho (y' \otimes x') \text{glue}_{\mathfrak{T}}(y'' \otimes x'')\\
&= \text{glue}_{\mathfrak{T}} \bigg( \displaystyle\sum_{(x)(y)} \rho (y' \otimes x')(y'' \otimes x'')\bigg)\\
&= \text{glue}_{\mathfrak{T}}(xy).
\end{align*}

This shows that $\text{glue}_\mathfrak{T}$ respects the multiplication of $\mathcal{S}_q^{SL_3}(\Sigma_1) \underset{-}{\otimes} \mathcal{S}_q^{SL_3}(\Sigma_2)$ and completes our proof.

\end{proof}

By applying Theorem 6 in the special case where $\Sigma_1$ and $\Sigma_2$ are both bigons $\mathfrak{B}$ we obtain the following corollary.

\begin{corollary}
We have that $$\mathcal{S}_q^{SL_3}(\mathfrak{T}) \cong \mathcal{O}_q(SL_3) \underset{-}{\otimes} \mathcal{O}_q(SL_3).$$
\end{corollary}

\pagebreak
\printbibliography
\end{document}